\documentclass[12pt,leqno]{article}
\usepackage{amsmath}
\usepackage{amssymb}
\usepackage{theorem}
\usepackage{epsf}

\numberwithin{equation}{section}

\newtheorem{theorem}{Theorem}[section]

\newtheorem{corollary}[theorem]{Corollary}
\newtheorem{lemma}{Lemma}[section]

\theorembodyfont{\rmfamily}
\newtheorem{defn}{Definition}[section]

\newenvironment{proof}{\medskip\noindent{\bf Proof.}}{\medskip}

\newcommand{\R}{\mathbb{R}}

\newcommand{\e}{\varepsilon}
\newcommand{\dist}{\,\mathrm{dist}\,}

\newcommand{\nint}{\int \kern-1.13em {\begin{turn}{-20}$\bigm/$%
\end{turn}}\!}
\newcommand\notint{{{\,\int \kern-1.01em \raise1pt\hbox{{\begin{turn}{-30}$/$%
\end{turn}\!\!}}}}}

\renewcommand{\SS}{\mathbb S}

\newcommand\qed{\hfill\vrule height8pt width6pt depth0pt}

\newcommand{\Min}{\mathop{\rm min}\nolimits}
\newcommand{\support}{\mathop{\rm support}\nolimits}
\renewcommand{\i}{\subset}
\newcommand{\prop}{\mathrm{Prop}}

\begin{document}

\title{A generalization of Reifenberg's theorem in $\R^3$}
\author{G.\ David,\and T.\ De Pauw \and T.\ Toro\thanks{Partially 
supported by the NSF under Grant DMS-0244834}}
\date{}
\maketitle

\begin{abstract}
In 1960 Reifenberg proved the topological disc property. He showed that 
a subset of $\R^n$ which is well approximated by $m$-dimensional affine 
spaces at each point and at each (small) scale is locally a bi-H\"older image
of the unit ball in $\R^m$. In this paper we prove that a subset of $\R^3$
which is well approximated by a minimal cone at each point and at each (small)
scale is locally a bi-H\"older deformation of a minimal cone. We also prove 
an analogous result for more general cones in $\R^n$.
\end{abstract}

\section{Introduction}\label{intro}

In 1960, thanks to the development of algebraic topology, Reifenberg
was able to formulate the Plateau problem for $m$-dimensional surfaces
of varying topological type in $\R^k$ (see \cite{R1}).
He proved that given a set $\Gamma\subset \R^k$ homeomorphic to
$\SS^{m-1}$ there exists a set $\Sigma_0$ with $\partial\Sigma_0=\Gamma$ 
which minimizes the ${\cal H}^m$ Hausdorff measure among all 
competitors in the appropriate class. Furthermore he showed
that for almost every $x\in\Sigma_0$ there exists a neighborhood of $x$
which is a topological disk of dimension $m$. A 
remarkable result in \cite{R1} is the Topological Disk Theorem. 
In general terms it says that if a set is close to an $m$-plane in the 
Hausdorff distance sense at all points and at all (small enough) scales, then 
it is locally biH\"older equivalent to a ball 
of $\R^m$. Using a monotonicity formula for the density, Reifenberg proved that 
some open subset of full measure of $\Sigma_0$ satisfies this condition. 
In 1964 he proved an Epiperimetric inequality for solutions to the Plateau 
problem described above. This allowed him to show that the minimizer $\Sigma_0$
is locally real analytic (see \cite{R2}, and \cite{R3}). Although the 
Topological Disk Theorem has never again been used as a tool to
study the regularity of minimal surfaces it has played a role in
understanding their singularities as well as the singular set of energy 
minimizing
harmonic maps (see \cite{HL}). Reifenberg's proof has been adapted to 
produce biLipschitz, quasi-symmetric and biH\"older parameterizations
both for subsets of Euclidean space and general metric spaces (under the
appropriate flatness assumptions). See \cite{To}, \cite{DT}, and \cite{CC}.

In 1976 J.Taylor \cite{Ta} classified the tangent cones for Almgren almost-minimal
sets in $\R^3$. She showed that there are three types of nonempty minimal cones
of dimension 2 in $\R^3$: the planes, sets that are obtained by
taking the product of a $Y$ in a plane with a line in the orthogonal direction, 
and sets composed of six angular sectors bounded by four half lines that start 
from the same point and make (maximal) equal angles. 
See the more precise Definitions \ref{toro-defn2.2} and 
\ref{toro-defn2.3}. By lack of better names, we shall call
these minimal cones sets of type 1, 2, and 3 respectively.

In this paper we generalize Reifenberg's Topological Disk Theorem to 
the case when the set is close in the Hausdorff distance sense to a 
set of type 1, 2 or 3 at every point and every (small enough) scale. 
Namely, if $E$ is a closed set in $\R^3$ which is
sufficiently close to a two-dimensional minimal cone at all
scales and locations, then there is, at least locally, a bi-H\"older
parameterization of $E$ by a minimal cone. Reifenberg's theorem
corresponds to the case of approximation by planes.
Let us first state the main 
result and comment later.

\begin{theorem} \label{T1.1}
For each $\e \in (0,10^{-15})$, we can find $\alpha \in (0,1)$ such the
following holds. Let $E \i \R^3$ be a compact set that contains the origin,
and assume that for each $x\in E\cap B(0,2)$ and each radius
$r>0$ such that $B(x,r)\i B(0,2)$, there is a minimal cone $Z(x,r)$
that contains $x$, such that
\begin{equation}\label{eqn1}  
D_{x,r}(E,Z(x,r)) \leq \varepsilon,
\end{equation}
where we use the more convenient variant of Hausdorff distance $D_{x,r}$
defined by
\begin{eqnarray}\label{eqn2}   
D_{x,r}(E,F) & = & {\frac{1}{r}} \, \max\Big\{
\sup \big\{ \dist(z,F) \,  ; \, z\in E \cap B(x,r) \big\},
\\
&& \hskip 1.6cm \sup \big\{ \dist(z,E) \,  ;  \, z\in F \cap B(x,r)
\big\}\Big\} \nonumber
\end{eqnarray}
whenever $E$, $F$ are closed sets that meet $B(x,r)$.
Then there is a minimal cone $Z$ through the origin and
an injective mapping $f : B(0,3/2) \to B(0,2)$, 
with the following properties: 
\begin{equation}\label{1.3}  
B(0,1) \i f(B(0,3/2)) \i B(0,2), 
\end{equation}
\begin{equation}\label{eqn3}  
E\cap B(0,1) \i f(Z \cap B(0,3/2)) \i E\cap B(0,2), 
\end{equation}
\begin{equation}\label{eqn4} 
(1+\alpha)^{-1} |x-y|^{1+\alpha}
\leq |f(x)-f(y)|
\leq (1+\alpha) |x-y|^{1/{1+\alpha}}
\hbox{ for } x,y \in B(0,3/2),
\end{equation}
\begin{equation}\label{1.6}  
|f(x)-x| \leq \alpha \hbox{ for } x \in B(0,3/2).
\end{equation}
We can even take $\alpha \leq C \varepsilon$ when $\varepsilon$ is small.
\end{theorem}

Notice that by (\ref{eqn3}) and (\ref{eqn4}) 
the restriction of $f$ to $Z \cap B(0,3/2)$ provides a biH\"older 
parameterization a piece of $E$ that contains $E\cap B(0,1)$. 
This may be enough information, but in some cases it may also be good
to know that this parameterization comes from a biH\"older local
homeomorphism of $\R^3$, as in the statement, because this yields
information on the position of $E$ in space. For instance, we can use
Theorem \ref{T1.1} to construct local retractions of space near $E$ onto $E$.

\medskip\noindent {\bf Remark 1.1}
We shall see that when $Z(0,2)$ is a plane, we can take $Z$ to be a plane;
when $Z(0,2)$ is a set of type 2 centered at the origin, 
we can take $Z$ to be a set of type 2 centered at the origin;
when $Z(0,2)$ is a set of type 3 centered at the origin, 
we can take $Z$ to be a set of type 3 centered at the origin.
In addition, our proof will yield that
\begin{equation}\label{1.7}  
B(0,17/10) \i f(B(0,18/10)) \i B(0,2), 
\end{equation}
\begin{equation}\label{1.8}  
E\cap B(0,18/10) \i f(Z \cap B(0,19/10)) \i E\cap B(0,2), 
\end{equation}
and that (\ref{eqn4}) and (\ref{1.6}) holds for $x, y \in B(0,18/10)$.

In fact, we shall omit the case of the plane (too easy and well known), 
and concentrate on the two other special cases. The general case will 
essentially follow. 

\medskip\noindent {\bf Remark 1.2} 
It would be a little too optimistic to think that $f$ is 
quasisymmetric. This is already false when $Z(0,2)$ is a plane,
because $E$ could be the product of $\R$ with a snowflake in $\R^2$.

\medskip\noindent {\bf Remark 1.3} 
Theorem \ref{T1.1} can probably be generalized to a number of 
situations, where approximation by minimal sets of types 1, 2, and 3
is replaced with approximation by various types of sets with a
hierarchical simplex structure. We did not find the optimal way to 
state this, and probably this would make the proof a little heavier.
What we shall do instead is state a slightly more general result in
Theorem \ref{T2.2}, prove Theorem \ref{T1.1} essentially as it is, and
add a few comments from place to place to explain how to generalize the
proof for Theorem \ref{T2.2}. A more general statement, if ever needed, 
is left out for future investigation.

\medskip
Our proof will use the hierarchical structure of $E$.
We shall see that 
$E\cap B(0,2)$ splits naturally into three disjoint subsets
$E_1$, $E_2$, and $E_3$, where $E_j$ is the set of points where $E$
looks like a set of type $j$ in small balls centered at $x$
(more precise definitions will be given in Section \ref{decomp}). 
If $Z = Z(0,2)$ is a plane, we do not need sets of type
2 or 3 in smaller balls, and we are in the situation of Reifenberg's
theorem. Two other main cases will remain, as in Remark 1.2. 
The case when $Z$ is set of type 2, 
with its spine passing through the origin (see Definition 2.2 below), 
and the case when $Z$ is set of type 3 centered at the origin. 
In the first case, we shall see that $E_3\cap B(0,199/100)$ 
is empty and $E_2$ is locally a Reifenberg-flat one-dimensional set;
in the second case, $E_3\cap B(0,3/2)$ is just composed of one point 
near the origin, and away from this point $E_2$ is locally a 
Reifenberg-flat. See the end of Section~\ref{decomp} for details

The sets $E_1$ and $E_2$ will play a special role in the proof. Even 
though we shall in fact construct $f$ directly as an infinite composition of
homeomorphisms in space (that move points at smaller and smaller dyadic 
scales), we shall pay much more attention to the definition of $f$ on
$E_2$, and then $E$, just as if we were defining $f$ first from the spine 
of $Z$ to $E_2$, then from $Z$ to $E$, then on the rest of
$B(0,3/2)$.

\medskip
The construction yields that the restriction of $f$ to $Z$
to each of the three or six faces of $Z$ is of class $C^1$ if
the approximation at small scales is sufficiently good
(for instance if we can take $\varepsilon = C r^{\beta}$
on balls of radius $r$), see Section \ref{estimates}.

Theorem 1.1 will be used in \cite{D1} and \cite{D2} to give a 
slightly different proof of J. Taylor's regularity result for minimal 
sets from \cite{Ta}.

\medskip
The plan for the rest of this paper is as follows.
In Section \ref{prelim} we define sets of type 2 and 3 and state
a uniform version of our main theorem. We also discuss the more 
flexible version of Theorem \ref{T1.1} mentioned in Remark 1.3.
In Section \ref{geometry} we record some of the simple geometrical
facts about minimal sets of types 1, 2, 3, and in particular lower 
bounds on their relative Hausdorff distances, that will be used in the proof.
In Section \ref{decomp}  we show that $E\cap B(0,2)$ is the disjoint 
union of $E_1$, $E_2$ and $E_3$, that $E_1$ is locally 
Reifenberg-flat, and that $E_3$ is discrete. 

In Section \ref{coverings} 
we define the partitions of unity which we use in the
construction of the biH\"older parameterization. 
In Section \ref{first} we construct a parameterization of $E_2$
when $E_3$ is empty. 
In Section \ref{param} we extend this to a parameterization of $E$. 
In Section \ref{case} we explain how to modify the construction when
there is a tetrahedral point in $E_3$.
In Section \ref{extension} we extend the parameterization to the 
whole ball. Finally, in Section \ref{estimates}, we prove that 
our parameterization of $E$ is $C^1$ on each face of $Z$ when the
numbers $D_{x,r}(E,Z(x,r))$ tend to $0$ sufficiently fast, uniformly
in $x$, as $r$ tends to $0$.

\vskip 5mm

{\noindent \sc Acknowledgments:} Part of this work was 
completed in Autumn, 2005. T. Toro was visiting the Newton Institute 
of Mathematical Sciences in 
Cambridge, U.K. and the University College London in 
London. She would like to thank these institutions for their hospitality.

\section{Two other statements}\label{prelim}

We start with a uniform version of Theorem \ref{T1.1}.

\begin{defn}\label{toro-defn2.1}
A set $E\i\R^3$ is $\e$-Reifenberg flat (of dimension 2)
if for each compact set $K\i\R^3$ there exists $r_K>0$ such that for 
$x\in E\cap K$ and $r<r_K$,
\begin{equation}\label{toro-eqn2.1}
\inf_{L\ni x}D_{x,r}(E,L)\le \e,
\end{equation}
where the infimimum is taken over all planes $L$ containing $x$.
\end{defn}

Note that this definition is only meaningful for $\e$ small. 
Reifenberg's Topological Disk Theorem gives local parameterizations
of $\e$-Reifenberg flat sets. We want to extend this to the
$\e$-minimal sets defined below, but first we give a
full description of sets of type 2 and 3. Recall that sets of type 1
are just planes.

\begin{defn}\label{toro-defn2.2}  
Define $\prop\i\R^2$ by
\begin{eqnarray}\label{toro-eqn2.2} 
\prop & = & \left\{(x_1,x_2): x_1\ge 0, x_2=0\right\} \nonumber\\
&&\hskip1cm
\bigcup\left\{(x_1,x_2): x_1\le 0, x_2=-\sqrt{3}x_1\right\} \\
&&\hskip2cm
\bigcup\left\{(x_1,x_2): x_1\le 0, x_2=\sqrt{3}x_1\right\}. \nonumber
\end{eqnarray}
Define $Y_0\i\R^3$ by $Y_0=\prop \times \R$. The spine of $Y_0$
is the line $L_0 = \left\{ x_1 = x_2 = 0 \right\}$. A set of type 2 is 
a set $Y = R(Y_0)$, where $R$ is the composition of a translation and
a rotation. The spine of $Y$ is the line $R(L_0)$.
\end{defn}

\begin{defn}\label{toro-defn2.3} 
Set $A_1 = (1,0,0)$, $A_2=\left(-\frac{1}{3},
\frac{2\sqrt{2}}{3},0\right)$, $A_3=\left(-\frac{1}{3},
-\frac{\sqrt{2}}{3}, \frac{\sqrt{6}}{3}\right)$, and
$A_4=\left(-\frac{1}{3}, -\frac{\sqrt{2}}{3}, -\frac{\sqrt{6}}{3}\right)$
(the four vertices of a regular tetrahedron).
Denote by $T_0$ the cone over the union of the six edges
$[A_i,A_j]$, $i \neq j$. The spine of $T_0$ is the union of the four 
half lines $[0,A_j)$. 
A set of type 3 is a set $T = R(T_0)$, where $R$ is the composition 
of a translation and a rotation. The spine of $T$ is the image by $R$ of the
spine of $T_0$.
\end{defn}

\begin{defn}\label{toro-defn2.4} 
A closed set $E\i\R^3$ is $\e$-minimal if for each compact set $K\i\R^3$
there exists $r_K>0$ such that for $x\in E\cap K$ and $r<r_K$,
\begin{equation}\label{toro-eqn2.3}
\inf_{ x\in Z} D_{x,r}(E,Z)\le\e,
\end{equation}
where the infimimum is taken over all minimal cones (i.e., sets of type 1, 2, 
or 3) containing $x$. That is, there exists a minimal cone $Z(x,r)$ such that 
$x\in Z(x,r)$ and $D_{x,r}(E,Z(x,r))\le\e$.
\end{defn}

Observe that we do not require $Z(x,r)$ to be centered at $x$, so
minimal cones are $\e$-minimal for every $\e$. 
Definition \ref{toro-defn2.4} is only meaningful for $\e$ small. 
For the next result, we shall assume that $\e<10^{-25}$. 
Here is a uniform version of our main theorem. 

\begin{theorem} \label{T2.1}
If $E$ is $\epsilon$-minimal for some $\e<10^{-25}$
then for each compact set $K\i\R^3$
we can find $r_K>0$ such that for $x\in E\cap K$ and $r<r_K$
there exists a minimal cone $Z$ through the origin and a
map $f:B(0,2r)\to \R^3$, with $f(0)=x$, satisfying
\begin{equation}\label{2.4}
(1-C\epsilon) \big(r^{-1}|y-z| \big)^{1+C\epsilon}\le r^{-1} |f(y)-f(z)|
\le (1+C\epsilon) \big(r^{-1}|y-z|\big)^{1-C\epsilon},
\end{equation}
\begin{equation}\label{eqn101}
B(x,r)\subset f(B(0,\frac{3r}{2}))\subset B(x,2r),
\end{equation}
\begin{equation}\label{eqn102}
E\cap B(x,r)\subset f(Z\cap B(0,\frac{3r}{2}))\subset E\cap B(x,2r).
\end{equation}
\end{theorem}

Thus $f$ is a local homeomorphism of the ambient space which
sends a minimal cone into the $\epsilon$-minimal set.
Theorem \ref{T2.1} is an immediate consequence of Theorem \ref{T1.1}.
Also, if we can take $\varepsilon$ to tend to $0$ in (\ref{toro-eqn2.3})
(uniformly in $x\in E \cap K$) when $r$ tends to $0$, then 
(\ref{2.4}) holds with constants that tend to $0$ when $r$ tends to $0$.

\medskip
Let us now try to generalize Theorem \ref{T1.1}. We want to extend the 
notion of sets of type 1, 2, or 3. Let us fix integers
$d \geq 2$ (the dimension of our sets) and $n \geq d+1$
(the ambient dimension).

Sets of type G1 are just $d$-planes in $\R^n$.

A generalized propeller in $\R^2$ will be any union $P$ of three co-planar
half 
lines with the same origin, and that meet with angles larger than $\pi/2$ 
at this point. Incidentally, $\pi/2$ was chosen a little at random here, 
and we shall not try to see whether it  can be replaced by smaller angles. 

A set of type G2 is a set $Y = R(P \times \R^{d-1})$, where $P$ is a 
generalized propeller in $\R^2$ and $R$ is an isometry of $\R^n$. When
$n>d+1$, we abused notation slightly, we should have written
$Y = R(P \times \R^{d-1} \times \{0\})$ to account for the
last $n-d-1$ coordinates. The spine of $R(P \times \R^{d-1} \times \{0\})$
is $R( \{0\} \times \R^{d-1} \times \{0\})$.

A two-dimensional set of type G3 in $\R^m$ is a set $T$ that we can 
construct as follows. We pick distinct points $A_j$, $1 \leq j \leq k$,
in the unit sphere of $\R^m$.
Then we choose a graph, as follows. We pick $G \i \{ 1, \cdots, k\}^2$ 
which is symmetric (i.e., $(i,j) \in G$ when $(j,i)\in G$), 
does not meet the diagonal, and is such that for each 
$i\in \{ 1, \cdots, k\}$, there are exactly three $j \neq i$ such that
$(i,j)\in G$. We call $\Gamma$ the union of the segments $[A_i,A_j]$,
$(i,j)\in G$, and then set 
$T = \{ t \gamma \, ; \, t \geq 0 \hbox{ and } \gamma \in \Gamma \}$
($T$ is the cone over $\Gamma$). Equivalently, for each $(i,j) \in G$,
we set $F_{i,j} = \{ tz \, ; \, t \geq 0 \hbox{ and } z\in [A_i,A_j]\}$,
and $T$ is the union of the $3k/2$ faces $F_{i,j}$. In particular,
$k$ is even. We add a few regularity constraints. First,
the interior of a face $F_{i,j}$ never meets another face, so the faces
only meet by sets of three along half lines $L_j = [0,A_j)$.
For each $j$, exactly three faces $F_{i,j}$ meet $L_j$, and we require
that they make angles larger than $\pi/2$ along $L_j$.
Thus the union of three half planes that contains these three faces
is a two-dimensional set of type G2. We also require that
\begin{equation}\label{2.7}  
|A_i - A_j| \geq \tau_0 \ \hbox{ for } i\neq j
\end{equation}
and
\begin{equation}\label{2.8}   
\dist(F_{i,j}\cap \partial B(0,1),F_{i',j'}\cap \partial B(0,1))
\geq \tau_0 \ \hbox{ when } F_{i,j}\cap F_{i',j'} = \emptyset,
\end{equation}
for some small positive constant $\tau_0$ that we pick in advance.


Another way to state the constraints is to consider
$\Gamma' =T \cap \partial B(0,1)$: we want $\Gamma'$ to be composed of
not too small arcs of circles that only meet at their ends, by sets of three,
and with angles greater than $\pi/2$. In addition, we only allow two
arcs to be very close when they share at least an end.

Let us also allow the case when more than one arc may connect
a given pair of points $A_j$ (this would imply changing our notation
a tiny bit, and taking a graph $G$ that is not necessarily contained in
$\{ 1, \cdots, k\}^2$), but demand that $k \geq 4$ (to avoid the case of a
set of type G2).


A set of type G3 (of dimension $d$ in $R^n$) is a set
$Z=R(T\times \R^{d-2})$, where $R$ is an isometry of $\R^n$
and $T$ is a two-dimensional set of type G3 in $\R^{n-d+2}$.
The spine of $Z$ is the union $L$ of the half lines $L_i$ when
$Z=T$ is two-dimensional, and $R(L\times \R^{d-2})$ when
$Z=R(T\times \R^{d-2})$.

Denote by $TGi$ the collection of sets of of type Gi
(of dimension $d$ in $\R^n$), and by $TG$ the union of the three
classes $TGi$. We claim that in Theorem \ref{T1.1}, we can replace the class 
of minimal sets with the class $TG$, as follows.

\begin{theorem} \label{T2.2}   
Let $n$ and $d$ be as above. Let $E \i \R^n$ be a compact set that contains 
the origin, and suppose that for each $x\in E\cap B(0,2)$ and $r>0$ such that 
$B(x,r)\i B(0,2)$, we can find $Z(x,r) \in TG$ that contains $x$, such that
$D_{x,r}(E,Z(x,r)) \leq \varepsilon$. If $\varepsilon >0$ is small 
enough, depending only on $n$, $d$, and $\tau_0$, there is a set $Z\in TG$
such that (\ref{1.3})--(\ref{1.6}) hold. Here again $\alpha$ depends 
only on $\varepsilon$, $n$, $d$, and $\tau_0$, and we can take 
$\alpha = C(n,d,\tau_0)\varepsilon$ for $\varepsilon$ small.
\end{theorem}

The reader may be surprised that we do not require any coherence
between the various sets $Z(x,r)$, but this coherence will follow
automatically from the fact that $D_{x,r}(E,Z(x,r)) \leq \varepsilon$
for many $x$ and $r$.
For instance, if $d=2$ and $Z(0,2)$ is of type G3, with a center at the origin,
our proof will show that for $r$ small, there is a point $x_r\in E$,
close to $0$, such that $Z(x_r,r)$ is of type G3, and the number
$k$ of half-lines in the spine of $Z(x_r,r)$ is the same as for
$Z(0,2)$. The point is that angles between the faces may vary little by 
little when $x$ and $r$ vary, but things like the type of $Z(x,r)$ and the 
number $k$ cannot jump randomly.

\medskip
There is  an analogue of Theorem \ref{T2.1} in the context of 
Theorem \ref{T2.2}, which the reader may easily state.

\medskip
The proof of Theorem \ref{T2.2} will almost be the same as for 
Theorem \ref{T1.1}; what we shall do is proceed with the proof of 
Theorem \ref{T1.1}, and indicate from time to time which 
modifications are required for Theorem \ref{T2.2}.

\section{Simple geometrical facts}\label{geometry}

We shall record in this section a few geometrical properties
of the minimal cones that will be used in the construction.
The statements will come with various numbers, like $1/3$ or
$1/4$ in the next lemma; these numbers are hopefully correct, 
but their precise values do not matter, in the sense that it 
is easy to adapt the value of later constants to make the proof
work anyway. In particular, this is what would happen with the
proof of Theorem \ref{T2.2}. 
This means that the reader may skip the proofs if she believes 
that the results of this section are true with different constants.
The geometrical facts below will be used extensively in 
Section \ref{decomp}, to establish the hierarchical structure of $E$.

Our first lemma says that the Hausdorff distances between 
sets of type 1, 2, and 3, are not too small. 

\begin{lemma}\label{toro-lem1}   
Let $Z$ be a minimal cone of type 3 centered at $x$. Then
\begin{equation}\label{toro-eqn2.7} 
D_{x,r}(Z,Y) \ge 1/3 \hbox{ for } r > 0
\end{equation}
whenever $Y$ is a set of type 1 or 2. Similarly,
\begin{equation}\label{toro-eqn2.8}  
D_{x,r}(Y,P) > 1/3 \hbox{ for } r > 0
\end{equation}
when $Y$ is a set of type 2 centered at $x$ and $P$ is a plane,
and
\begin{equation}\label{toro-eqn2.8+}  
D_{x,r}(Y,P) > 1/4 \hbox{ for } r > 0
\end{equation}
when $Y$ is a set of type 2 or 3 whose spine contains $x$ and $P$ 
is a plane.
\end{lemma}

\begin{proof}
We start with the proof of (\ref{toro-eqn2.7}).
By translation, dilation, and rotation invariance, we may 
assume that $x=0$, $r=1$, and $Z$ is the set $T_0$ of 
Definition \ref{toro-defn2.3}.
It will be good to know that if $P$ is any plane through the
origin and $\pi$ denotes the orthogonal projection onto $P$, then
\begin{equation}\label{2.9}  
B(0,1/3) \cap P \i \pi(T_0 \cap B(0,1)).
\end{equation}

Let $Q$ denote the solid tetrahedron with vertices $A_j$, $1 \leq j \leq 4$,
(i.e., the convex hull of the $A_j$), where the $A_j$ are as in 
Definition \ref{toro-defn2.3}. It is clear that $B(0,1/3)$ 
lies on the right of the left face of $Q$ (where the first coordinate 
is $-1/3$). By symmetry,
$B(0,1/3)$ lies in $Q$, because it lies on the right side of each face.
Also observe that $T_0 \cap Q$ separates the (open) faces of $Q$ 
from each other inside $Q$. For instance, the component 
in $Q \setminus T_0$ of the lower 
face is a pyramidal domain bounded by three faces of $T_0$.

Let $z\in B(0,1/3) \cap P$ be given, and let $\ell$ denote the
line through $z$ perpendicular to $P$. Then $\ell$ meets $Q$
because $B(0,1/3) \i Q$. First suppose that it does not touch 
any edge of $Q$. Then it enters $Q$ through one
(open) face and leaves it through another one, so it meets 
$T_0\cap Q$. If $\ell$ touches a edge of $Q$, it meets 
$T_0\cap Q$ trivially (because the edges are contained in $T_0$).
Thus $\ell$ meets $T_0 \cap B(0,1)$ (except perhaps if 
$\ell$ contains one of the $A_j$), and hence $z\in \pi(T_0 \cap B(0,1))$.
The remaining case when $\ell$ contains some $A_j$ is easily obtained
by density (and anyway we don't need it); (\ref{2.9}) follows.

Return to (\ref{toro-eqn2.7}). If $Y$ is of type 2, choose $P$
orthogonal to the spine of $Y$; then $\pi(Y)$ is a propeller
(a set like $\prop$ in (\ref{toro-eqn2.2})), and we can find 
$z\in \partial B(0,1/3) \cap P$ such that $\dist(z,\pi(Y)) \geq 1/3$.
Set $z_n = (1-2^{-n}) \, z$ for $n>1$.
By (\ref{2.9}) we can find $y_n \in T_0 \cap B(0,1)$ such that $\pi(y_n) = z_n$.
Then $\dist(y_n,Y) = \dist(z_n,\pi(Y)) \geq 1/3 - 2^{-n}$,
and (\ref{toro-eqn2.7}) follows from the definition (\ref{eqn2})
The case when $Y$ is a plane is even easier; we choose $P$
perpendicular to $Y$, so that $\pi(Y)$ is a line. So 
(\ref{toro-eqn2.7}) holds in all cases.

\medskip
Let us now prove (\ref{toro-eqn2.8}). We may assume that $x=0$ and $Y$
is the set $Y_0 = \prop \times \R$ in (\ref{toro-defn2.2}).
Let $P$ be a plane, suppose that $D_{x,r}(Y,P)\le 1/3$, and let us
find a contradiction. 

First we want to show that $P$ is almost horizontal.
Call $P_0$ the horizontal plane $\R^2 \times \{ 0 \}$, and denote by
$b_1$, $b_2$, and $b_3$ the three points of 
$\partial B(0,\frac{9}{10}) \cap Y \cap P_0 =
\partial B(0,\frac{9}{10}) \cap \big[ \prop\times\{ 0 \}\big]$.
Each $b_j$ lies in $Y \cap B(0,1)$, 
so we can find $b'_j$ in $P \cap \overline B(b_j,1/3)$ because 
$D_{x,r}(Y,P)\le 1/3$. Denote by $\pi$ the orthogonal projection on 
$P_0$, and set $p_j = \pi(b'_j)$. The $p_j$ lie
in the disks $D_j = P_0 \cap \overline B(b_j,1/3)$.  

Call $d$ the smallest possible distance from a $p_j$ to the line through 
the other two $p_l$. Then $d$ corresponds to the case when, for instance, 
$p_j$ lies at the extreme right of the left-most disk $D_j$ and the two other 
$p_l$ lie at the extreme left of the $D_l$. 
We get that $d = \frac{9}{10}\frac{3}{2} - 2 \cdot \frac{1}{3} = 
\frac{41}{60}$.

Recall that $P$ goes through the $b'_j$.
The fact that $d>0$ (or that the $p_j$ are not aligned) already implies
that $P$ is a graph over $P_0$. We want to show that it is a 
$1$-Lipschitz graph, or equivalently that
\begin{equation}\label{toro-eqn2.9} 
|v_3| \leq |\pi(v)|
\end{equation}
when $v$ is a unit vector in the (direction of) the plane $P$,
and $v_3$ denotes its last coordinate.
Since the triangle with vertices $b'_j$ is nondegenerate, we
can write $v = \alpha(b'_j-c_j)$, where $\alpha\in \R \setminus \{ 0 \}$,
$1 \leq j \leq 3$, and $c_j$ lies in the opposite side 
$[b'_k,b'_l]$. Since $c_j$ is a convex combination of 
$b'_k$ and $b'_l$, the size of its last coordinate is at most $\frac{1}{3}$;
then the size of the last coordinate of $b'_j-c_j$ is at most $\frac{2}{3}$.
On the other hand, $|\pi(b'_j-c_j)| \geq d = \frac{41}{60}$ by 
definition of $d$. Now (\ref{toro-eqn2.9}) follows because 
$\frac{2}{3} < \frac{41}{60}$.

Finally consider the points $b_4=(0,0,-{\frac{9}{10}})$ and 
$b_5=(0,0,{\frac{9}{10}})$. Since they lie in $Y \cap B(0,1)$,
we can find $b'_4$ and $b'_5$ in $P$, with $|b'_j-b_j|\leq 1/3$.
Set $v=b'_5-b'_4$. Then $v_3 \geq \frac{9}{5}-\frac{2}{3}=\frac{17}{15}$, 
while $|\pi(v)| \leq \frac{2}{3}$,
so $v$ does not satisfy (\ref{toro-eqn2.9}). This contradiction 
completes our proof of (\ref{toro-eqn2.8}).

For the proof of (\ref{toro-eqn2.8+}), we may assume that $x=0$ and 
that the spine of $Y$ contains the positive third axis 
$\{ x_1 = x_2 =0 \hbox{ and } x_3 \geq 0 \}$.
Then $Y$ coincides in $\{ x_3 \geq 0 \}$ with a set of type 2 
whose spine is the third axis, and we may even assume that (after a
rotation around the third axis) this set is $P_0$. 
We shall only need to know that $Y$ contains
$Y_0 \cap \{ x_3 \geq 0 \}$. 

This time we assume that $D_{x,r}(Y,P)\le 1/4$. We can follow the 
proof of (\ref{toro-eqn2.9}), except that now our disks $D_j$ are smaller
and $d = {\frac{9}{10}}\frac{3}{2} - 2 \cdot \frac{1}{4} = \frac{17}{20}$.
We get that $|v_3| \leq \frac{2}{3} \cdot \frac{20}{17} \, |\pi(v)|
= \frac{40}{54}  \, |\pi(v)|$ instead of (\ref{toro-eqn2.9}).

Now we take $b_4 = 0$ and $b_5=(0,0,{\frac{99}{100}})$.
As before, we can find $b'_4$ and $b'_5$ in $P$, with $|b'_j-b_j|\leq 1/4$ 
(because $D_{x,r}(Y,P)\le 1/4$). Set $v=b'_5-b'_4$. 
Then $v_3 \geq {\frac{99}{100}}-\frac{1}{2}={\frac{49}{100}}$,
and $|\pi(v)| \leq \frac{1}{2}$, so $v_3 \geq \frac{98}{100} |\pi(v)|
> \frac{40}{54}  \, |\pi(v)|$, a contradiction; 
(\ref{toro-eqn2.8+}) and Lemma \ref{toro-lem1} follow.
\qed\end{proof}

Notice that even for sets of type $TGi$, Lemma \ref{toro-lem1}
is very easy to prove by compactness, when we allow very small
constants instead of $1/3$ and $1/4$.

The next lemmas will make it easier to apply Lemma \ref{toro-lem1}.

\begin{lemma}\label{L2.2}  
Let $T$ be a minimal cone of type 3, and let $B$ be a ball such that 
$\frac{5}{4}B$ (the ball with the same center, and $\frac{5}{4}$ times
the radius) does not contain the center of $T$. Then there 
is a minimal cone $Y$ of type 1 or 2 such that $Y\cap B = T\cap B$.
\end{lemma}

\begin{proof}
Again this would be very easy, even for sets of type $TGi$,
if we allowed a very large constant instead of $\frac{5}{4}$.
We shall even prove the lemma with $\frac{5}{4}$ replaced
with the better constant $\alpha^{-1}$, where $\alpha = (2/3)^{1/2}$
(observe that $\frac{5}{4} > \alpha^{-1}=(3/2)^{1/2}$).
By translation and rotation invariance, we can assume that
$T$ is the set $T_0$ from Definition \ref{toro-defn2.3}. 

Let $P$ denote the plane through $0$, $A_2$, and $A_3$, and set
$x_t = (-t,0,0)$ for $t > 0$; we claim that  
$\dist(x_t,P) = \alpha t$.
Indeed, an equation of $P$ is  $2 \sqrt 2 x + y + \sqrt 3 z=0$ 
(just check that $0$, $A_2$, and $A_3$ satisfy the equation).
Then $\dist(x_t,P) = (2\sqrt 2 t)/\sqrt{12} = \alpha t$, as announced.

Let $x \neq 0$ be given, and set $B_x = B(x,\alpha|x|)$. We just 
checked that when $x = x_t$, $B_x$ does not meet the three faces on the 
left of $T_0$. For a general $x$, $B_x$ cannot meet the three faces at 
the same time, because the optimal position for this to happen would 
precisely be when $x$ lies on the negative first axis.
By symmetry, $B_x$ never meets the three faces of $T$ that bound a given 
connected component of $\R^3 \setminus T$. So the worse that it can do is meet 
the three faces of $T$ that share a single edge $[0,A_j)$ (if $B_x$ 
meets two opposite faces, it also meets the three faces that bound 
some connected component). 

Suppose first that $B_x$ meets three such faces, and let $Y$ denote the set 
of type 2 that contains these faces. We claim that
\begin{equation}\label{3.6b} 
T \cap B_x = Y \cap B_x.
\end{equation}
The direct inclusion is clear, because $B_x$ only meets the three faces
of $T$ that are contained in $Y$. For the converse inclusion, we just
need to check that $F \cap B_x \i F' \cap B_x$ when $F$ is a face
of $Y$ and $F'$ is the face of $T$ that is contained in $F$.
Let $y\in F \cap B_x$ be given. By assumption, $B_x$ meets $F'$, so
we can pick $z \in F' \cap B_x$. Since $B_x$ is open, we can even pick
$z$ in the interior of $F'$. Notice that the segment $[y,z]$ lies in 
$F \cap B_x$ by convexity. In addition, $(y,z)$ does not meet the spine of $Y$
(because $z$ lies in the interior of $F$). Now $B_x$ does not meet
the other edges of $T$ (those that are not contained in the spine of 
$Y$), because if it meets an edge $e$, it also meets the three faces 
of $T$ that touch $e$. So $(y,z)$ does not meet the boundary of $F'$,
and $y\in F'$. Our claim follows.

Next suppose that $B_x$ only meets two faces of $T$. Still denote by
$Y$ the set of type 2 that contains them. As before, 
$T \cap B_x \i Y \cap B_x$ trivially. To prove the converse inclusion,
consider a face $F$ of $Y$, and first assume that the face $F'$ of
$T$ that it contains meets $B_x$. Pick $z\in F' \cap B_x$. For each
$y \in F \cap B_x$, the segment $[x,y]$ is contained in $B_x$,
hence it does not meet the boundary of $F'$ (because if $B_x$
meets any edge of $T$, it meets at least three faces), so it is 
contained in $F'$, and $y\in F'$. We just proved that 
$F \cap B_x \i F' \cap B_x $. 

Call $F_1$, $F_2$, and $F_3$ the three face of $Y$, and denote
by $F'_j$ the face of $T$ which is contained in $F_j$. Exactly
two of the $F'_j$ meet $B_x$, and the third one does not; let us 
assume that this last one is $F'_3$. We just need to check that 
$F_3$ does not meet $B_x$ either. Suppose it does. Let $y_j$ be 
a point of $F_j \cap B_x$, $j\in \{1,2,3\}$. Call $y'_j$ the 
orthogonal projection of $y_j$ onto the plane $P$ through
the center of $B_x$ and perpendicular to the spine $L$ of $Y$;
then $y'_j \in F_j \cap B_x$ as well. 
Call $p$ the point of $L\cap P$; by convexity of $B_x$ and because
$p$ is a convex combination of the three $y'_j$, $B_x$ contains $p$. 
We reached the desired contradiction, because $p$ lies on an edge of 
$P$ and $B_x$ does not meet edges in the present case. 
So we proved (\ref{3.6b}) in this second case.

Finally assume that $B_x$ only meets one face $F'$ of $T$. Then
it does not meet any edge of $T$, so if $Y$ denotes the plane that
contains $F'$, $T\cap  B_x$ contains $Y\cap B_x$ (as before, $B_x$
does not meet the boundary of $F'$ in $Y$). In this case also we have
(\ref{3.6b}). 

We just proved that for $x \neq 0$, we can find a set $Y$ of type 1 or 
2 such that (\ref{3.6b}) holds. Now let $B$ be as in the statement, 
and call $x$ its center and $r$ its radius. 
We know that $B(x,\frac{5r}{4})$ does not contain the origin,
so $r \leq 4|x|/5 \leq \alpha |x|$, $B \i B_x$, and the analogue of
(\ref{3.6b}) holds for $B$.
This completes our proof of Lemma \ref{L2.2}.
\qed\end{proof}

\begin{lemma}\label{L2.3}  
Let $Z$ be a minimal cone of type 3 centered at $z$, and let 
$T$ be a minimal cone. If $D_{x,r}(T,Z) < 1/3$, then $T$ is of type 3
and its center lies in $B(z,5r/3)$.
\end{lemma}

\begin{proof}
Let $Z$ and $T$ be as in the statement. If $T$ is of type 1 or 2,
we can apply Lemma~\ref{toro-lem1} directly to get a contradiction.
So $T$ is of type 3, and let $t$ denote its center. By Lemma~\ref{L2.2}, 
we can find a set $Y$ of type 1 or 2 that coincides with $T$ in 
$B = B(z,\frac{4}{5}|z-t|)$. If $|z-t| \geq 5r/3$, 
$B$ contains $B(z,\frac{4}{3}r)$, $Y$ coincides with $T$ in 
$B(z,\frac{4}{3}r)$, and it is easy to see that $D_{x,r}(Y,Z) < 1/3$
because $D_{x,r}(T,Z) < 1/3$. But this is impossible, by Lemma~\ref{toro-lem1},
so $|z-t| < 5r/3$. Lemma \ref {L2.3} follows.
\qed\end{proof}

\begin{lemma}\label{L3.4} 
Let $Z$ be a set of type 2 or 3. 
Suppose that the ball $B(x,r)$ meets at least two faces of $Z$, 
or that $Z$ does not coincide with any plane in 
$B(x,r)$. Then the distance from $x$ to the spine
of $Z$ is at most $6r/5$.
\end{lemma}

\begin{proof}
For sets of type G2 or G3 (and with $6/5$ replaced with a large 
constant), this is a simple consequence of (\ref{2.8}). We now return
to the standard case.
Let $L$ denote the spine of $Z$. We can assume that $B=B(x,r)$ 
does not meet $L$. If $Z$ does not 
coincide with any plane in $B$, it meets at least one face $F$; call $P$
the plane that contains $F$. Since $B$  does not meet $L$, it
does not meet the boundary of $F$ in $P$, so $P\cap B \i Z$.
Then $B$ meets some other face of $F'$ of $Z$. Call $P'$ the plane 
that contains $F'$, and set $D=P\cap P'$ and $d=\dist(x,D)$.
Then $r \geq d\cos(30^\circ) = d \sqrt 3/2$ (the worse case is
when $x$ lies in the middle of the $120^\circ$ sector bounded
by $P$ and $P'$, and $B$ is tangent to both planes). 

If $Z$ is of type 2, we get the result directly, because $L=D$ and 
$\sqrt 3/2 \geq 5/6$. So we can assume that $Z$ is of type 3. 
If $B(x,d)$ meets $L$, then $\dist(x,L) \leq d \leq 6r/5$ as needed.
Otherwise, $P\cap B(x,d) \i F$ , because $B(x,d)$ meets $F$ and
does not meet the boundary of $F$ in $P$. But $P\cap B(x,d)$
goes all the way to the point of $D$ that minimizes the distance
to $x$, so this point lies in $L$ and $\dist(x,L) \leq d$.
Lemma \ref{L3.4} follows.
\qed\end{proof}

\section{The decomposition of $E$}\label{decomp}

From now on, $E$ is a compact set in $\R^3$ that satisfies the 
assumptions of Theorem \ref{T1.1}. In particular, (\ref{eqn1}) 
holds for some $\varepsilon<10^{-25}$.
The following definitions will make it easier to define the $E_j$'s
and work with them. For $x\in E$ and $r>0$, set
\begin{equation}\label{eqn5}
a(x,r)  =  \inf \big\{ D_{x,r}(E,P) \, ; \, P
\hbox{ is a plane through } x \big\}, 
\end{equation}
\begin{equation}\label{eqn6}
\qquad b(x,r)  =  \inf \big\{ D_{x,r}(E,Y) \, ; \, 
\begin{array}[t]{l}Y \hbox{ is a set of type 2 whose}
	\hbox{ spine contains } x \big\},\end{array}
\end{equation}
and
\begin{equation}\label{eqn7}
c(x,r) = \inf \big\{ D_{x,r}(E,T) \, ; \, T
\hbox{ is a set of type 3 centered at } x \big\}.
\end{equation}
It is not always true that either $a(x,r)$, $b(x,r)$, or $c(x,r)$ is small
(when $x\in E$ and $B(x,r) \i B(0,2)$), 
because, for instance $Z(x,r)$ may be a minimal cone of type 3 
centered anywhere in $B(x,r)$. The pairs where 
either $a(x,r)$, $b(x,r)$, or $c(x,r)$ is small
are interesting to study.

\begin{defn}\label{defn8}
Let $x \in E\cap B(0,2)$ be given. We say that:
\begin{itemize}
\item{} $x$ is of
\underbar{type 1} when $a(x,r) \leq 2\varepsilon$ for all $r$ small,
\item{} $x$ is of \underbar{type 2} when $b(x,r) \leq 1500 \varepsilon$
for all $r$ small, 
\item{} $x$ is of \underbar{type 3} when
$c(x,r) \leq 150\varepsilon$ for all $r$ small.
\end{itemize}
For $j=1,2,3$, we denote by $E_j$ the set of points
$x \in E\cap B(0,2)$ that are of type $j$.
\end{defn}

We should say something now about the choice of constants.
The constants 2, 1500, and 150 are just chosen to make later
lemmas easier to use. They depend on the values of the various 
constants in Section \ref{geometry}, but we could reorganize things 
with different constants. Other strange constants will appear
in the proof; the reader should not pay too much attention
to their precise values, hopefully this will make them easier to track.
For Theorem \ref{T2.2}, all these constants need to be replaced
with  other ones; let us indicate in advance how they need to be chosen,
to make it easier for the reader to check things out.

Our first constants are $a_1 = 2$, $a_2 = 1500$, and $a_3 = 150$ 
in Definition \ref{defn8}; $a_1$ comes from the proof of 
Lemma \ref{lem9}, and we can keep $a_1 = 2$; the constraints
on (the new values of)  $a_2$ and $a_3$ will be explained soon.

The constant $a_4=10^{-3}$ in Lemma \ref{lem9} just needs to be small,
depending on $a_1$. In Lemma~\ref{lem10}, $a_5 = 5$ and $a_6 = 57$ 
just come from the geometry, $a_7 = 10$ is just $2a_5$, we can keep 
$\frac{1}{2}$ as it is, and the value of $a_8=7$ follows from the geometry.
Then $a_9=10^{-3}$ just needs to be small, depending on  
$a_5$ and $a_6$. Also, $a_{10}=11$ just needs to be a little larger than 
$a_7$. The constant $a_{11}=15$ in Corollary \ref{cor11} just needs to
be larger than $2a_8$. 

The constants $a_{12}=10^{-3}$, $a_{13}=10^{-3}$, and $a_{14}=132$
in Lemma~\ref{lem12} are just geometric constants that come from 
Lemma~\ref{L2.3}. Our first constraint on $a_3$ comes from 
Lemma~\ref{lem12} and is that $a_3 \geq a_{14}$.

No constraint comes from Lemma \ref{lem13}; the constant in 
(\ref{3.12}) is just the same as $a_{11}$ in Corollary \ref{cor11}.

In Lemma \ref{lem14}, we just need $a_{15}=25$ needs to be large
enough (it needs to be larger than 10 in our proof, because we want
the pair $(z_1,r/10)$ to satisfy the assumptions in the lemma,
a little above (\ref{eqn16})); other constraints will come later. 
The value of $a_{16}=600$ follows
by geometric computations, and so does $a_{17}=150$ in 
(\ref{eqn15}). We need to pick $a_3 \geq a_{17}$ in 
Definition \ref{defn8}.

The first constraint on $a_2$ in Definition \ref{defn8}
is that $a_2$ be larger than the constant in (\ref{3.20a}).

The constant $a_{18}=17$ in Lemma~\ref{lem17} is just a geometric 
constant. In (\ref{3.20}) and the few lines above, $10^{-3}$ needs
to be replaced with $a_9$ from Lemma~\ref{lem10}. A few lines later,
we only get that $b(x,r) \leq a_{19}\varepsilon$, instead of 
$500\varepsilon$; this gives a new constraint on $a_2$, namely that
$a_2 \geq a_{19}$ (so that we can deduce that $x\in E_2$).
The rest of the proof of (\ref{eqn19}) goes smoothly.

In Lemma~\ref{lem21}, $a_{20}=24$ can be replaced with 
$\frac{2}{3}(2a_{18}+2)$ and $a_{21}=600$ with the geometric constant 
$a_{16}\,$; in the proof, we need to replace $198/100$ with a geometric 
constant $C$, just below (\ref{3.23}); the rest of the proof is the 
same, but gives a second constraint on $a_{15}$ in Lemma~\ref{lem14},
i.e., it needs to be larger than $a_{20}$.

For Lemma~\ref{lem24}, we need to replace $15$ with 
$a_{11}$, as in Corollary \ref{cor11}, and $200$ with
some geometric constant $a_{23}$ that depends on 
Lemmas \ref{L2.2} and \ref{L3.4}. We also get the last constraint 
that $a_{15}$ be large enough, depending on $a_{23}$.

This completes our list of constraints on constants for the results 
of this section; the reader may check that they are compatible.

\medskip
The main point of this section is to establish a few properties of the sets
$E_j$. We shall prove that $E\cap B(0,2) = E_1 \cup E_2 \cup E_3$ 
(a disjoint union), that $E_3\cap B(0,199/100)$ has at most one point, 
and that $E_2$ is locally Reifenberg-flat of dimension 1 
(away from $E_3$ and $\partial B(0,2)$). 

The results of this section are also valid with minor modifications for 
$\varepsilon$-minimal sets (as in the previous section), 
but the reader should not bother too much; most of the statement
in this section can also be easily deduced from Theorem \ref{T2.1}.

Let us first check that
\begin{equation}\label{3.4a}   
\hbox{$E_1$, $E_2$, and $E_3$ are disjoint.}
\end{equation}

Suppose for instance that $x \in E_2 \cap E_3$.
For $r$ small enough, we can find a set $Y$ of type 2 and a set
$Z$ of type 3 centered at $x$, such that 
$D_{x,2r}(E,Y) \leq 1500 \varepsilon$
and $D_{x,2r}(E,Z) \leq 150 \varepsilon$. If $y\in Y \cap B(x,r)$,
we can find $e\in E$ such that $|e-y| \leq 3000 \varepsilon r$.
Then $e\in E \cap B(x,2r)$, so we can find $z\in Z$ such that
$|z-e| \leq 300 \varepsilon r$. 
Altogether, $|z-y| \leq 3300 \varepsilon r$. Similarly, if
$z\in Z\cap B(x,r)$, we can find $e\in E$ such that 
$|e-y| \leq 300 \varepsilon r$, 
then $e\in E \cap B(x,2r)$, we can find $y\in Y$ such that
$|y-e| \leq 2000 \varepsilon r$, and $|z-y| \leq 3300 \varepsilon r$.
So $D_{x,r}(Y,Z) \leq 3300 \varepsilon < 1/3$. 
Lemma~\ref{toro-lem1} says that this is impossible. The two other 
cases are just as easy.

Our first lemma shows that $E_1$ is a large open set.

\begin{lemma}\label{lem9} 
Let $x\in E$ and $r>0$ be such that $B(x,r) \i B(0,2)$
and $a(x,r) \leq 10^{-3}$. Then
$a(y,t) \leq 2\varepsilon$ for $y\in B(x,2r/3)$ and $0 < t < r/4$.
In particular, $E\cap B(x,2r/3) \i E_1$.
\end{lemma}

\begin{proof} Let $x$ and $r$ be as in the statement. 
Let $P_0$ be a plane through $x$ such that $D_{x,r}(E, P_0)\leq 10^{-3}$.
By (\ref{eqn1}), there is a minimal cone $Z=Z(x,r)$ through $x$,
such that $D_{x,r}(E,Z) \leq \varepsilon$.
It could be that $Z$ is of type 2 or 3, but even so let us check that 
\begin{equation}\label{3.5} 
\hbox{there is a plane $P$ through $x$ 
such that } Z \cap B(x,97r/100) = P \cap B(x,97r/100).
\end{equation}

First suppose that $Z$ is a set of type 3 centered at 
some $y \in B(x,99r/100)$. Set $\rho = r/200$, and let us
check that $D_{y,\rho}(Z,P_0) < 1/3$. If 
$z\in Z\cap B(y,\rho)$, we can find $e\in E$ such that
$|z-e| \leq \varepsilon r$. Obviously, $e\in B(x,r)$, so
we can find $p\in P_0$ such that $|p-e| \leq 10^{-3} r$;
then $|p-z| \leq \varepsilon r + 10^{-3} r < \rho/3$. Similarly,
if $p\in P_0 \cap B(y,\rho)$ we can find $e\in E$ such that
$|p-e| \leq 10^{-3} r$. Then $e\in B(x,r)$, so we can find $z\in Z$
such that $|z-e| \leq \varepsilon r$; altogether 
$|p-z| \leq \varepsilon r + 10^{-3} r < \rho/3$, and 
$D_{y,\rho}(Z,P_0) < 1/3$. This is impossible, 
by Lemma~\ref{toro-lem1}.

Next suppose that $Z$ is a set of type 3, whose
center $z\not\in B(x,99r/100)$, but whose spine meets 
$B(x,98r/100)$ at some point $y$.
As before, $D_{y,\rho}(Z,P_0) < 1/3$ for $\rho = r/200$.
Since $|z-y| \geq r/100$, Lemma \ref{L2.2} says that
$Z$ coincides with a set $Y$ of type 2 in $B(y,4r/500)= B(y,8\rho/5)$.
Then $D_{y,\rho}(Y,P_0) < 1/3$ too. But the spine of $Y$ goes through $y$, 
so Lemma~\ref{toro-lem1} says that this is impossible.
The same argument excludes the case when $Z$ is a set of type 2 whose spine 
meets $B(x,98r/100)$.

We are left with the case when $Z$ is of type 1, or else its spine does 
not meet $B(x,98r/100)$. Let $F$ denote the face of $Z$ that contains $x$, 
and $P$ the plane that contains $F$. The boundary of $F$ is contained in 
the spine of $Z$, so it does not meet $B(x,98r/100)$; hence
\begin{equation} \label{3.6a}  
F \cap B(x,98r/100) = P\cap B(x,98r/100).
\end{equation}
Every point of $P \cap B(x,98r/100)$ lies in $Z$ by (\ref{3.6a}),
so it is $\varepsilon r$-close to $E$, and then
$(\varepsilon + 10^{-3})r$-close to $P_0$.
This forces every point of $P_0 \cap B(x,r)$ to be 
$2 \cdot 10^{-3}r$-close to $P$. This stays true (with a constant
larger than $2$) when we deal with sets of type $Gi$ and approximations 
with $d$-planes in $\R^n$.

If (\ref{3.5}) fails, (\ref{3.6a}) says that there is another face 
$F_1$ of $Z$ that meets $B(x,97r/100)$ at some point $y$. 
We know that $y$ is $\varepsilon r$-close to $E$, 
hence $(\varepsilon + 10^{-3}) \, r$-close to $P_0$.
Hence, $\dist(y,P) \leq (\varepsilon + 3 \cdot 10^{-3})r$.
[In all these estimates, we use the fact that $y\in B(x,99r/100)$, so 
the successive points that we implicitly use never lie out of 
$B(x,r)$.]

On the other hand, $y$ lies in some other face of $Z$, 
and the angles between faces of $Z$ are not too small, 
so the distance from $y$ to the spine of $Z$ is less than 
$10^{-2}r$. For sets of type $Gi$, we deduce this from (\ref{2.8}). 
This is impossible, because the spine of $Z$ does not meet 
$B(x,98r/100)$. So (\ref{3.5}) holds. 

\medskip
We now return to the lemma, and let $y\in B(x,2r/3)$ be as in the statement. 
Let us first estimate $a(y,t)$ for $t\in (r/10,4r/10)$.
First observe that $B(y,t) \subset B(x,967r/1000)$.
By (\ref{eqn1}) and (\ref{3.5}), we can find $q_y \in P$ such that 
$|y-q_y|\leq\e r$. Set $P'=P+(y-q_y)$. We can use $P'$ 
to compute $a(y,t)$, because it goes through $y$, so 
$a(y,t) \leq D_{y,t}(E,P')$. If $e \in E \cap B(y,t)$, 
$\dist(e,P') \leq \dist(e,P)+|y-q_y| \leq \dist(e,P)+\varepsilon r
\leq 2\varepsilon r$, by (\ref{eqn1}) and (\ref{3.5}).
If $p'\in P'\cap B(y,t)$, $p=p' -(y-q_y)$ lies in $Z \cap B(x,97r/100)$
by (\ref{3.5}), and $\dist(p',E) \leq \dist(p,E)+\varepsilon r 
\leq 2\varepsilon r$ by (\ref{eqn1}). So $a(y,t) \leq D_{y,t}(E,P')
\leq 2\varepsilon r/t \leq 20\varepsilon $.

The pair $(y,t)$ satisfies the hypothesis of the lemma, namely, 
$a(y,t) \leq 10^{-3}$, so we can iterate the previous argument 
(this time, keeping $y$ at the center). This yields that 
$a(y,s) \leq 20\varepsilon$ for $r/100 \leq s \leq 4r/10$, and 
(after many iterations of the argument) for every $s < 4r/10$.

Finally, for each $s < 4r/10$ there is a set $Z(y,s)$ as in 
(\ref{eqn1}), and since $a(y,s) \leq 20\varepsilon$ the proof 
of (\ref{3.5}) shows that $Z(y,s)$ coincides with a plane $P$ on 
$B(y,97s/100)$.
We can use $P$ in the definition of $a(y,96s/100)$, and we get that
$a(y,96s/100) \leq 100\varepsilon/96 \leq 2\varepsilon$.
Since this holds for $s < 4r/10$, Lemma \ref{lem9} follows.
\qed
\end{proof}

Next we focus on $E_3$ and $c(x,r)$. We start 
by showing how being close to a tetrahedron at a small scale determines 
the behavior of the set at larger scales.

\begin{lemma}\label{lem10}  
Let $x\in E$ and $r > 0$ be such that $B(x,11r) \i B(0,2)$ and
$c(x,r)<10^{-3}$. Let $Z_{11} = Z(x,11r)$ be a minimal set through
$x$ such that $D_{x,11r}(E,Z_{11})\le\e$ (as in (\ref{eqn1})). 
Then $Z_{11}$ is of type 3, and its center $x_0$ is such that 
$|x-x_0| \leq 5 c(x,r)  r + 57 \varepsilon r$.
In addition, 
\begin{equation}\label{eqn3,5} 
c(x,10r) \leq \frac{1}{2}\, c(x,r) + 7 \varepsilon.
\end{equation}
\end{lemma}

\begin{proof}
The reader should not pay too much attention to the value of the
various constants. Nevertheless the fact that $\frac{1}{2} < 1$ will be 
useful. 
By assumption, there is a minimal set $Z$ of type 3, centered at $x$, 
and such that $D_{x,r}(E,Z)\le c(x,r)\le 10^{-3}$. Set
$\rho = 3 c(x,r) r + 34 \varepsilon r$, and let us check that 
\begin{equation}\label{4.8} 
D_{x,\rho}(Z,Z_{11}) < 1/3.
\end{equation}

If $z\in Z\cap B(x,\rho)$, we can find $e\in E$ such that
$|e-z| \leq c(x,r) r$. Obviously $e\in B(x,11r)$, so we can find
$z'\in Z_{11}$ such that $|z'-e| \leq 11\varepsilon r$; altogether 
$\dist(z,Z_{11}) \leq (c(x,r) + 11 \varepsilon)r < \rho/3$. Similarly, 
if $z'\in Z_{11} \cap B(x,\rho)$, we can find $e\in E$ such that
$|e-z'| \leq 11\varepsilon r$, and then $z\in Z$ 
such that $|z-e| \leq c(x,r)r$, so 
$\dist(z',Z) \leq (c(x,r) + 11 \varepsilon) r < \rho/3$. This
proves (\ref{4.8}).

By Lemma \ref{L2.3}, $Z_{11}$ is of type 3 and its center $x_0$ lies in
$B(x,5\rho/3)$. That is, $|x-x_0| \leq 5 c(x,r)  r + 170 \varepsilon r/3
\leq 5 c(x,r)  r + 57 \varepsilon r$.

We cannot use $Z_{11}$ directly to estimate $c(x,10r)$, because
it is not centered at $x$, but we may use $Z'=Z_{11} + (x-x_0)$. 
If $y\in E \cap B(x,10r)$, we can find $z\in Z_{11}$ such that
$|z-y| \leq 11 \varepsilon r$, and then $z'=z+(x-x_0) \in Z'$
is such that $|z'-y| \leq 5 c(x,r)  r + 68 \varepsilon r$.
Similarly, if $z'\in Z' \cap B(x,10r)$, then
$z = z'-(x-x_0) \in Z \cap B(x,11r)$, so we can find 
$y\in E$ such that $|z-y| \leq 11 \varepsilon r$, and then 
$|z'-y| \leq 5 c(x,r)  r + 68 \varepsilon r$. Altogether,
$D_{x,10r}(E,Z') \leq (10 r)^{-1}(5 c(x,r)  r + 68 \varepsilon r)
\leq \frac{1}{2} c(x,r) + 7 \varepsilon$. 
Lemma \ref{lem10} follows.
\qed
\end{proof}

\begin{corollary}\label{cor11}  
If $x\in E_3$, then $c(x,r) \leq 15\varepsilon$ for every $r>0$
such that $B(x,11r/10) \subset B(0,2)$. 
\end{corollary}

\begin{proof} Indeed for $k$ large, $c(x,10^{-k}r) \leq 150 \varepsilon$
by definition of $E_3$. Lemma \ref{lem10} says that 
$c(x,10^{-k}r) \leq \frac{1}{2} \, c(x,10^{-k+1}r) + 7 \varepsilon
\leq 157 \varepsilon$ for $k$ large. 
Multiple iterations of Lemma \ref{lem10} eventually lead to the result.
In fact define a sequence 
$\lambda_n$ by $\lambda_0=150\varepsilon$ 
and $\lambda_{n+1}=\frac{1}{2} \lambda_n + 7 \varepsilon$, note that
it converges to $14\varepsilon < 15\varepsilon$.
\qed
\end{proof}

The influence of tetrahedral points also transmits to smaller scales.

\begin{lemma}\label{lem12} 
Let $x\in E$ and $r>0$ be such that $B(x,r) \subset B(0,2)$
and $c(x,r) \leq 10^{-4}$. Then there is a point
$\xi \in E \cap B(x,10^{-3}r)$ such that
$c(\xi,\rho) \leq 132 \varepsilon$ 
for $0 \leq \rho \leq r/3$.
In particular, $\xi \in E_3$.
\end{lemma}

\begin{proof} 
Pick $Z=Z(x,r)$, so that $D_{x,r}(E,Z) \leq \varepsilon$.
Also choose a minimal cone $Y$ of type 3 centered at $x$
and such that $D_{x,r}(E,Y) \leq c(x,r)\leq 10^{-3}$.
Set $\rho_0 = 3 c(x,r) r + 4 \varepsilon r$; then 
$D_{x,\rho_0}(Y,Z) < 1/3$, by the proof of (\ref{4.8}).
As before, Lemma \ref{L2.3} says that $Z$ is of type 3, 
with a center $z_0$ in $B(x,5\rho_0/3)$. That is, 
$|z_0-x| \leq 5 c(x,r) r + 20 \varepsilon r/3$.
By definition of $Z$, we can find $x_1 \in E$ so that 
$|x_1-z_0| \leq \varepsilon r$, and hence
\begin{equation}\label{3.9}   
|x_1-x| \leq 5 c(x,r) r + 8 \varepsilon r.
\end{equation}
The cone $Z' = Z + (x_1-z_0)$ is centered at $x_1$, so we can use 
it to compute $c(x_1,r/2)$. The same computation as in Lemma \ref{lem10}
yields
\begin{equation}\label{eqn3,7} 
c(x_1,r/2) \leq D_{x_1,r/2}(E,Z') 
\leq (2/r) \big(|x_1-z_0|+ \varepsilon r \big)
\leq 4 \varepsilon.
\end{equation}

Let us iterate the construction with the ball $B(x_1,r/2)$.
We get a second point $x_2 \in E$, with
$|x_{2}-x_1| \leq [5c(x_1,r/2)+8\varepsilon] r/2 
\leq 28\varepsilon r/2$, and such that $c(x_2,r/4) \leq 4\varepsilon$.
By iteration we can find points $x_k$, $k\geq 1$,
such that $c(x_k,2^{-k}r) \leq 4 \varepsilon$ and
$|x_{k+1}-x_k| \leq 28 \varepsilon 2^{-k} r$ for $k\geq 1$.

Call $\xi$ the limit of the $x_k$; thus 
\begin{eqnarray}\label{eqn3,8} 
|x- \xi| &\leq &|x - x_1| + \sum_{k \geq 1} |x_{k+1}-x_k| 
\nonumber\\
&\leq &[5c(x,r)+8\varepsilon]r + \sum_{k \geq 1} 28 \varepsilon 
2^{-k} r 
\leq [5 c(x,r) + 36 \varepsilon]r. 
\end{eqnarray}

We still need to check that $c(\xi,\rho) \leq 132 \varepsilon$ 
for $0 \leq \rho \leq r/3$. Pick $k$ such that
$2^{-k-1}r\le 11 \rho /10 \leq 2^{-k}r$; thus $k \geq 1$ and hence
$c(x_k,2^{-k}r) \leq 4 \varepsilon$.
Also, $|\xi-x_k| \leq \sum_{l \geq k} 28 \varepsilon 2^{-l} r 
\leq 56 \varepsilon 2^{-k} r$.
Call $W$ a set of type 3 centered at $x_k$ and which is
$4 \varepsilon 2^{-k}r$-close to $E$ in $B(x_k,2^{-k}r)$;
then the translation of $W$ by $\xi-x_k$ is
$60 \varepsilon 2^{-k}r$-close to $E$ in $B(\xi,\rho)$,
which leads to $c(\xi,\rho) \leq 60 \varepsilon 2^{-k}r/\rho
\leq (60 \varepsilon)(22/10) < 132 \varepsilon$, as needed.
The fact that $\xi \in E_3$ is just the definition of $E_3$
\qed
\end{proof}

\begin{lemma}\label{lem13}  
There is at most one point in $E_3 \cap B(0,199/100)$.
\end{lemma}

\begin{proof} 
Suppose that $E_3\cap B(0,199/100)$ contains two different points 
$x$ and $y$. Notice that 
\begin{equation}\label{3.12}   
c(x,r) \leq 15\varepsilon
\ \hbox{ for } 0 < r \leq 2 \cdot 10^{-3},
\end{equation}
by Corollary \ref{cor11} and because $B(x,11r/10) \subset B(0,2)$.

Set $\rho = |x-y|$, and first assume that $\rho<10^{-3}$. 
We take $r=2\rho$ in (\ref{3.12}) and get that 
$c(x,2\rho) \leq 15\varepsilon$.
Thus there is a set $X$ of type 3, centered at $x$, such that
$D_{x,2\rho}(X,E) \leq 15\varepsilon$. Similarly, there is a 
set $Y$ of type 3, centered at $y$, such that
$D_{y,2\rho}(Y,E) \leq 15\varepsilon$. Notice that 
$B(y,2\rho)$ contains $B(x,\rho)$, so 
$D_{x,\rho}(X,Y) \leq 60\varepsilon$ by the proof of (\ref{3.4a}).
By Lemma \ref{L2.3}, the center of $Y$ lies in
$B(x,5\rho/6)$, a contradiction.

We are left with the case when $\rho = |x-y| \geq 10^{-3}$.
Since $c(x,10^{-3}) \leq 15 \varepsilon$ by (\ref{3.12}),
we can pick $X$ of type 3, centered at $x$, and such that
$D_{x,10^{-3}}(X,E) \leq 15 \varepsilon$.
Also set $Z = Z(0,2)$; thus $D_{0,2}(Z,E) \leq \varepsilon$ 
by (\ref{eqn1}), and the proof of (\ref{3.4a}) yields
$D_{x,10^{-4}}(X,Z) \leq 3 \cdot 10^{4}\varepsilon$. Then
Lemma~\ref{L2.3} says that $Z$ is of type 3, with a center
in $B(x,2 \cdot 10^{-4})$. By the same argument, the center
of $Z$ lies in $B(y,2 \cdot 10^{-4})$. But these balls are
disjoint because $|x-y| \geq 10^{-3}$, our last case is impossible, 
and Lemma \ref{lem13} follows. \qed
\end{proof}

\medskip
The situation with respect to $E_3$ is reasonably clear now:
it contains at most one point $x_0$ in $B(0,199/100)$, and 
Corollary \ref{cor11}
says that $c(x_0,r) \leq 15\varepsilon$ for all $r$ such that
$B(x_0,11r/10) \subset B(0,2)$. Furthermore, Lemma \ref{lem12} 
essentially forbids $E$ to look like a set of type 3 away from $x_0$, 
because this would create a new point of $E_3$. 
Next we focus on $b(x,r)$ and $E_2$.

\begin{lemma}\label{lem14}  
Let $x\in E$ and $r>0$ be such that $B(x,r) \subset B(0,2)$,
and assume that we can find a minimal cone $Y$ of type 2,
whose spine $L$ contains $x$, and such that
$D_{x,r}(E,Y) \leq 25\varepsilon$. Then $B(x,2r/3)$
does not meet $E_3$ and $D_{x,2r/3}(L,E_2) \leq 600 \varepsilon$.
\end{lemma}

\begin{proof} Let $x$ and $r>0$ be as in the lemma.
We shall need to check that
\begin{equation}\label{3.13a}   
z \in E_1  \ \hbox{ for $z\in E\cap B(x,2r/3)$ such that }
\dist(z,L) \geq 400 \varepsilon r.
\end{equation}  

Let $z \in E\cap B(x,2r/3)$ be given, with 
$\dist(z,L) \geq 400 \varepsilon r$. Set $\rho=\Min(\dist(z,L)/2,r/5)$; 
in particular, $16\rho/10 \leq 16 r/50 < r/3$, so 
$B(z,16\rho/10) \i B(x,r) \i B(0,2)$. Notice that 
$\dist(z,Y) \leq 25\varepsilon r$ because 
$D_{x,r}(E,Y) \leq 25\varepsilon$. Call $P$ the plane that contains 
the face of $Y$ that gets close to $z$;
then $Y \cap B(z,16\rho/10) = P \cap B(z,16\rho/10)$
because $16\rho/10 \leq 8 \dist(z,L)/10$.
Also recall that $25\varepsilon r \leq 25\rho/200$, 
either trivially because $\rho = r/5$ or else because 
$\rho = \dist(z,L)/2 \geq 200\varepsilon r$.

Now set $\widetilde Z=Z(z,16\rho/10)$. Let us check that 
\begin{equation}\label{3.13}   
\dist(w,P) \leq 16\varepsilon\rho /10 + 25\rho/200
\leq 13\rho/100
\ \hbox{ for } w\in \widetilde Z \cap B(z,14\rho/10). 
\end{equation}
Indeed we can find $e\in E$ such that $|e-w| \leq 16\varepsilon\rho /10$
(by (\ref{eqn1})); then $e\in B(z,15\rho/10)\i B(x,r)$, so we can find
$y\in Y$ such that $|y-e| \leq 25\varepsilon r \leq 25\rho/200$.
Thus $|y-w| \leq 16\varepsilon\rho /10 + 25\rho/200 < 2\rho/10$, so
$y\in Y \cap B(z,16\rho/10) = P\cap B(z,16\rho/10)$, 
and (\ref{3.13}) holds.

We deduce from (\ref{3.13}) and elementary  geometry that
$\widetilde Z$ coincides with a plane in $B(z,11\rho/10)$.

In the case of sets of type $Gi$, we need to replace $400$ in 
(\ref{3.13a}) with a much larger constant, so that
$\dist(w,P) \leq \alpha\rho$ for $w\in \widetilde Z \cap B(z,14\rho/10)$,
with a very small constant $\alpha$. If $\widetilde Z$ is of type
$G2$ or $G3$, this prevents its spine from meeting 
$B(z,13\rho/10)$ (by direct inspection: sets of type $G2$
and $G3$ do not stay close to planes near their spine); 
then Lemma \ref{L3.4} says that $Z$ coincides with a plane in 
$B(z,C^{-1}\rho)$. This is enough to continue the argument
(but get worse constants).

Next, $\widetilde Z$ goes through $z$, so we can use this plane to estimate
$a(z,\rho)$. We get that
$a(z,\rho) \leq \frac{16}{10} D_{z,16\rho/10}(E,\widetilde Z) \leq
\frac{16}{10}\,\varepsilon$. Since we know that $B(z,\rho) \i B(x,r) 
\subset B(0,2)$, we can
apply Lemma \ref{lem9} and get that $a(z,t) \leq 2\varepsilon$
for $t<\rho/4$; hence $z\in E_1$, and (\ref{3.13a}) holds.

Now we want to show that
\begin{equation}\label{eqn15}  
\dist(\ell,E_2) \leq 150 \varepsilon r
\ \hbox{ for } \ell \in L\cap B(x,2r/3).
\end{equation}

Since $\ell \in L \i Y$, we can find $z_0\in E$ such that
$|z_0-\ell| \leq 25 \varepsilon r$. Set $Z=Z(z_0,r/2)$. Let us first 
use the fact that $Z$ and $Y$ are both quite close to $E$ to show that 
\begin{equation}\label {3.15}  
D_{z_0,r/4}(Y,Z) \leq 102\varepsilon.
\end{equation}

First let $y \in Y \cap B(z_0,r/4)$ be given. Observe that
$|y-x| \leq |y-z_0|+|z_0-\ell|+|\ell-x| \leq r/4 + 25 \varepsilon r +
2r/3 = 11r/12 + 25 \varepsilon r$, so $y\in B(x,r)$ and (by definition 
of $Y$) we can find $y'\in E$ such that $|y'-y| \leq 25 \varepsilon r$. 
Then $y'\in B(z_0,r/2)$, so we can find $y'' \in Z$
such that $|y''-y'| \leq \varepsilon r/2$. Thus
$\dist(y,Z) \leq 51\varepsilon r/2$. Similarly, if
$y\in Z \cap B(z_0,r/4)$ we can find $y'\in E$ such that
$|y'-y| \leq \varepsilon r/2$, then $y'\in B(x,r)$, so there exists
$y'' \in Y$ such that $|y''-y'| \leq 25 \varepsilon r$.
Thus $\dist(y,Y) \leq 51\varepsilon r/2$; (\ref{3.15}) follows.

Recall that $|z_0-\ell| \leq 25 \varepsilon r$, with 
$\ell \in L \cap B(x,2r/3)$ (and $L$ is the spine of $Y$).
By (\ref {3.15}) and the second half of Lemma~\ref{toro-lem1}, 
$Z$ is of type 2 or 3. We claim that
\begin{equation}\label {3.16}  
\hbox{$Z$ coincides with a set of type 2 in $B(z_0,r/6)$.}
\end{equation}

This is clear if $Z$ is of type 2. Otherwise, denote its
center by $z_1$. If $B(z_1,r/500) \i B(z_0,r/4)$, (\ref {3.15})
says that $D_{z_1,10^{-3}r}(Y,Z) \leq 10^5 \varepsilon < 1/3$.
This is impossible, by Lemma~\ref{L2.2}.
So $z_1$ lies out of $B(z_0,t)$, with $t=r/4-r/500$, and  
Lemma \ref{L2.2} says $Z$ coincides with a set of type 2 in 
$B(z_0,4t/5)$;  (\ref {3.16}) follows.

Once we know (\ref {3.16}), it is easy to deduce from (\ref {3.15}) 
that every point of $L\cap B(z_0,r/10)$ lies within $100 \varepsilon r$
of the spine $L_Z$ of $Z$. This uses the fact that for sets of type 2,
the Hausdorff distance between spines is controlled by the Hausdorff 
distance between the sets. For sets of type $G2$, we can check 
this directly, or deduce it from Lemmas \ref{toro-lem1} and \ref{L3.4},
but with worse constants.
We apply this to our point $\ell \in L$, and we find $\ell_1 \in L_Z$ 
such that $|\ell_1 - \ell| \leq 100 \varepsilon r$. 
By definition of $Z$, we can then find $z_1 \in E$ such that 
$|z_1-\ell_1| \leq \varepsilon r/2$ (see (\ref{eqn1})).
Thus 
\begin{equation}\label{3.17a} 
|z_1-z_0| \leq |z_1-\ell_1|+|\ell_1 - \ell|+|\ell-z_0|
\leq 126 \varepsilon r.
\end{equation}

Call $Y'$ the set of type 2 provided by (\ref {3.16}),
and then set $Y_1 = Y'+(z_1-\ell_1)$. Notice that the spine
$L_1$ of $Y_1$ goes through $z_1$ (because (\ref {3.16})
says that $\ell_1$ lies in the spine of $Y'$), so we can use it to
estimate $b(z_1,r/10)$. By (\ref {3.16}), $Y'$ is 
$\varepsilon r/2$-close to $E$ in $B(z_1,r/7)$, hence
$Y_1$ is $\varepsilon r$-close to $E$ in $B(z_1,r/8)$, so
$b(z_1,r/10) \leq 10\varepsilon$.

The pair $(z_1,r/10)$ satisfies the assumptions of Lemma \ref{lem14}; 
thus we can apply the argument above, with $x$, $\ell$, and $z_0$ 
all replaced with $z_1$, and $r$ replaced with $r/10$. 
We get a new point $z_2 \in E$,
with $|z_2-z_1| \leq 101 \varepsilon r/10$ and 
$b(z_2,r/100) \leq 10\varepsilon$ (we can drop the
$25\varepsilon r$ that come from $|\ell-z_0|$ in (\ref{3.17a})).

Then we iterate and find a sequence $\{ z_k \}$ in $E$, with
\begin{equation}\label{eqn16} 
b(z_k,10^{-k}r) \leq 10\varepsilon 
\ \hbox{ and } \ |z_{k+1}-z_k| \leq 101 \cdot 10^{-k}\varepsilon r.
\end{equation}

Set $w=\lim_{k\to \infty} z_k$. Then
$|w-z_k| \leq 120 \cdot 10^{-k}\varepsilon r$ for $k \geq 1$.
Let us check that 
\begin{equation}\label{3.20a} 
b(w,t) \leq 1500\varepsilon  \ \hbox{ for $t$ small.}
\end{equation}

Choose $k$ such that 
$\frac{9}{10}10^{-k-1}r \leq t \leq \frac{9}{10}10^{-k}r$.
By (\ref{eqn16}), there is a set $T$ of type 2, with a spine
through $z_k$, which is $10^{-k+1}\varepsilon r$-close to $E$
in $B(z_k,10^{-k}r)$. We can use $T+ (w-z_k)$ to estimate
$b(w,t)$. Since $|w-z_k| \leq 120 \cdot 10^{-k}\varepsilon r$,
we get that $b(w,t) \leq t^{-1}(130 \cdot 10^{-k}\varepsilon r)
\leq 1500\varepsilon$, as needed for (\ref{3.20a}).

By (\ref{3.20a}), $w\in E_2$.
Observe that $|w-\ell| \leq |w-z_1|+|z_1-\ell| \leq
120 \varepsilon r/10 + |z_1-\ell_1|+|\ell_1 - \ell| 
\leq 113 \varepsilon r$, by (\ref{eqn16}) and the estimates above 
(\ref{3.17a}). So (\ref{eqn15}) holds.

We are ready to estimate $D_{x,2r/3}(L,E_2)$.
We just checked that $\dist(\ell,E_2) \leq 150 \varepsilon r$
for $\ell \in L\cap B(x,2r/3)$. Conversely, if 
$z \in E_2\cap B(x,2r/3)$, (\ref{3.13a}) forbids 
$\dist(z,L)$ to be more than $400 \varepsilon r$ (because
$E_1$ is disjoint from $E_2$, see (\ref{3.4a})).
Hence $D_{x,2r/3}(L,E_2) \leq 600 \varepsilon$, as needed.

We still need to check that $E_3$ does not meet $B(x,2r/3)$.
If $z\in E_3 \cap B(x,2r/3)$, Corollary~\ref{cor11}
says that $c(z,\rho)  \leq 15\varepsilon$ for $\rho > 0$ such that
$B(z,11\rho/10) \i B(0,2)$. We pick $\rho = r/4$, and find a set
$Z'$ of type 3, centered at $z$, and such that 
$D_{z,r/4}(Z',E) \leq 15\varepsilon$. But $D_{x,r}(E,Y) \leq 
25\varepsilon$, hence $D_{z,r/5}(Y,Z) \leq 140\varepsilon$
by the usual comparison argument. This is impossible, by 
Lemma~\ref{toro-lem1}, so $E_3$ does not meet $B(x,2r/3)$.
This completes our proof of Lemma \ref{lem14}.\qed
\end{proof}

Next we study the local Reifenberg-flatness of $E_2$.

\begin{lemma}\label{lem17}  
If $x\in E \cap B(0,2) \setminus E_1$, then for every $r>0$ such that
$B(x,r) \subset B(0,2)$, $Z(x,r)$ is a set type 2 or 3, whose spine
passes at distance at most $17 \varepsilon r$ from $x$.
\end{lemma}

\begin{proof} 
For $x\in E$ and $r>0$  such that $B(x,r)  \i B(0,2)$, we have a
minimal cone $Z(x,r)$. Set $\delta(x,r) = + \infty$  when 
$Z(x,r)$ is of type 1. Otherwise, denote by $L(x,r)$ the spine
of $Z(x,r)$, and set $\delta(x,r) = r^{-1} \dist(x,L(x,r))$.
First observe that 
\begin{equation}\label{3.17} 
\delta(x,r) < 2/3 \ \hbox{ when } a(x,r/2) \geq 10^{-4}.
\end{equation}
Indeed, if $\delta(x,r) \geq 2/3$, $Z(x,r)$ is of  type 1 or $L(x,r)$ 
does not meet $B(x,2r/3)$, hence Lemma \ref{L3.4} says that 
$Z(x,r)$ coincides with a plane $P$ in $B(x,10r/18)$.
Then we can use $P$ to check that $a(x,r/2) < 10^{-4}$.
So (\ref{3.17}) holds. If we were dealing with sets of type $Gi$,
we would prove (instead of (\ref{3.17}))
that $\delta(x,r) < 2/3$ when $a(x,r/C) \geq a_4/10$,
where $a_4$ comes from Lemma \ref{lem9}.

Next we check that
\begin{equation}\label{eqn18} 
\delta(x,3r) \leq \frac{1}{3} \, \delta(x,r) + 11 \varepsilon.
\end{equation}
when $a(x,r/2) > 10^{-4}$ and $B(x,3r) \subset B(0,2)$.
For sets of type $Gi$, we would write
$\delta(x,2Cr) \leq \frac{1}{2C} \, \delta(x,r) + C' \varepsilon$
when $a(x,r/C) > a_{4}/10$ and $B(x,2Cr) \subset B(0,2)$, 
with $C$ as above.

By (\ref{3.17}), we can find $z_0 \in L(x,r)$, with
$|x-z_0| \leq r\delta(x,r) \leq 2r/3$. 
Set $\rho =  16 \varepsilon r$; trivially, $B(z_0,\rho) \i B(x,9r/10)$.
For $z\in Z(x,r) \cap B(z_0,\rho)$, there is a point $z_1 \in E$ such that
$|z_1-z| \leq \varepsilon r$, and then a point $z_2 \in Z(x,3r)$
such that $|z_2-z_1| \leq 3\varepsilon r$, so
$\dist(z,Z(x,3r)) \leq 4\varepsilon r$. Similarly, 
$\dist(z,Z(x,r)) \leq 4\varepsilon r$ for 
$z\in Z(x,3r) \cap B(z_0,\rho)$. 
So $D_{z_0,\rho}(Z(x,r),Z(x,3r)) \leq 1/4$, and 
Lemma \ref{toro-lem1} says that $Z(x,3r)$ cannot coincide with a
plane in $B(z_0,\rho)$. Hence, $Z(x,3r)$ is of type 2 or 3, and
Lemma \ref{L3.4} says that its spine meets $B(z_0,2\rho)$.
Then $\dist(x,L(x,3r)) \leq |x-z_0| + 2 \rho 
\leq r\delta(x,r)+32 \varepsilon r$, and 
$\delta(x,3r) \leq \frac{1}{3} \delta(x,r) + 11 \varepsilon$, as 
needed for (\ref{eqn18}).

Let us also check that $a(x,3r/2) > 10^{-4}$ in the proof of (\ref{eqn18}).
Indeed, otherwise there is a plane $P$ such that 
$D_{x,3r/2}(P,E) \leq 10^{-4}$. At the same time, 
$D_{x,r}(Z(x,r),E) \leq \varepsilon$ and $z_0 \in B(x,2r/3)$,
so $D_{z_0,r/10}(Z(x,r),P) \leq 16 \cdot 10^{-4}$ by the usual
argument, and this is impossible because the spine of $Z(x,r)$ 
goes through $z_0$ and by Lemma \ref{toro-lem1}.
So we can apply (\ref{eqn18}) again to the pair $(x,3r)$, at least
if $B(x,9r) \i B(0,2)$.

Now let $x$ and $r$ be as in the statement of Lemma \ref{lem17}.
Since $x\notin E_1$, Lemma \ref{lem9} says that $a(x,\rho) \geq 10^{-3}$ 
for $\rho$ small. Multiple applications of (\ref{eqn18}) (starting from
$\delta(x,3^{-k}r)$, with $k$ very large) show that 
$\delta(x,r) < 17 \varepsilon$ 
(because $\frac{3}{2} \cdot 11\varepsilon < 17 \varepsilon$). 
Lemma \ref{lem17} follows.
\qed\end{proof}

\medskip
We are now ready to prove that
\begin{equation}\label{eqn19}   
E \cap B(0,2) = E_1 \cup E_2 \cup E_3 \ \hbox{ (a disjoint union)}.
\end{equation}

We already know from (\ref{3.4a}) that the union is disjoint. 
Let $x\in E\cap B(0,2)$ be given, assume that $x\notin E_1$,
and let us check that $x$ lies in $E_2$ or $E_3$.

Suppose first that there are arbitrarily small $r$ such that
$c(x,r) \leq 10^{-3}$; then multiple applications of Lemma \ref{lem10}
give arbitrarily small $\rho$ such that $c(x,10^k \rho) \leq 15\varepsilon$
for every $k \geq 0$ such that $B(x,\frac{11}{10} \cdot 10^k \rho) \subset B(0,2)$.
[The proof is the same as for Corollary \ref{cor11}.]
Then we also have that $c(x,r) \leq 150\varepsilon$ for every
$r>0$ such that $B(x,11\rho) \subset B(0,2)$. [Use any $\rho < r$ and 
pick $k$ such that $10^{k-1}\rho\leq r < 10^k \rho$.]
In this case $x\in E_3$ and we are happy.
So may assume that 
\begin{equation}\label{3.20}   
c(x,r) > 10^{-3} \ \hbox{ for $r$ small enough.}
\end{equation}

Let $r$ be small, and set $Z = Z(x,10r)$. Since $x\notin E_1$, 
Lemma \ref{lem17} says that $Z$ is of type 2 or 3, with a spine $L$ 
that goes through $\overline B(x,170 \varepsilon r)$.
Pick $\ell \in L \cap \overline B(x,170 \varepsilon r)$.
If $Z$ coincides with a set $Y$ of type 2 in $B(x,2r)$, 
then $Y' = Y + (x-\ell)$ is a set of type 2 whose spine
goes through $x$, we can use it to evaluate $b(x,r)$, and we get that
$b(x,r) \leq 10\varepsilon + 170 \varepsilon < 500\varepsilon$.
If this is the case for every small $r$, then $x \in E_2$ and we are 
happy.

So it is enough to show that for $r$ small, $Z$ coincides with a set
of type 2 in $B(x,2r)$. Recall that $Z$ is of type 2 or 3; if
it is of type 2, we are happy. If it is of type 3, with a center
$z_0 \notin B(x,5r/2)$, Lemma \ref{L2.2} says that $Z$ coincides
with a set of type 1 or 2 in $B(x,2r)$, and we are happy too (type 1
is excluded because $L$ goes through $\overline B(x,170 \varepsilon r)$).
So we can assume that $z_0 \in B(x,5r/2)$. [With sets of type $Gi$,
we only get that $z_0 \in B(x,Cr)$, but this is enough.]

Pick $\xi \in E$ such that $|\xi-z_0| \leq 10\varepsilon r$; we can 
use $Z+(\xi-z_0)$ to evaluate $c(\xi,r)$, and we get that
$c(\xi,r) \leq 20 \varepsilon$. If $r$ is small enough, 
$B(x,10^5r)\i B(0,2)$, we can apply Lemma \ref{lem10} four times, 
and we get that $c(\xi,10^4r) \leq 20 \varepsilon$. This means that 
there is a set $T$ of type 3, centered at $\xi$, and such that 
$D_{\xi,10^4r}(T,E) \leq 20 \varepsilon$. We can use $T + (x-\xi)$
to evaluate $c(x,5000r)$, and we get that
$c(x,5000r) \leq (5000 r)^{-1}(3r+2 \cdot 10^5 \varepsilon r) < 10^{-3}$.
This is impossible, by (\ref{3.20}), so we proved that $x\in E_2$, and
our proof of (\ref{eqn19}) is complete.

\medskip
Before we return to the local Reifenberg-flatness of $E_2$ away from $E_3$
(if $E_3 \neq \emptyset$), observe that
\begin{equation}\label{eqn20}  
E_2 \cup E_3 \hbox{ is closed in } B(0,2).
\end{equation}
Indeed, let $x\in B(0,2)$ be the limit of a sequence in $E_2 \cup E_3$.
If $x\in E_1$, $a(x,r) \leq 10\varepsilon$ for $r$ small, and
Lemma \ref{lem9} says that $E \cap B(x,2r/3) \subset E_1$. 
This is impossible, by
definition of $x$ and because the sets $E_i$ are disjoint. Recall 
that $E$ is closed, so $x \in E \setminus E_1$, and (\ref{eqn19}) 
says that $x\in E_2 \cup E_3$; (\ref{eqn20}) follows.

\begin{lemma}\label{lem21}   
If $x\in E_2$ and $0 < r \leq \frac{1}{2}
\Min\big[\dist(x,E_3),\dist(x,\partial B(0,2))\big]$, then there is
a set $Y=Y(x,r)$ of type 2, whose spine $L=L(x,r)$ goes
through $x$, such that
\begin{equation}\label{eqn22} 
D_{x,3r/2}(E,Y) \leq 24 \varepsilon
\ \hbox{ and }\
D_{x,r}(E_2,L) \leq 600 \varepsilon.
\end{equation}
\end{lemma}

\begin{proof} We already know from Lemma \ref{lem17} 
that $Z=Z(x,2r)$ is a set type 2 or 3,
whose spine $L$ passes at distance at most $34 \varepsilon r$ from $x$.
Let us check that
\begin{equation}\label{3.23}  
\hbox{there is a set $Y_1$ of type 2 such that }
Y_1 \cap B(x,155r/100) = Z \cap B(x,155r/100).
\end{equation}

This is clear if $Z$ is of type 2. If $Z$ is of type 3, and its center
lies out of $B(x,198r/100)$, Lemma \ref{L2.2} says that there is a 
set $Y_1$ of type 1 or 2 that coincides with $Z$ in 
$B(x,\frac{4}{5}\cdot\frac{198r}{100})$. Obviously, $Y_1$ is of type 2 
because the spine of $Z$ passes through $B(x,r)$,
and (\ref{3.23}) follows because
$\frac{155r}{100} < \frac{4}{5}\cdot\frac{198 r}{100}$. Finally suppose that
$Z$ is of type 3 and its center $z$ lies in $B(x,198r/100)$.
By definition of $Z$, we can find $y\in E \cap \overline B(z,2\varepsilon r)$; 
then we can use $Z+(y-z)$ to compute $c(y,r/200)$, and we get that
$c(y,r/200) \leq 200(2\varepsilon + 2 \varepsilon) < 10^{-4}$.
Lemma~\ref{lem12} gives a point $\xi \in E_3 \cap B(y,r/200)$. This
is impossible because $\dist(x,E_3) \geq 2r$ so (\ref{3.23}) holds.

Recall that $\dist(x,L) \leq 34 \varepsilon r$. Pick
$\ell \in L$ such that $|x -\ell| \leq 34 \varepsilon r$ and 
set $Y = Y_1 + (x-\ell)$; then $Y$ is a set of type 2 whose spine
goes through $x$. Let us check that
\begin{equation}\label{3.24}  
D_{x,3r/2}(E,Y) \leq 24\varepsilon.
\end{equation}
If $e\in E \cap B(x,3r/2)$, we can find $y\in Z$ such that
$|y-e| \leq 2 \varepsilon r$, then $y \in B(x,155r/100)$,
(\ref{3.23}) says that $y\in Y_1$ and $y'=y-(x-\ell) \in Y$.
So $\dist(e,Y) \leq |y'-e| \leq |y'-y| + |y-e| = |x -\ell| + |y-e| \leq
36 \varepsilon r$. Similarly, if $y' \in Y \cap B(x,3r/2)$,
then $y = y' +(x-\ell)$ lies in $Y_1 \cap B(x,155r/100)$, so we can
find $y\in E$ such that $|y-e| \leq 2 \varepsilon r$. Then
$\dist(y',E) \leq |y'-e| \leq 36 \varepsilon r$; (\ref{3.24})
follows.

Now the ball $B(x,3r/2)$ satisfies the hypotheses of Lemma \ref{lem14},
and this lemma says that $D_{x,r}(E_2,L) \leq 600 \varepsilon$. Since 
(\ref{3.24}) gives the first part of (\ref{eqn22}), 
Lemma \ref{lem21} follows. \qed
\end{proof}

The second part of (\ref{eqn22}) gives the local Reifenberg flatness of $E_2$
inside of $B(0,2)$ and away from $E_3$.
The next lemma says that near a point of $E_3$, $E_2$ looks like
the spine of a set of type 3 (that is, four half lines).

\begin{lemma}\label{lem24} 
Let $x\in E_3$ and $r>0$ be such that $B(x,2r) \subset B(0,199/100)$.
Then there is a set $T=T(x,r)$ of type 3 centered at $x$,
such that
\begin{equation}\label{eqn25} 
D_{x,4r/3}(E,T) \leq 15\varepsilon
\ \hbox{ and }\
D_{x,r}(E_2,L) \leq 220 \varepsilon,
\end{equation}
where $L$ denotes the spine of $T$.
\end{lemma}

\begin{proof} We proceed as in Lemma \ref{lem21}. 
Corollary \ref{cor11} says that $c(x,4r/3) \leq 15\varepsilon$, 
so we can find a set $T$ of type 3, centered at $x$, such that 
$D_{x,4r/3}(E,T) \leq 15 \varepsilon$.

Call $L$ the spine of $T$. Let us first check that
\begin{equation}\label{eqn26}   
\dist(y,L) \leq 100 \varepsilon r
\ \hbox{ for } y\in E_2 \cap B(x,r).
\end{equation}
Let $y\in E \cap B(x,r)$ be such that $\dist(y,L) \geq 100 \varepsilon r$, 
and call $Z=Z(y,100 \varepsilon r)$ the set promised by (\ref{eqn1}). 
Every point of $Z \cap B(y,100 \varepsilon r)$ lies within 
$100 \varepsilon^2 r$ of $E$, and hence within $21 \varepsilon r$ of $T$, 
by definition of $T$.
Since $\dist(y,L) \geq 100 \varepsilon r$, 
Lemma \ref{L3.4} says that $Z$ coincides with a 
plane $P$ in $B(y,60 \varepsilon r)$. 
In addition, $Z$ (and hence also $P$) goes through $y$, so we can use 
$P$ to compute $a(y,50\varepsilon r)$; we get that
$a(y,50\varepsilon r) \leq 2\varepsilon \leq 10^{-4}$, and
Lemma \ref{lem9} says that $y\in E_1$. This proves (\ref{eqn26}).

We still need to show that
\begin{equation}\label{eqn27}  
\dist(\xi,E_2) \leq 220 \varepsilon r
\ \hbox{ for } \xi\in L \cap B(x,r).
\end{equation}

First consider $\xi\in L \cap B(x,r) \setminus B(x,220\varepsilon r)$.
By definition of $T$ we can find 
$y\in E \cap \overline B(\xi,20\varepsilon r)$.
Then set $Z=Z(y,220\varepsilon r)$; let us check that
\begin{equation}\label{3.32}   
D_{\xi,90\varepsilon r}(Z,T) \leq \frac{21}{90} < 1/4.
\end{equation}

Every point of $Z \cap B(\xi,90 \varepsilon r)$
lies within $220 \varepsilon^2 r$ of $E$, hence within 
$220 \varepsilon^2 r + 20 \varepsilon r < 21 \varepsilon r$ of $T$
(by definition of $T$, and because $\xi \in B(x,r)$).
Similarly, if $t \in T \cap B(\xi,90 \varepsilon r)$,
we can find $e\in E$ such that $|e-t| \leq 20 \varepsilon r$,
then $|e-y| \leq |e-t|+|t-\xi|+|\xi-y| < 130 \varepsilon r$,
so $e \in E \cap B(y,130\varepsilon r)$ and we can find $z\in Z$
such that $|z-e| \leq 220 \varepsilon^2 r < \varepsilon r$; 
(\ref{3.32}) follows.

Recall that $T$ is of type 3 and $\xi$ lies on its spine;
then Lemma \ref{toro-lem1} says that 
$D_{\xi,90\varepsilon r}(T,P) \geq 1/4$ for any plane, so
(\ref{3.32}) forbids $Z$ to coincide with a plane in 
$B(\xi,120\varepsilon r)$.

By Lemma \ref{L3.4}, $Z$ to be of type 2 or 3, with a spine that meets 
the slightly larger ball $B(\xi,144\varepsilon r)$.
[For  sets of type $Gi$, replace $144$ with a geometric constant.]
Pick $\ell$ in the spine of $Z$, with 
$|\ell - \xi| \leq 144\varepsilon r$. Obviously
$\ell \in B(y,164 \varepsilon r)$ (because $|\xi - y| \leq 20 
\varepsilon r$), so we can pick $z_1\in E$
such that $|z_1-\ell| \leq 220 \varepsilon^2 r$. Let us check that
\begin{equation}\label{3.33} 
\hbox{ there is a set $Y$ of type 2 that coincides with
$Z$ in $B(z_1,24\varepsilon r)$.}
\end{equation}

This is clear if $Z$ is of type 2, so let us assume that it is
of type 3. Let $z_0$ denote its center. By Lemma \ref{L2.2},
it is enough to show that $|z_0 - z_1| \geq 32\varepsilon r$. 
Otherwise, 
\begin{equation}\label{3.34} 
|z_0 - y| \leq |z_0 - z_1| + |z_1 - \ell| + |\ell -y|
\leq 32\varepsilon r + 220 \varepsilon^2 r + 164 \varepsilon r
\leq 197 \varepsilon r. 
\end{equation}
[For sets of type $Gi$, we can leave $24$ as it is, replace
$32$ with $24C$, and also replace $197$ with a larger
geometric constant $a_{22}$; the next line forces us to replace
$220$ with a constant $a_{23} > a_{22}$.]

Since $z_0 \in Z \cap B(y,220\varepsilon r)$, 
we can find $e\in E$, with $|e-z_0| \leq 220 \varepsilon^2 r$.
The set $Z+(e-z_0)$ is of type 3, with a center at $e$, so we can use 
it to compute $c(e,\varepsilon r)$. Since 
$D_{y,220\varepsilon r}(Z,E) \leq \varepsilon$,
the usual argument shows that $c(e,\varepsilon r) 
\leq 440 \varepsilon$; then Lemma~\ref{lem12} gives 
a point $\eta \in E_3 \cap B(e,\varepsilon r)$. Notice that 
$|\eta-\xi| \leq |\eta-e| + |e-z_0| + |z_0-y| + |y-\xi| \leq 
\varepsilon r + 220 \varepsilon^2 r
+ 197\varepsilon r + 20 \varepsilon r \leq 219 \varepsilon r$,
by (\ref{3.34}) in particular. Since
$|\xi-x| \geq 220 \varepsilon r$, we get that $\eta \neq x$. At the same 
time, $\eta \in B(x,2r)$ because $\xi \in B(x,r)$,
so we have two points of $E_3$ in $B(x,2r) \i B(0,199/100)$.
Lemma \ref{lem13} says that this is impossible, so 
$|z_0 - z_1| \geq 32\varepsilon r$ and (\ref{3.33}) holds.

Recall from just above (\ref{3.33}) that 
$|z_1-\ell| \leq 220 \varepsilon^2 r$, and $\ell$ lies on the spine
of $Z$. Hence the spine of $Y' = Y + (z_1-\ell)$ goes through $z_1$,
and we can use $Y'$ to evaluate $b(z_1,23\varepsilon r)$. Recall
that (\ref{3.33}) says that $Y$ is $220 \varepsilon^2 r$-close
to $E$ near $B(z_1,23\varepsilon r)$; we get that
$b(z_1,23\varepsilon r) \leq 440 \varepsilon/23 < 25\varepsilon$. 
Thus we can apply Lemma \ref{lem14}, and there is a point 
$w\in E_2 \cap B(z_1,\varepsilon r)$. [For sets of type $Gi$,
we then need the constant $a_{15}$ in Lemma 4.5 to be large enough
compared to the replacement $a_{23}$ for $220$; this is not a serious 
problem, because $a_{23}$ is a geometric constant coming from 
Section \ref{geometry}.] Now
$|w-\xi| \leq |w-z_1| + |z_1 - \ell| + |\ell -\xi| \leq 
\varepsilon r + 220 \varepsilon^2 r + 144\varepsilon r
\leq 146 \varepsilon r$.

\smallskip
So we proved (\ref{eqn27}) when 
$\xi\in L \cap B(x,r) \setminus B(x,220\varepsilon r)$.
We can apply this to a point of 
$L \cap \partial B(x,250\varepsilon r)$, and we get the
existence of a point in $E_2 \cap B(x,500 \varepsilon r)$.
Also, the argument above applies to any ball $B(x,\rho)$ centered at $x$
such that $B(x,2\rho) \subset B(0,199/100)$, so by taking small values
of $\rho$, we find a sequence in $E_2$ that converges to $x$.

Now $x$ lies in the closure of $E_2$, so $\dist(\xi,E_2) \leq 220\varepsilon r$
for every $\xi \in B(x,220\varepsilon r)$.
This completes our proof of (\ref{eqn27}), and Lemma \ref{lem24} follows.
\qed
\end{proof}

\medskip
This ends our general description of $E_1$ and $E_2$.
As was suggested in Remark 1.1, 
we shall rapidly restrict to two simple situations. 
In addition to (\ref{eqn1}), we shall assume that either
\begin{equation}\label{eqn28}     
\hbox{$Z(0,2)$ is a set of type 2, whose spine $L$ contains $0$}
\end{equation}
(recall that $Z(x,r)$ is the set given by (\ref{eqn1})), or else
\begin{equation}\label{eqn29}      
\hbox{the origin lies in $E_3$, and $Z(0,2)$ is a set
of type 3 centered at $0$}.
\end{equation}
In these two cases, we shall establish the slightly more
precise (\ref{1.7}) and (\ref{1.8}) instead of (\ref{1.3})
and (\ref{eqn3}). 

Let us rapidly say what to do in the other cases. 
First suppose that we can find $x_0 \in E_3 \cap B(0,199/100)$.
By Lemma \ref{lem13}, $x_0$ is unique.
 
If $x_0=0$, Corollary \ref{cor11} says that $Z(0,2)$ is a set of type 3 
centered near the origin, and we may as 
well suppose that it is centered at $0$, as in (\ref{eqn29}) 
(and replace $\varepsilon$ with a slightly larger constant). 

If $x_0$ lies in $B(0,195/100)$, we can again assume that $Z(0,2)$ 
is a set of type 3 centered at $x_0$, and keep the proof below almost 
exactly as it is (just center the balls $B_{i_0}$ at $x_0$); we even get 
(\ref{1.7}) and (\ref{1.8}).

If $Z(0,2)$ is a set of type 3 and its center lies in $B(0,194/100)$,
Lemma \ref{lem12} gives a point in $E_3 \cap B(0,195/100)$, and
we are in one of the two previous cases. Otherwise, 
Lemma \ref{L2.2} says that $Z(0,2)$ coincides with a set $Z'$ of type 1 or 2 
in $B(0,{4 \cdot \frac{194}{500}}) \supset B(0,155/100)$. 
If $Z'$ is of type 2, the construction below restricted to $B(0,\frac{3}{2})$
yields the function $f$ promised in Theorem~\ref{T1.1}. This works 
independently of whether $0$ lies in the spine of $Z'$. If $Z'$ is of type 1
for $y\in B(0,\frac{154}{100})\cap E$, it is easy to see that 
$a(y,\frac{1}{100})<10^{-3}$. In this case, Lemma \ref{lem9} ensures that 
$a(y,t)\le 2\varepsilon$, for all $y\in B(0,\frac{154}{100})\cap E$ and $t$ 
small. Thus $B(0,\frac{154}{100})\cap E=\emptyset$, and we are in the standard 
Reifenberg situation.



Anyway we will just find it as convenient to restrict to
the case when $Z(0,2)$ is of type 2 whose spine contains the origin,
as in (\ref{eqn28}).

\medskip
Let us record here that 
\begin{equation}\label{4.37}      
E_3 \cap B(0,199/100) \hbox{ is empty when (\ref{eqn28}) holds.}
\end{equation}
Indeed, if $x\in E_3 \cap B(0,199/100)$, Corollary~\ref{cor11}
says that $c(x,1/200) \leq 15\varepsilon$, and Lemma~\ref{toro-lem1}
says that this is incompatible with the fact that 
$D_{0,2}(Z(0,2),E) \leq \varepsilon$ for a set $Z(0,2)$ of type~2.

Also recall that when (\ref{eqn29}) holds, Lemma~\ref{lem13} 
says that $0$ is the only point of $E_3 \cap B(0,199/100)$.

\section{Coverings and partitions of unity}\label{coverings}

Our general plan is to follow the same scheme as in the
standard proof of Reifenberg's theorem; 
thus we shall construct our parameterization by successive deformations 
of the set $Z(0,2)$. These deformations will occur near $E \cap B(0,2)$,
and appropriate partitions of unity play a key role.

We construct one such partition for each scale
$2^{-n}, n\geq 0$. We also construct the mapping
$f$ in a hierarchical way, which means that we define it first
on the spine of $Z(0,2)$, then on $Z(0,2)$ itself, and finally 
on the rest of $B(0,2)$. Our partitions of unity reflect this.

In the more general situation of Theorem \ref{T2.2}, only one minor 
modification is needed in this section. Since we have less control 
on the various constants that arise in the previous section, the
small security constant $10^{-20}$ that is used below needs to be 
replaced with a smaller constant, that depend on $n$, $d$, and
$\tau_0$.

For $n \geq 0$ given, we construct a ``good'' covering of $E$ at scale 
$2^{-n}$.
If (\ref{eqn29}) holds, we first cover $E_3 \cap B(0,199/100) = \{ 0 \}$
with a unique ball $B_{i_0} = B(0,2^{-n-20})$. For accounting
reasons, we set $I_3 = I_3(n) = \{ i_0 \}$. If (\ref{eqn28}) holds,
we simply take $I_3 = \emptyset$ and choose no ball.

Next we want to cover
\begin{equation}\label{eqn3.1} 
E'_2 = E_2 \cap B(0,198/100) \setminus \frac{7}{4}B_{i_0}
\end{equation}
(or just $E'_2=E_2\cap B(0,198/100)$ if (\ref{eqn28}) holds).
Here and below, $\lambda B$ is a notation for the 
ball with the same center as $B$ and $\lambda$ times the radius.
Select a maximal subset of $E'_2$, with the property that different
points of $E'_2$ lie at mutual distances at least $2^{-n-40}$.
Call $x_i$, $i\in I_2 = I_2(n)$ (with $I_2\cap I_3=\emptyset$), 
the points of this set, and set
$r_i=2^{-n-40}$ and $B_i = B(x_i,r_i)$ for $i\in I_2$. By maximality,
\begin{equation}\label{eqn3.2} 
\hbox{ the balls $\overline B_i$, $i\in I_2$, cover } E'_2.
\end{equation}
Then we take care of $E_1$. Set
\begin{equation}\label{eqn3.3}  
V_2 = \frac{15}{8}B_{i_0} \cup \bigcup_{i \in I_2} \frac{7}{4} B_i
\ \hbox{ and } \
E'_1 = E_1 \cap B(0,197/100) \setminus V_2
\end{equation}
(forget about $B_{i_0}$ if (\ref{eqn28}) holds), and pick a maximal
subset of $E'_1$ with the property that different points of $E'_1$
lie at mutual distances at least $2^{-n-60}$. Call $x_i$,
$i\in I_1 = I_1(n)$, with $I_1\cap (I_2\cup I_3)=\emptyset$, 
the points of this set, and set
$r_i=2^{-n-60}$ and $B_i = B(x_i,r_i)$ for $i\in I_1$. Then
\begin{equation}\label{eqn3.4}  
\hbox{ the balls $\overline B_i$, $i\in I_1$,
cover } E'_1.
\end{equation}
Let us go one more step, and set
\begin{equation}\label{eqn3.5}  
V_1 = \frac{31}{16}B_{i_0} \cup
\bigcup_{i \in I_2}  \frac{15}{8}B_{i} \cup
\bigcup_{i \in I_1} \frac{7}{4} B_i
\ \hbox{ and } \
E'_0 = \R^3 \setminus V_1,
\end{equation}
pick a maximal set in $E'_0$ with points at mutual distances at least
$2^{-n-80}$, call its points $x_i$, $i\in I_0 = I_0(n)$, with 
$I_0\cap (I_1\cup I_2\cup I_3)=\emptyset$ and set
$r_i=2^{-n-80}$ and $B_i = B(x_i,r_i)$ for $i\in I_0$. Again
\begin{equation}\label{eqn3.6}  
\hbox{ the balls $\overline B_i$, $i\in I_0$, cover } E'_0.
\end{equation}

Let us check that
\begin{equation}\label{eqn3.7}  
E_2 \cap B(0,198/100) \subset \frac{7}{4}B_{i_0} \cup
\bigg[ \bigcup_{i\in I_2} \overline B_i \bigg].
\end{equation}
Indeed, if $x\in E_2 \cap B(0,198/100)$ does not lie in
$\frac{7}{4}B_{i_0}$, then it lies in $E'_2$, and (\ref{eqn3.2})
gives the result. Similarly,
\begin{equation}\label{eqn3.8}  
E \cap B(0,197/100) \subset \frac{15}{8}B_{i_0}
\cup \bigg[\bigcup_{i \in I_2} \frac{7}{4} B_i \bigg]
\cup \bigg[ \bigcup_{i \in I_1} \overline B_i \bigg]
\end{equation}
because, if $x\in E \cap B(0,197/100)$ lies out of
$\frac{15}{8}B_{i_0} \cup \Big[ \bigcup_{i \in I_2} \frac{7}{4} B_i \Big]$,
then it lies in $E'_1$ by (\ref{eqn19}) and (\ref{eqn3.3}); 
thus (\ref{eqn3.8}) follows from (\ref{eqn3.4}). Finally,
\begin{equation}\label{eqn3.9} 
\frac{31}{16}B_{i_0} \cup
\bigg[ \bigcup_{i \in I_2}  \frac{15}{8}B_{i} \bigg]\cup
\bigg[ \bigcup_{i \in I_1} \frac{7}{4} B_i \bigg]
\cup \bigg[ \bigcup_{i \in I_0} \overline B_i \bigg]
= \R^3,
\end{equation}
again because the union of the first three pieces is $V_1$, and
by (\ref{eqn3.6}).

Now we define a partition of unity. Set
$I = I(n) = I_3 \cup I_2 \cup I_1 \cup I_0$. For each $i\in I$,
pick a smooth function $\widetilde\theta_i$ such that
\begin{equation}\label{eqn3.10}  
\widetilde\theta_i(x)=1 \hbox{ in } 2B_i, \
\widetilde\theta_i(x)= 0 \hbox{ out of } 3B_i, \hbox{ and }
0 \leq \widetilde\theta_i(x) \leq 1 \hbox{ everywhere}.
\end{equation}
We may choose the $\widetilde\theta_i$ as translations and
dilations of a same model. In any case they are chosen so that
\begin{equation} \label{5.11}     
|\nabla \widetilde\theta_i| \leq C 2^{n}\ \ \hbox{and}\ \ \
|\nabla^2 \widetilde\theta_i| \leq C 2^{2n}.
\end{equation}

Set $\Theta = \sum_{i\in I} \widetilde\theta_i \,$.
Since  by (\ref{eqn3.9}), $\{2B_i\}_{i\in I}$ covers
$\R^3$, $1\le\Theta(x)$. Moreover note that 
$\Theta(x) \leq C <\infty$ because the choice of balls ensures 
that they only overlap a bounded number of times.
Hence 
we can set
\begin{equation}\label{eqn3.11}    
\theta_i(x) = \widetilde\theta_i(x) /\Theta(x)
\ \hbox{ for } x\in \R^3.
\end{equation}

The usual computations yield
\begin{equation}\label{eqn3.12}     
\sum_{i\in I(n)} \theta_i(x) = 1,
\ |\nabla \theta_i| \leq C 2^{n}, \hbox{ and }
\ |\nabla^2 \theta_i| \leq C 2^{2n}
\ \hbox{ for } x\in \R^3.
\end{equation}

The point of our complicated choice of balls is that we
shall have a good control on the supports of the $\theta_i$.
For instance, if (\ref{eqn29}) holds,
\begin{equation}\label{eqn3.13}    
\dist(\frac{3}{2}B_{i_0}, 100B_i) \geq 2^{-n-30}
\hbox{ when } i\in I_2 \cup I_1 \cup I_0.
\end{equation}
Indeed we chose $x_i\not \in\frac{7}{4}B_{i_0}$,
then $2^{-n-30}$ and the radius of $100B_i$ are both
much smaller than one fourth of the radius of $B_{i_0}$
(which is $2^{-n-22}$). Similarly,
\begin{equation}\label{eqn3.14}  
\dist(\frac{3}{2}B_{j},100B_i) \geq 2^{-n-50}
\hbox{ when $j\in I_2$ and } i\in I_1 \cup I_0
\end{equation}
and
\begin{equation}\label{eqn3.15}   
\dist(\frac{3}{2}B_{j},100B_i) \geq 2^{-n-70}
\hbox{ when $j\in I_1$ and } i\in I_0.
\end{equation}
We also claim that
\begin{equation}\label{eqn3.16}   
\dist(\frac{7}{4}B_{i_0},100B_i) \geq 2^{-n-50}
\hbox{ when } i\in I_1 \cup I_0.
\end{equation}
Indeed, for $i\in I_1 \cup I_0$, we chose $x_i\not\in\frac{15}{8}B_{i_0}$.
Similarly,
\begin{equation}\label{eqn3.17}   
\dist(\frac{7}{4}B_{j},100B_i) \geq 2^{-n-70}
\hbox{ when $j\in I_2$ and } i\in I_0.
\end{equation}
Finally,
\begin{equation}\label{eqn3.18}   
\dist(\frac{15}{8}B_{i_0},100B_i) \geq 2^{-n-30}
\hbox{ when } i\in I_0,
\end{equation}
this time because $x_i\not \in\frac{31}{16}B_{i_0}$.

The inequalities (\ref{eqn3.13})-(\ref{eqn3.18}) have 
implications on the supports of the $\theta_i$.
For instance, (\ref{eqn3.13}) says that on $\frac{3}{2}B_{i_0}$ all the
$\theta_i$, $i\neq i_0$, vanish (and so $\theta_{i_0}(x)=1$).
Also,
\begin{equation}\label{eqn3.19}    
\theta_i(x)=0 \ \hbox{ when } \ \dist\Big(x,
\frac{7}{4}B_{i_0} \cup \bigcup_{j\in I_2} \frac{3}{2}B_{j}\Big)
\leq 2^{-n-70} \hbox{ and } i\in I_1 \cup I_0,
\end{equation}
by (\ref{eqn3.16}) and (\ref{eqn3.14}). By (\ref{eqn3.7}), this forces
\begin{equation}\label{eqn3.20}     
\theta_i(x)=0 \ \hbox{ when } \
\dist(x, E_2 \cap B(0,198/100)) \leq 2^{-n-70}
\hbox{ and } i\in I_1 \cup I_0,
\end{equation}
and hence, by (\ref{eqn3.12}),
\begin{equation}\label{eqn3.21}    
\sum_{i\in I_3 \cup I_2} \theta_i(x) = 1
\ \hbox{ when } \
\dist(x, E_2 \cap B(0,198/100)) \leq 2^{-n-70}.
\end{equation}
Finally,
\begin{equation}\label{eqn3.22}    
\theta_i(x)=0 \ \hbox{ when } \
\dist \Big( x, \frac{15}{8}B_{i_0} \cup \bigcup_{j\in I_2} \frac{7}{4}B_{j}
\cup \bigcup_{j\in I_1} \frac{3}{2}B_{j} \Big) \leq 2^{-n-70}
\ \hbox{ and } i\in I_0,
\end{equation}
by (\ref{eqn3.18}), (\ref{eqn3.17}), and (\ref{eqn3.15}), so (\ref{eqn3.8}) 
says that
\begin{equation}\label{eqn3.23}    
\theta_i(x)=0 \ \hbox{ when } \
\dist(x,E \cap B(0,197/100)) \leq 2^{-n-70}
\ \hbox{ and } i\in I_0,
\end{equation}
and then, by (\ref{eqn3.12}),
\begin{equation}\label{eqn3.24}    
\sum_{i\in I_3 \cup I_2 \cup I_1} \theta_i(x) = 1
\ \hbox{ when } \
\dist(x,E \cap B(0,197/100)) \leq 2^{-n-70}.
\end{equation}

\section{A first parameterization of $E_2$}\label{first}

We are now ready to start a construction of the mapping $f$.
To simplify the description, we start with the case when
(\ref{eqn28}) holds, i.e., when $Z(0,2)$ is of type 2. 
As we shall see in Section \ref{case},
the modifications that are needed when (\ref{eqn29}) holds
are mostly cosmetic.

We still want to proceed by layers, so we shall start with
a parameterization of a big part of $E_2$, which we shall mostly
be interested to define on
\begin{equation}\label{eqn4.1} 
\Gamma = L \cap \overline B(0,197/100)),
\end{equation}
where $L$ denotes the spine of $Z(0,2)$.

Recall from Lemma \ref{lem21}, (\ref{4.37}), and (\ref{eqn20}), 
that $E_2$ is a closed locally Reifenberg-flat 
set of dimension 1 in $B(0,199/100)$. This allows us to 
reproduce a standard Reifenberg argument and produce a 
parameterization for $E_2$. As some details of the construction
are needed later, we include it here.

Our parameterization $f^\ast$ is the limit of a sequence of mappings 
$f^\ast_n$ defined on $\Gamma$ and constructed by induction. We start with 
$f_0^\ast(x)=x$, and then set
\begin{equation}\label{eqn4.2}  
f^\ast_{n+1} = g^\ast_{n} \circ f^\ast_{n}
\ \hbox{ for } n\geq 0.
\end{equation}
In turn the deformation $g^\ast_{n}$ is defined by
\begin{equation}\label{eqn4.3} 
g^\ast_{n}(x) = \sum_{i\in I(n)} \theta_i(x) \psi_i^\ast(x),
\end{equation}
where $\theta_i$ is as in the partition of unity defined in the previous 
section, and
$\psi_i^\ast$ is a suitable deformation, which 
pushes points closer to $E_2$. Even though we mostly care about
the values of $g^\ast_{n}$ and the $\psi_i^\ast$ in a neighborhood
of $\Gamma_n = f^\ast_{n}(\Gamma)$, it will be just as easy to define 
them everywhere. We take $\psi_i^\ast(x)=x$ everywhere when 
$i\in I_1(n) \cup I_0(n)$; this does not matter too much, 
because we shall prove that $\Gamma_n$ stays 
very close to $E_2$ and hence (\ref{eqn3.20}) ensures that
the $\theta_i$'s, $i\in I_1(n) \cup I_0(n)$, vanish there.
We still need to define $\psi_i^\ast$ when $i\in I_2$.
The following notation will also be used in later sections.

Recall that when $i\in I_2(n)$, $x_i$ lies in
$E'_2 = E_2 \cap B(0,198/100) \setminus \frac{7}{4} B_{i_0}$
(but in this section there is no $B_{i_0}$), so Lemma \ref{lem21}
(applied to the pair $(x_i,10r_i)$, with $r_i=2^{-n-40}$)
tells us that there is a set $Y_i = Y(x_i,10 r_i)$
of type 2, whose spine $L_i$ goes through $x_i$, and such that
\begin{equation}\label{eqn4.4} 
D_{x_i,15 r_i}(E,Y_i) \leq 24\varepsilon
\ \hbox{ and } \
D_{x_i,10 r_i}(E_2,L_i) \leq 600\varepsilon.
\end{equation}
For the moment, we just need the second half of (\ref{eqn4.4}).
 
When $i\in I_2(n)$, we simply take for $\psi_i^\ast$ the orthogonal 
projection onto $L_i$. Thus we have defined $\psi_i^\ast$
for all $i\in I(n)$, and $g^\ast_{n}(x)$ and $f^\ast_{n}(x)$
are defined for all $x$ and $n$.
Several estimates are needed before we can continue our construction. 
First we claim that
\begin{equation}\label{6.5} 
D_{x_i,10 r_i}(L_i,L_j) \leq C \varepsilon
\ \hbox{ when $i,j \in I_2(n)$ and $8B_i$ meets $8B_j$.}
\end{equation}
Indeed suppose that $8B_i$ meets $8B_j$, and pick $x$ in the intersection.
Let us first estimate $D_{x,r_i}(L_i,L_j)$.
For $y\in L_i \cap B(x,r_i)$, (\ref{eqn4.4}) says that we can find
$z\in E_2$, with $|z-y|\leq 6000\varepsilon r_i$. Then $y\in 10B_j$,
and we can find $t\in L_j$ such that $|t-z|\leq 6000\varepsilon r_j
=6000\varepsilon r_i$. So $\dist(y,L_j) \leq 12000\varepsilon r_i$.
Similarly, $\dist(y,L_i) \leq 12000\varepsilon r_i$ for $y\in L_j \cap B(x,r_i)$.
So $D_{x,r_i}(L_i,L_j) \leq C \varepsilon$. 
Define $L_i$ by a point $y_i \in L_i \cap B(x,r_i)$ and a unit vector $v_i$, 
and proceed similarly with $L_j$; we can choose $y_j$ such that 
$|y_j-y_i| \leq C \varepsilon r_i$, and it is easy to see that 
$|v_j - v_i| \leq C \varepsilon$ or $|v_j + v_i| \leq C \varepsilon$.
Then (\ref{6.5}) follows. Next
\begin{equation}\label{eqn4.5} 
|\psi_i^\ast(x) - \psi_j^\ast(x)| \leq C \varepsilon 2^{-n}
\ \hbox{ when $i,j \in I_2(n)$ and $x \in 8B_i \cap 8B_j$.}
\end{equation}
Indeed, $8B_i$ meets $8B_j$, so $D_{x_i,10 r_i}(L_i,L_j) \leq C \varepsilon$
by (\ref{6.5}). Again it is a matter of elementary geometry that the
orthogonal projections $\psi_i^\ast(x)$ and $\psi_j^\ast(x)$ are 
close to each other (for instance, write $L_i = y_i + \R v_i$ and 
$L_j = y_j + \R v_j$ as above and compute). For the same reason,
the (constant) derivatives of $\psi_i^\ast$ and $\psi_j^\ast$ are
$C\varepsilon$-close, i.e.,
\begin{equation}\label{eqn4.6} 
|D\psi_i^\ast - D\psi_j^\ast| \leq C \varepsilon
\ \hbox{ when $i,j \in I_2(n)$ and $8B_i$ meets $8B_j$.}
\end{equation}

For each $x$ such that
\begin{equation}\label{eqn4.7} 
\dist(x,E_2 \cap B(0,198/100)) \leq 2^{-n-70},
\end{equation}
(\ref{eqn3.20}) tells us that all the $\theta_i(x)$ with $i\notin I_2(n)$,
vanish at $x$. This allows us to pick $j(x) \in I_2(n)$ such
that $\theta_{j(x)}(x) \neq 0$, use (\ref{eqn4.5}) and (\ref{eqn3.10})
to show that
$|\psi_i^\ast(x) - \psi_{j(x)}^\ast(x)| \leq C \varepsilon 2^{-n}$
whenever $\theta_{i}(x) \neq 0$, and get that
\begin{equation}\label{eqn4.8}  
|g_n^\ast(x) - \psi_{j(x)}^\ast(x)| \leq C_1 \varepsilon 2^{-n},
\end{equation}
because $g_n^\ast(x)$ is an average of the 
$\psi_i^\ast(x)$. Here $C_1$, just like the next constants $C_j$, 
is a simple geometric constant.
Since $\psi_{j(x)}^\ast(x) \in L_{j(x)} \cap 3B_{j(x)}$,
(\ref{eqn4.4}) says that
\begin{equation}\label{eqn4.9}  
\dist (g_n^\ast(x),E_2) \leq C_2\varepsilon 2^{-n-1}.
\end{equation}

Let $y\in E_2$ minimize the distance to $x$; notice that
$y\in 6B_{j(x)}$ because $x\in 3B_{j(x)}$ and hence 
$\dist(x,E_2) \leq |x-x_{j(x)}| \leq 3r_{j(x)}$. Then by (\ref{eqn4.4})
\begin{equation}\label{6.11} 
|\psi_{j(x)}^\ast(x)-x| = \dist(x,L_{j(x)}) 
\leq |x-y|+ \dist(y,L_{j(x)}) 
\leq \dist(x,E_2) + 6000 \varepsilon r_{j(x)}.
\end{equation}
Hence by (\ref{eqn4.8}) and (\ref{6.11})
\begin{equation}\label{6.12} 
|g_n^\ast(x)-x| \leq \dist(x,E_2) + C_3 \varepsilon 2^{-n}.
\end{equation}
We shall remember that all this happens for 
$x\in B(0,2)$ such that (\ref{eqn4.7}) holds.

We are ready to show by induction that if $z\in \Gamma = L \cap 
\overline B(0,197/100)$, then
\begin{equation}\label{eqn4.10} 
\dist (f_n^\ast(z),E_2) \leq C_2\varepsilon 2^{-n},
\end{equation}
and
\begin{equation}\label{eqn4.11} 
|f_{n+1}^\ast(z)-f_n^\ast(z)| \leq (C_2+C_3) \varepsilon 2^{-n}
\end{equation}
for every $n \geq 0$.

Indeed we know that (\ref{eqn4.10}) holds for $n=0$ (if we did not choose
$C_2$ too small). Also, suppose that (\ref{eqn4.10}) holds for $n$ and 
(\ref{eqn4.11}) holds for all $m < n$ (a vacuous condition when $n=0$). 
Then $|f_n^\ast(z)-z| \leq 4 (C_2+C_3) \varepsilon$ (by repeated use 
of (\ref{eqn4.11})), so $f_n^\ast(z)$ lies inside or very close to 
$B(0,197/100)$. Hence $\dist(f_n^\ast(z),E_2 \cap B(0,198/100)) 
\leq C_2\varepsilon 2^{-n} \leq 2^{-n-70}$ (by (\ref{eqn4.10})
and if $\varepsilon$ is small enough; our bound of $10^{-25}$
should be more than enough). That is, 
\begin{equation}\label{6.15} 
x=f_n^\ast(z) \hbox{ satisfies (\ref{eqn4.7}).}
\end{equation}
Then $\dist (f_{n+1}^\ast(z),E_2) =
\dist (g_n^\ast(f_n^\ast(z)),E_2) \leq C_2\varepsilon 2^{-n-1}$
by (\ref{eqn4.9}). That is, (\ref{eqn4.10}) holds for $n+1$. In addition, 
$|f_{n+1}^\ast(z)-f_n^\ast(z)| = |g_n^\ast(x)-x| \leq 
\dist(x,E_2) + C_3 \varepsilon 2^{-n} \leq (C_2+C_3) \varepsilon 2^{-n}$
by (\ref{6.12}) and (\ref{eqn4.10}).
Thus (\ref{eqn4.10}) and (\ref{eqn4.11}) hold for every $n$.

\smallskip
By (\ref{eqn4.11}), $\{ f_n^\ast \}$ converges uniformly on $\Gamma$
to some limit $f^\ast$, and
\begin{equation}\label{eqn4.12} 
||f^\ast -  f_n^\ast||_{\infty} \leq 2 (C_2+C_3) \varepsilon 2^{-n}.
\end{equation}

Next we examine the Lipschitz properties of
$\Gamma_n = f_n^\ast(\Gamma)$.

\begin{lemma}\label{lem4.13}  
The restriction of $f_n^\ast$ to $\Gamma$ is of class $C^2$, 
with a derivative that does not vanish.
For each $i\in I_2(n)$, $\Gamma_n \cap 5B_i$ is contained
in a $C_4 \varepsilon$-Lipschitz graph $G_{n,i}$ over $L_i$,
that meets $B(x_i,C_5 \varepsilon 2^{-n})$.
\end{lemma}

\begin{proof} Let us prove all this by  induction. 
The case when $n=0$ is clear, because $F_0^{\ast}=id$, 
$\Gamma$ is a straight line, and its closeness to all $L_i$,
$i\in I_2(0)$, is an easy consequence of (\ref{eqn28}) and (\ref{eqn4.4}).
So let us assume that we know the lemma for $\Gamma_n$, 
and prove it for $\Gamma_{n+1}=g_n^\ast(\Gamma_n)$.

We start with the Lipschitz description.
Let $i\in I_2(n+1)$ be given. Recall that 
$x_i \in E'_2 = E_2 \cap B(0,198/100)$, by (\ref{eqn3.1}). 
Then (\ref{eqn3.7}) ensures that
$x_i \in \overline B_j$ for some $j\in I_2(n)$. Here the situation
is a little simpler because $I_3(n)=\emptyset$, but the argument
will be very similar when $I_3(n) \neq\emptyset$. By the induction
assumption, there is a $C_4 \varepsilon$-Lipschitz graph $G_{n,j}$
over $L_j$ that contains $\Gamma_n \cap 5B_j$. Notice that
$5B_i \i 4B_j$, so $g_n^\ast(\Gamma_n \setminus 5B_j)$ does not 
meet $5B_i$, by (\ref{eqn4.11}). Thus we may restrict our attention
to points coming from $\Gamma_n \cap 5B_j \i G_{n,j}\cap 5B_j$.

Let $z\in G_{n,j}\cap 5B_j$ be given, and denote by $\tau$
a unit tangent vector to $G_{n,j}$ at $z$. Observe that
\begin{equation}\label{6.14}  
\dist(z,E_2) \leq C' \varepsilon 2^{-n},
\end{equation}
where $C'$ depends on $C_4$ and $C_5$. Indeed 
$\dist(z,L_j) \leq C_5 \varepsilon 2^{-n} + 11C_4 \varepsilon 2^{-n-40}$,
because $G_{n,j}$ is a $C_4 \varepsilon$-Lipschitz graph that goes through 
$B(x_j, C_5 \varepsilon 2^{-n})$, and $x_j \in L_j$ by definition
of $Y_j$ and $L_j$. Now (\ref{6.14}) follows from (\ref{eqn4.4}).

Let $k\in I(n)$ be such that $\theta_k(z) \neq 0$. 
We know that $k\in I_2(n)$, by (\ref{6.14}) and (\ref{eqn3.20}).
More precisely, we should make sure that $C'\varepsilon$ in (\ref{6.14}) 
is less
than $2^{-70}$. We claim that this is the case, but the reader 
should not worry, otherwise we could have chosen a smaller constant 
than $10^{-25}$ as an upper bound for $\varepsilon$. Then
by (\ref{eqn4.5}) and (\ref{eqn4.6})
\begin{equation}\label{eqn4.14} 
|D\psi_k^\ast- D\psi_j^\ast| \leq C \varepsilon,
\ \hbox{ and } \ 
|\psi_k^\ast- \psi_j^\ast| \leq C \varepsilon 2^{-n}
\hbox{ in $5B_j$.} 
\end{equation}
Since  
$Dg_n^\ast(z) = \sum_{k\in I(n)} \theta_k(z) D\psi_k^\ast(z)
+ \sum_{k\in I(n)} D\theta_k(z) \psi_k^\ast(z)$
by (\ref{eqn4.3}),
\begin{equation}\label{eqn4.16} 
|Dg_n^\ast(z)- D\psi^\ast_j|
\leq \sum_{k\in I(n)} \theta_k(z) |D\psi_k^\ast - D\psi^\ast_j|
+ \sum_{k\in I(n)} |D\theta_k(x)| |\psi_k^\ast(z) - \psi^\ast_j(z)|
\leq C_6 \varepsilon
\end{equation}
because $\sum_k \theta_k =1$ and $\sum_k D\theta_k =0$, 
and by (\ref{eqn3.12}) and (\ref{eqn4.14}). Here $C_6$
is a geometric constant, and in particular does not depend on our 
future choice of $C_5$.

Since $G_{n,j}$ is a $C_4 \varepsilon$-Lipschitz graph over $L_j$
and $\psi^\ast_j$ is the orthogonal projection onto $L_j$,
$|\tau-D\psi^\ast_j(\tau)| \leq C_4 \varepsilon$, hence
\begin{equation}\label{eqn4.17} 
|Dg_n^\ast(z)(\tau) - \tau|
\leq |Dg_n^\ast(z)(\tau) - D\psi^\ast_j(\tau)| + C_4 \varepsilon
\leq (C_4 + C_6) \varepsilon
\end{equation}
by (\ref{eqn4.16}). In particular,
$|Dg_n^\ast(z) (\tau)| \geq 1-(C_4 + C_6) \varepsilon\geq 99/100$.

Set $v=Dg_n^\ast(z)(\tau)$; since $v\neq 0$, it is a tangent
vector to $g_n^\ast(G_{n,j})$ locally . Recall that $\tau$ is the unit 
tangent vector  to $G_{n,j}$ at $x$.

we want to show that it makes a small angle with $L_j$. 
Call $\pi^\perp_j$ the projection in the direction orthogonal to $L_j$; 
then $\pi^\perp_j(v) = \pi^\perp_j (Dg_n^\ast(z)(\tau))
= \pi^\perp_j (Dg_n^\ast(z)(\tau) - D\psi^\ast_j(\tau))$, hence
$|\pi^\perp_j(v)| \leq C_6 \varepsilon$ by (\ref{eqn4.16}).
Since $|v| \geq 99/100$, $v$ makes a small angle with $L_j$. 
Also,
\begin{equation}\label{6.21} 
D_{x_j,10r_j}(L_i,L_j) \leq C \varepsilon,
\end{equation}
by the proof of (\ref{6.5}). [The small difference is that now $j\in I_2(n)$ and 
$i\in I_2(n+1)$, so $B_j$ is twice as big as $B_i \,$; otherwise, we still have that 
$x_i \in B_i \cap \overline B_j$, which is more than enough for the 
proof to go through.] 

Thus $L_i$ makes a small angle with $L_j$
and $v$ makes a small angle with $L_i$. We deduce from this 
(and the fact that $g_n^\ast$ is smooth on $G_{n,j} \cap 5B_j$, with a 
non-vanishing derivative) that
\begin{equation}\label{eqn4.19} 
g_n^\ast(G_{n,j} \cap 5B_j)
\hbox{ is (contained in) a $C_7 \varepsilon$-Lipschitz graph
$G_{n+1,i}$ over }L_i \, ,
\end{equation}
with $C_7$ depending only on $C_6$ and other geometric constants.
Recall that we already checked that every point of $\Gamma_{n+1} \cap 5B_i$
lies in $g_n^\ast(G_{n,j} \cap 5B_j)$,
so we just need to take $C_4$ larger than $C_7$  to get the
part of the induction that says that $\Gamma_{n+1} \cap 5B_i$
is contained in a $C_4 \varepsilon$-Lipschitz graph over $L_i$. 

\smallskip
Next we check that $G_{n+1,i}$ meets $B(x_i,C_5 \varepsilon 2^{-n-1})$.
By the induction hypothesis, we can find 
$z \in G_{n,j} \cap B(x_j,C_5 \varepsilon 2^{-n})$.
Set $z_1 = g_n^\ast(z)$; thus $z_1 \in g_n^\ast(G_{n,j} \cap 5B_j)
\i G_{n+1,i}$, by (\ref{eqn4.19}). Call $z'=\psi_j^\ast(z)$ the
orthogonal projection of $z$ onto $L_j$. By (\ref{eqn4.14}),
$|\psi_k^\ast(z)-z'| \leq C \varepsilon 2^{-n}$
when $\theta_k(z) \neq 0$, and hence
$|z_1-z'| = |g_n^\ast(z)-z'| \leq C \varepsilon 2^{-n}$,
with a constant $C$ that does not depend on $C_5$.
This proves that $z_1 \in B(x_j,(C_5+C) \varepsilon 2^{-n}) \i B_j$
(because $z' \in B(x_j,C_5 \varepsilon 2^{-n})$, just like $z$ 
as $x_j\in L_j$),
and that $\dist(z_1,L_j) \leq C \varepsilon 2^{-n}$.
By (\ref{6.21}), $\dist(z_1,L_i) \leq C \varepsilon 2^{-n}$.
By definition of a graph, there is a point $z_2 \in G_{n+1,i}$, with 
$\psi^\ast_{j}(z_2) = x_i$. Then $\dist(z_2,L_i) \leq C\varepsilon 2^{-n}$
(because $\dist(z_1,L_i) \leq C \varepsilon 2^{-n}$ and the slope between 
$z_1$ and $z_2$ is small). Then $z_2 \in B(x_i,C \varepsilon 2^{-n})$.
The constant $C$ that we get does not depend on $C_5$, so we can choose 
$C_5 \geq C$, and $G_{n+1,i}$ meets $B(x_i,C_5 \varepsilon 2^{-n-1})$
as needed.

\smallskip
Finally we need to check that the derivative of the restriction of 
$f^\ast_{n+1}$ to $\Gamma$ does not vanish. Let $t\in \Gamma$ be 
given, and set $z=f^\ast_{n}(t)$. By (\ref{eqn4.11}),
$|f_n^\ast(t) - t| \leq 2(C_2+C_3) \varepsilon$, so 
$f_n^\ast(t)$ lies well inside $B(0,198/100)$ (because
$z\in \overline B(0,197/100)$ by (\ref{eqn4.1}).
By (\ref{eqn3.7}), $z \in \overline B_j$ for some $j\in I_2(n)$.
We can use the graph $G_{n,j}$ provided by our induction assumption.
Call $\tau_0$ a unit tangent vector to $\Gamma$  and set
$\tau = Df^\ast_{n}(t)(\tau_0)$. By induction assumption, 
$\tau \neq 0$. Then $\tau$ is a tangent vector to $G_{n,j}$.
We can follow the argument between (\ref{6.14})
and (\ref{eqn4.17}), and we get that $Dg_n^\ast(f_n^{\ast}(t))(\tau) \neq 0$.
Thus, $Df^\ast_{n+1}(t)(\tau_0) \neq 0$, as needed.
This completes our proof of Lemma~\ref{lem4.13}.
\qed
\end{proof}

Next we claim that if $n \geq 0$, $i\in I_2(n)$, and
$x_i \in B(0,196/100)$,
\begin{equation}\label{eqn4.20}  
\Gamma_n \cap 5B_i = G_{n,i} \cap 5B_i,
\end{equation}
where $G_{n,i}$ is still as in Lemma \ref{lem4.13}.

First recall from Lemma \ref{lem4.13} that $\Gamma_{n}$ is a $C^1$ curve,
with just two extremities. If none of the two extremities lies in 
$5B_i$, then only two things can happen. Either $\Gamma_n \cap 5B_i$
is empty, or else $\Gamma_n$ enters $5B_i$ at some point. Since $\Gamma_n \cap 
5B_i \i G_{n,i}$ and the derivative of $f_n^\ast$ on $\Gamma$ does 
not vanish, $\Gamma_n$ has to follow $G_{n,i}$ without turning back, 
until it eventually leaves $5B_i$. Then it runs along $G_{n,i} \cap 5B_i$
the whole way through, and (\ref{eqn4.20}) holds.

It will be easier to check by induction that (\ref{eqn4.20})
holds for $i\in I_2(n)$ such that
$\dist(x_i,\R^3 \setminus B(0,197/100))\geq \sum_{m=0}^{n} 2^{-m-40}$.
The two extremities of $\Gamma_{n}$ lie within $2(C_2+C_3) \varepsilon$
of $\partial B(0,197/100)$, by (\ref{eqn4.1}) and (\ref{eqn4.12}), so
they lie out of $5B_i$ and we just need to exclude the case when 
$\Gamma_n \cap 5B_i$ is empty. 

When $n=0$, $\Gamma_n=\Gamma$ meets $B_i$. Indeed
$x_i \in E_2\cap B(0,197/100)$, hence Lemma \ref{lem17} says that
$Z(x_i,10^{-2})$ is of type 2 or 3, with a spine that passes at distance
at most $17 \varepsilon 10^{-2}$ from $x_i \,$; this forces $L$ to pass 
through $\frac{1}{2}B_i$, because otherwise $Z(0,2)$ coincides with 
a plane in $\frac{5}{12}B_i$ (by Lemma \ref{L3.4}), and then 
$D_{x_i,10^{-3}}(Z(x_i,10^{-2}),Z(0,2)) \geq 1/5$ (by Lemma \ref{toro-lem1}); 
this is impossible because $Z(x_i,10^{-2})$ and $Z(0,2)$ are both so close 
to $E$ near $B(x_i,10^{-3})$ (by (\ref{eqn1})). So $L$ meets $\frac{1}{2}B_i$ and 
$\Gamma$ meets $B_i$ and (\ref{eqn4.20}) holds for $n=0$.

Finally let us assume our claim for some $n\geq 0$ and prove it for $n+1$.
Let $i\in I_2(n+1)$ be as above. By (\ref{eqn3.1}), 
$x_i \in E'_2 = E_2 \cap B(0,198/100)$; by (\ref{eqn3.7}), 
$x_i \in \overline B_j$ for some $j\in I_2(n)$.
Thus $|x_j-x_i| \leq 2^{-n-40}$ (because $x_i \in \overline B_j$), hence
\begin{equation}\label{eqn4.21}  
\dist(x_j,\R^3 \setminus B(0,197/100)) \geq
\dist(x_i,\R^3 \setminus B(0,197/100)) - 2^{-n-40}
>  \sum_{m=0}^{n} 2^{-m-40}.
\end{equation}
By induction assumption, (\ref{eqn4.20}) holds for $n$ and $j$, and in 
particular we can find $\xi\in \Gamma_n \cap B(x_j,C_5 \varepsilon 2^{-n})$
(see Lemma \ref{lem4.13}). Now $g_n^\ast(\xi) \in \Gamma_{n+1} \cap 5B_i$, 
by (\ref{eqn4.11}). Thus $\Gamma_{n+1}\cap 5B_i\not =\emptyset$
 and (\ref{eqn4.20}) holds for $B_i$, as needed.
This completes our proof of (\ref{eqn4.20}) by induction.

\medskip 
Next we want to show that $f^\ast$ is biH\"older.
Let us check that
\begin{equation}\label{eqn4.22} 
(1-C\varepsilon) \dist(y,z)
\leq \dist(g_n^\ast(y),g_n^\ast(z))
\leq (1+C\varepsilon) \dist(y,z)
\end{equation}
when $y$, $z\in \Gamma_n$ are such that $|y-z| \leq 2^{-n-40}$.
By (\ref{eqn4.10}), $\dist(y,E_2) \leq C_2 \varepsilon 2^{-n}$; more
trivially, $y$ lies in $B(0,197/100)$ or very close to it.
Then (\ref{eqn3.7}) says that $y$ lies within $C\varepsilon 2^{-n}$ of
$B_i$ for some $i \in I_2(n)$.
Let $G_{n,i}$ be the Lipschitz graph of Lemma \ref{lem4.13}; since both
$y$ and $z$ lie on $G_{n,i}\cap 3B_i$, we can replace $\dist(y,z)$
with the length of the arc of $\Gamma_n$ between these points, and
make a relative error less than $C\varepsilon$. We can proceed
similarly with $\dist(g_n^\ast(y),g_n^\ast(z))$ and the arc-length
on $g_n^\ast(G_{n,j} \cap 5B_j)$ (for the appropriate $j$), which is also a
$C\varepsilon$-Lipschitz graph by (\ref{eqn4.19}). Finally, the ratio
between the arc-lengths can be computed in terms of the derivative
$|Dg_n^\ast(z)(\tau(z))|$ on $\Gamma_n$. Thus using (\ref{eqn4.17}),
we get (\ref{eqn4.22}).
 
Our biH\"older estimate will follow from (\ref{eqn4.22}) by
a rather mechanical argument, which will be repeated later in this 
paper.
Let $y$, $z\in \Gamma$ be given, and  set $y_n= f_n^\ast(y)$ and 
$z_n= f_n^\ast(z)$. Assume that $|z-y|<1$; we shall 
take care of the other case later.
As long as $|y_n-z_n| \leq 2^{-n-40}$ we can apply (\ref{eqn4.22}) and get that
\begin{equation}\label{eqn4.23} 
(1-C\varepsilon)^{n+1} |y-z| \leq |y_{n+1}-z_{n+1}|
\leq (1+C\varepsilon)^{n+1} |y-z|.
\end{equation}
How long can this last? If (\ref{eqn4.23}) holds for $n-1$ and
$|y_n-z_n| \leq 2^{-n-40}$, we have that
\begin{equation}\label{6.27} 
(1-C\varepsilon)^{n} |y-z| \leq |y_n-z_n| \leq 2^{-n-40}
\end{equation}
hence $n \log_2(1-C\varepsilon) + \log_2(|y-z|) \leq -n-40 \leq -n$,
or equivalently 
\begin{equation}\label{6.28a} 
\log_2(\frac{1}{|y-z|}) \geq n + n \log_2(1-C\varepsilon)
\geq n (1-C\varepsilon).
\end{equation}

This cannot happen for $n$ large, thus there is a smallest $n_0$ such that
$|y_{n_0}-z_{n_0}| > 2^{-n_0-40}$, and we even know that
(\ref{6.28a}) holds for $n_0-1$, so
\begin{equation}\label{6.29a} 
n_0 \leq 1 + (1-C\varepsilon)^{-1} \log_2(\frac{1}{|y-z|})
\leq 1 + (1+C'\varepsilon)\log_2(\frac{1}{|y-z|}).
\end{equation}

Now we can use (\ref{eqn4.11}) to say that
\begin{equation}\label{eqn4.24} 
\big| |y_n-z_n|-|y_{n_0}-z_{n_0}| \big| \leq C \varepsilon 2^{-n_0}
\end{equation}
for $n > n_0$, which leads to
\begin{equation}\label{6.28} 
\big| |f^\ast(y)-f^\ast(z)|-|y_{n_0}-z_{n_0}| \big| \leq C \varepsilon 2^{-n_0}
\end{equation}
and then, since $|y_{n_0}-z_{n_0}| > 2^{-n_0-40}$, to
\begin{equation}\label{eqn4.25}
(1-C\varepsilon) |y_{n_0}-z_{n_0}| \leq |f^\ast(y)-f^\ast(z)|
\leq (1+C\varepsilon) |y_{n_0}-z_{n_0}|
\end{equation}
and 
\begin{equation}\label{eqn4.26} 
(1-C\varepsilon)^{n_0+1} |y-z| \leq |f^\ast(y)-f^\ast(z)|
\leq (1+C\varepsilon)^{n_0+1} |y-z|,
\end{equation}
by (\ref{eqn4.23}). That is,
$\left| \log \frac{|f^\ast(y)-f^\ast(z)|}{|y-z|} \right|
\leq C \varepsilon (n_0+1)
\leq C \varepsilon \left[ 2 + (1 + 
C'\varepsilon)\log_2(\frac{1}{|y-z|})\right]
\leq C'' \varepsilon \big[ (1+ \log_2(\frac{1}{|y-z|}) \big]$,
by (\ref{6.29a}). Then we take exponentials and obtain that
\begin{equation}\label{eqn4.27} 
(1-C\varepsilon) \, |y-z|^{1+C\varepsilon}
\leq |f^\ast(y)-f^\ast(z)| \leq
(1+C\varepsilon) \, |y-z|^{1-C\varepsilon}.
\end{equation}
For $y,z\in\Gamma$ with $1<|z-y|<4$ (\ref{eqn4.11}) or (\ref{eqn4.12}) imply
\begin{equation}\label{eqn4.27A} 
(1-C\varepsilon) \, |y-z|
\leq |f^\ast(y)-f^\ast(z)| \leq
(1+C\varepsilon) \, |y-z|,
\end{equation}
which yields (\ref{eqn4.27}).
Thus 
$f^\ast$ is bi-H\"older on $\Gamma$, with exponents that
are as close to $1$ as we want.

Observe that it is clear from (\ref{eqn4.10}) that 
$f^\ast(\Gamma)\i E_2$. Also,
\begin{equation}\label{eqn4.28}
f^\ast(\Gamma) \hbox{ contains } E_2 \cap B(0,196/100).
\end{equation}
Indeed let $z\in E_2 \cap B(0,196/100)$; for every large $n$,
$z$ lies in $\overline B_i$ for some $i\in I_2(n)$, 
$x_i\in B(0,196/100)$ too, and (\ref{eqn4.20}) says that $B_i$ 
meets $\Gamma_n$ (recall from Lemma \ref{lem4.13} that $G_{n,i}$
meets $B(x_i,C\varepsilon 2^{-n})$). Thus $\dist(z,\Gamma_n) \leq 2^{-n-20}$, 
and $z \in f^\ast(\Gamma)$ by compactness.

\section{A parameterization of $E$ when $E_3 = \emptyset$}\label{param}

We still assume that (\ref{eqn28}) holds, so that in particular
$E_3 \cap B(0,199/100)= \emptyset$ by (\ref{4.37}), and we want to define 
the restriction to a big part of $Z(0,2)$ of the parameterization $f$. 
As before, we shall define a mapping everywhere, but we shall only care 
about $f(z)$ when $z$ lies in the set
\begin{equation}\label{7.1}  
Y = Z(0,2) \cap B(0,195/100).
\end{equation}
We still want to get $f$ as the limit of mappings $f_n$, where $f_0(x)=x$,
and
\begin{equation}\label{eqn5.1}  
f_{n+1} = g_{n} \circ f_{n} \ \hbox{ for } n\geq 0,
\hbox{ with }
g_n(x) = \sum_{i\in I(n)} \theta_i(x) \psi_i(x)
\end{equation}
and where the $\psi_i$, $i\in I(n)$, are suitable deformations that will
be defined soon. We can again take $\psi_i(x)=x$ for $i\in I_0(n)$ 
(but $f_n(Y)$ will avoid the support of $\psi_i$ anyway), 
but this time we need to do something  when $i\in I_1(n)$.
Let us first define $\psi_i$ in this case.

Recall that when $i\in I_1(n)$, $B_i$ is a ball of radius
$r_i=2^{-n-60}$ centered at some $x_i\in E'_1 \subset E \cap B(0,197/100)$;
see near (\ref{eqn3.3}). Also recall from (\ref{eqn3.14}) and (\ref{eqn3.16}) that
\begin{equation}\label{eqn5.2}  
\dist\Big(100B_i,\frac{7}{4}B_{i_0} \cup \bigcup_{j\in I_2} \frac{3}{2} B_j
\Big) \geq 2^{-n-50},
\end{equation}
where we only mention $i_0$ here to convince the reader that things
will also work when (\ref{eqn29}) holds. Then (\ref{eqn3.7}) says that
\begin{equation}\label{eqn5.3} 
\dist(100B_i,E_2) \geq 2^{-n-50},
\end{equation}
where we also use the fact that $x_i \in B(0,197/100)$.
Let us check that
\begin{equation}\label{eqn5.4} 
Z(x_i,100r_i) \hbox{ coincides with a plane $P_i$ in }
B(x_i,60r_i).
\end{equation}

If $Z=Z(x_i,100r_i)$ is of type 3 and its center lies in $B(x_i,98r_i)$, 
we can choose $\xi \in E$ close to the center, apply Lemma \ref{lem12}
to $B(\xi,r_i)$, and we get that $B(\xi,r_i)$ meets $E_3$. This
is impossible, by (\ref{eqn5.2}) (or (\ref{eqn5.3}) and the fact that
points of $E_3$ are always limits of points of $E_2$).

Otherwise, and if (\ref{eqn5.4}) fails, Lemma \ref{L3.4} says that 
the spine of $Z$ meets $\overline B(x_i,80r_i)$. 
Pick $\xi \in E \cap B(x_i,81r_i)$, close to the spine of $Z$. 
Then $Z$ coincides with a set of type $2$ in $B(\xi,2r_i)$, we can 
apply Lemma \ref{lem14} to $B(\xi,r_i)$, and we get that $B(\xi,r_i)$
meets $E_2$. This is impossible too, by (\ref{eqn5.3}). So (\ref{eqn5.4})
holds. For the generalization in Theorem \ref{T2.2}, we would just 
need to replace 100 with a larger constant and make $2^{-20}$ and
$\varepsilon$ somewhat smaller.

By definition, $Z(x_i,100r_i)$ passes through $x_i$. We take
\begin{equation}\label{eqn5.5}  
\psi_i = \pi_i,  \hbox{ the orthogonal projection onto $P_i$,
when } i\in I_1.
\end{equation}

Now we need to define $\psi_i$ when $i\in I_2(n)$. There is a
constraint, because since we want to have $f=f^\ast$ on $\Gamma$,
we shall demand that
\begin{equation}\label{eqn5.6} 
\psi_i(x) = \psi_i^\ast(x)
\hbox{ for } x\in \Gamma_n \cap \support(\theta_i).
\end{equation}
This will suffice, because if we have (\ref{eqn5.6}) for all $i\in I_2(n)$, 
then $g_n=g_n^\ast$ on $\Gamma_n$, since $\Gamma_n$ does not meet the
support of the $\theta_i$, $i\in I_1(n)\cup I_0(n)$, by (\ref{eqn3.20}) and
(\ref{eqn4.10}).

We want to use a standard map associated to a propeller in the plane. 
Let $\prop \i \R^2$ denote the same standard propeller centered at the origin 
as in Definition \ref{toro-defn2.2}.
Call $D_j$, $j=1,2,3$ its three branches. Choose a Lipschitz map $h$ 
such that $h(x)$ is the orthogonal projection of $x$ onto (the line through)
$D_j$ when $\dist(x,D_j) \leq |x|/10$, say. We shall see
that the values of $h$ in the other regions do not really
matter, because images of points of $Y$ will not fall there.
To simplify the notation later, choose $h$ smooth away from the origin,
such that $h(\lambda x)=\lambda h(x)$ for $\lambda >0$, and also
invariant under symmetries with respect to the $D_j$ and rotations
of $120^\circ$.
It is slightly unpleasant that $h$ is not $C^2$ at the origin, but
the restrictions to the three sectors where $\dist(x,D_j) \leq |x|/10$
are (they are even affine).

So far $h$ was defined on $\R^2$; we extend it to a map defined on $\R^3$ 
by setting $h(x,t)=(h(x),t)$, with the obvious notation.
This is the map associated to the set $Y_0 = \prop \times \R$ of
Definition~\ref{toro-defn2.2}. If $Z$ is any set of type $2$, we define
a mapping $h_Z$ by $h_Z = \tau \circ h \circ \tau^{-1}$, where
$\tau$ is an isometry of $\R^3$ that sends $Y_0$ to $Z$.

Return to $i\in I_2(n)$. Let $Y_i$ be the set of type 2
that was introduced in (\ref{eqn4.4}) and was already used to define $\psi^\ast_i$.
Recall from Lemma \ref{lem4.13} that $\Gamma_n \cap 5B_i$
is contained in a $C_4\varepsilon$-Lipschitz graph $G_{n,i}$
over $L_i$, the spine of $Y_i$.
To illustrate the construction, choose coordinates 
of $\R^3$ so that $x_i$ is the origin
and $L_i$ is the first axis. Parameterize $G_{n,i}$ near $B_i$ by
$x \to (x,\zeta(x))$. Let us record for future use that
\begin{equation}\label{eqn5.7}  
|\zeta(0)| \leq C \varepsilon 2^{-n}\ \hbox{ and } \
|\zeta'(x)| \leq C_4 \varepsilon,
\end{equation}
by Lemma \ref{lem4.13}. 
We already decided that $\psi_i^\ast(x,\zeta(x)) = (x,0)$
on $\Gamma_n$, and we want to extend this. So we set
\begin{equation}\label{eqn5.8}  
\eta_i(x,y) = (x,y-\zeta(x))
\hbox{ for } x\in \R \hbox{ and } y \in \R^2,
\end{equation}
and then
\begin{equation}\label{eqn5.9} 
\psi_i(z) = h_i(\eta_i(z))
\hbox{ for } z\in \R^3,
\end{equation}
where we set $h_i = h_{Y_i}$.
Notice that (\ref{eqn5.6}) holds because in $3 B_i$, $\eta_i$ maps points of
$\Gamma_n$ to $L_i$ just by the same formula as the orthogonal
projection $\psi_i^\ast$, and then $h_i$ leaves them alone.

This completes our definition of the $\psi_i$, $g_n$, and $f_n$.
Now we need to do the same sort of estimates as in the last section.
Hopefully, the similarity will help reduce the amount of work, but we
may need to be a little more careful with the analogue of Lemma
\ref{lem4.13}.

Our set $Y = Z(0,2) \cap B(0,195/100)$ is composed of three
faces. Let $F$ be one of these faces, and set $F_n = f_n(F)$.
From time to time, we need to recall which initial face
things are coming from, and in this case we add an index 
$l$, $1 \leq l \leq 3$,
as an exponent. Thus there are three faces $F^l$, and for each $n$
three sets $F_n^l$.

First let us check that points do not move too much. Observe that
for $g_n$ and in (\ref{eqn5.1})
\begin{equation}\label{eqn5.10} 
|g_{n}(z)-z| \leq  2^{-n-37}  \hbox{ for } z \in \R^3,
\end{equation}
just because $|\psi_i(x)-x| \leq 6 r_i \leq 2^{-n-37}$ in
$3B_i$ (the support of $\theta_i$), and then $g_{n}(z)$
is a convex combination of the $\psi_i(x)$. This will be useful
to make sure that our local descriptions of the $F_n \cap 5B_i$
are not upset by points coming from outer space.

Now we want to prove, at the same time and by induction on $n$,
the following collection of estimates. Set
\begin{equation}\label{eqn5.11}
\rho_n = \frac{195}{100} - \sum_{k=0}^n 2^{-k-30}. 
\end{equation}
First we want to prove that
\begin{equation}\label{eqn5.12}  
\dist(z,E) \leq C_1 \varepsilon 2^{-n}
\hbox{ for } z\in F_n \cap B(0,\rho_n).
\end{equation} 

Next consider $i\in I_1(n)$ such that $x_i \in B(0,\rho_n)$.
We will show that two of the sets $F_n^l \cap 5B_i$,
$1 \leq l \leq 3$, are empty, and for the third one there
is a $C_2\varepsilon$-Lipschitz graph $A_{n,i}$ over $P_i$ such that
$A_{n,i}$ meets $B(x_i,C_3 \varepsilon 2^{-n})$ and
\begin{equation}\label{eqn5.13} 
F_n^l \cap 5B_i = A_{n,i} \cap 5B_i.
\end{equation}

Similarly, if $i\in I_2(n)$, $x_i \in B(0,\rho_n)$, and
$1 \leq l \leq 3$, we can find a face $V_i^l$ of $Y_i$,
a (closed) Lipschitz domain $S_i^l$ in the plane $P_i^l$
that contains $V_i^l$, with
\begin{equation}\label{eqn5.14}  
D_{x_i,10r_i}(V_i^l,S_i^l) \leq C_5\varepsilon 
\end{equation}
(we leave out the name $C_4$, to avoid confusion with 
Lemma \ref{lem4.13}), and a $C_6\varepsilon$-Lipschitz graph 
$T_{n,i}^l$ over $S_i^l$, such that $T_{n,i}^l$ meets 
$B(x_i,C_7\varepsilon 2^{-n})$ and
\begin{equation}\label{eqn5.15} 
F_n^l \cap 5B_i = T_{n,i}^l \cap 5B_i.
\end{equation}
Moreover, the three faces $V_i^l$, $1 \leq l \leq 3$, are different,
and the three $T_{n,i}^l$ are bounded by $\Gamma_n \cap 5B_i$
(which is a Lipschitz graph, by Lemma \ref{lem4.13} and (\ref{eqn4.20})).

Let us first check that all the properties above hold for $n=0$. 
Here the $F^l_n$ are the three faces of $Y = Z(0,2) \cap B(0,195/100)$,
so (\ref{eqn5.12}) follows directly from the fact that 
$D_{0,2}(Z(0,2),E) \leq \varepsilon$.
The Lipschitz description in (\ref{eqn5.13}) for $i\in I_1(0)$
follows from (\ref{eqn5.4}), Lemma \ref{toro-lem1}, and the same 
fact. Similarly, for $i\in I_2(0)$, the same fact and Lemma \ref{toro-lem1} 
force $Z(0,2)$ to coincide with a set of type 2 in $5B_i$, with a spine that
gets very close to $x_i$ (see Lemma \ref{lem17}), and the description in 
(\ref{eqn5.13})-(\ref{eqn5.14}) follows.
So now we assume that $n\geq 0$, and that the properties hold for $m\leq n$.

Let $z\in F_n^l \cap B(0,\rho_n)$ be given.
By (\ref{eqn5.12}) and (\ref{eqn3.8}), $z$ lies $C_1 \varepsilon 2^{-n}$ close 
to some $B_i$, $i\in I_1(n)$ or to some $\frac{7}{4}B_i$, $i\in I_2(n)$.
Choose such an index $i$, and call it $j(z)$.
Also call $I_z$ the set of indices
$i\in I(n)$ such that $z\in 3B_i$. We already know
from (\ref{eqn5.12}) and (\ref{eqn3.23}) that $I_z \i I_2(n) \cup I_1(n)$.

Let $i\in I_z$ be given, and let us assume that
$x_i\in B(0,\rho_n)$, so that we can use the descriptions
above. Note that in this case $F_n^l \cap 5B_i$ is not empty
(because it contains $z$), so we have (\ref{eqn5.13}) or (\ref{eqn5.15}).
Since the plane $P_i$ (when (\ref{eqn5.13}) holds), or $P_i=P_i^l$
(when (\ref{eqn5.15}) holds) is reasonably well determined by
$A_{n,i}\cap 5B_i\cap B(z,2^{-n-80})$ or
$T_{n,i}\cap 5B_i\cap B(z,2^{-n-80})$, we have that
\begin{equation}\label{eqn5.16}  
D_{z,2^{-n}}(P_i,P_j) \leq C \varepsilon
\ \hbox{ for } i,j\in I_z.
\end{equation}
With this proof, $C$ depends on $C_3$ and $C_6$, but since we also
know that $P_i$ comes from $Z(x_i,100r_i)$ as in (\ref{eqn5.4}) or $Y_i$ as in
(\ref{eqn4.4}), and similarly for $P_j$, we can easily improve this first
estimate and get (\ref{eqn5.16}) with a geometric constant (i.e., that does not
depend on our future choices of constants $C_k$).

Then let us look more carefully at the $\psi_j$ and check that
\begin{equation}\label{eqn5.17}  
|\psi_i(z)-\psi_j(z)| \leq C \varepsilon 2^{-n}
\ \hbox{ for } i,j\in I_z.
\end{equation}
The main point will be that for $i\in I_z$,
\begin{equation}\label{eqn5.18}  
|\psi_i(z)-\pi_i(z)| \leq C \varepsilon 2^{-n},
\end{equation}
where $\pi_i$ denotes the orthogonal projection onto $P_i$.
When $i \in I_1(n)$, we even have that $\psi_i(z)=\pi_i(z)$,
by (\ref{eqn5.5}). When $i \in I_2(n)$, and since $z\in 3B_i$,
(\ref{eqn5.15}) says that $z\in T_{n,i}^l \cap 5B_i$.

Also, $T_{n,i}^l \cap 5B_i$ is a piece of Lipschitz graph over $P_i$,
bounded by $\Gamma_n \cap 5B_i \,$; its image under $\eta_i$ is still 
a Lipschitz graph, but it is now bounded by an arc of $\eta_i(\Gamma_n)$, 
which is an arc over $L_i$. See (\ref{eqn5.8}) and the few lines above it. 
Since the Lipschitz constant for $\eta_i(T_{n,i}^l \cap 5B_i)$ is very 
small (by (\ref{eqn5.7}) and (\ref{eqn5.8}) in particular), and
this set is bounded by an arc of $L$, it is contained in the small
sector around $V_i^l$ (the face that is contained in $P_i^l$),
where $h_i=h_{Y_i}$ coincides with the projection $\pi_i^l$ on $P_i^l$.
Thus (\ref{eqn5.9}) says that
\begin{equation}\label{eqn5.19}  
\psi_i(z) = \pi_i^l(\eta_i(z))
\hbox{ for } z\in T_{n,i}^l \cap 5B_i.
\end{equation}
Now $|\eta_i(z)-z| \leq C \varepsilon 2^{-n}$ on $5B_i$ by
the definition (\ref{eqn5.8}) and the fact that
$\zeta$ is $C\varepsilon$-Lipschitz and nearly vanishes at the origin
(see (\ref{eqn5.7})). Then
(\ref{eqn5.18}) follows from (\ref{eqn5.19}). Once we have 
(\ref{eqn5.18}), (\ref{eqn5.17}) follows from (\ref{eqn5.16}).

Next we claim that
\begin{equation}\label{eqn5.20}  
\dist(g_n(z),E) \leq C \varepsilon 2^{-n},
\end{equation}
again with a geometric constant $C$. When we can find $i\in I_z$ such that
$i\in I_1(n)$, this comes directly from (\ref{eqn5.17}), (\ref{eqn5.18}), and
(\ref{eqn5.4}) (which says that $\pi_i(z)$ lies very close to $E$).
Otherwise, $i=j(z)$ (say) lies in $I_2(n)$, and (\ref{eqn5.20}) follows from
(\ref{eqn5.17}), (\ref{eqn5.19}), and the fact that $\pi_i^l(\eta_i(z))$ lies on
the face $V_i^l \i Y_i$ (by our proof of (\ref{eqn5.19})).

\medskip\noindent{\bf Proof of (\ref{eqn5.12}).}
We are ready to prove (\ref{eqn5.12}) for $n+1$, with a geometric
constant $C_1$. Indeed, let $y\in F_{n+1} \cap B(0,\rho_{n+1})$ be given.
By construction (i.e., because $F_{n+1}=f_{n+1}(F)$ and by (\ref{eqn5.1})),
$y=g_n(z)$ for some $z\in F_{n}$. Since $\rho_{n+1}=\rho_n - 2^{-n-32}$
and by (\ref{eqn5.10}), $z$ lies in $B(0,\rho_{n})$, and so do all the $x_i$
such that $i\in I_z$. Then (\ref{eqn5.20}) holds, and
$\dist(y,E) = \dist(g_n(z),E) \leq C \varepsilon 2^{-n-1}$, which proves
(\ref{eqn5.12}).

\medskip
Let us also check that 
\begin{equation}\label{eqn5.21} 
|g_n(z)-z| \leq C\varepsilon2^{-n}
\ \hbox{ when $z\in F_n \cap B(0,\rho_n)$ and 
$x_{j(z)} \in B(0,\rho_n)$.}
\end{equation}
The constant $C$ depends on $C_5$ and $C_6$, but this will be all right.
Set $i=j(z)$. First suppose  that $i\in I_1(n)$. 
By (\ref{eqn5.12}), $\dist(z,E) \leq C_1 \varepsilon 2^{-n}$;
then $\dist(z,P_i) \leq C\varepsilon2^{-n}$ too, by (\ref{eqn5.4}).
[We could  also have used the description (\ref{eqn5.13}) as in our 
second case.] If $i\in I_2(n)$, we use the induction hypothesis and find that
$z \in T^l_{n,i}\cap  5B_i$, as in (\ref{eqn5.15}), then
$\dist(z,P_i^l) \leq C\varepsilon2^{-n}$. In both cases, 
(\ref{eqn5.18}) says that $|\psi_i(z)-z| \leq C\varepsilon2^{-n}$. 
By (\ref{eqn5.17}), $|\psi_j(z)-z| \leq C\varepsilon2^{-n}$ for 
all other $j\in I_z$; (\ref{eqn5.21}) follows because $g_n(z)$ is an 
average of $\psi_j(z)$.

\medskip\noindent{\bf Proof of (\ref{eqn5.13}).}
Next we want to prove (\ref{eqn5.13}) for $n+1$, so let $i\in I_1(n+1)$
be given and assume that $x_i \in B(0,\rho_{n+1})$; we want to look
at the sets $F_{n+1}^l \cap 5B_i$.
Since $x_i \in E'_1 \i E\cap B(0,197/100)$, (\ref{eqn3.8}) says that
$x_i$ lies in $\frac{7}{4}B_j$ for some $j\in I_2(n)$ or in
$\overline B_j$ for some $j\in I_1(n)$. So let us choose $j$, so that
\begin{equation}\label{eqn5.22} 
j \in I_1(n) \hbox{ and } x_i \in \overline B_j
\ \hbox{ or }
j \in I_2(n) \hbox{ and } x_i \in \frac{7}{4}B_j.
\end{equation}
In both cases, we shall apply the induction assumption to get
a good description of $F_n \cap 5B_j$, and use it to get information on
$F_{n+1} \cap 5B_i$. First notice that
\begin{equation}\label{7.23}  
|x_j| \leq |x_i| + 2 r_j \leq \rho_{n+1} + 2^{-n-39} < \rho_n,
\end{equation}
so we can use the induction assumption to describe $F_n \cap 5B_j$.

Next check that for $l=1,2,3$,
\begin{equation}\label{eqn5.23}  
F_{n+1}^l \cap 5B_i = g_n(F_n^l \cap 5B_j) \cap 5B_i.
\end{equation}

Clearly $F_{n+1}^l \cap 5B_i$ contains $g_n(F_n^l \cap 5B_j) \cap 5B_i$.
Conversely, if $y\in F_{n+1}^l \cap 5B_i$, then $y=g_n(z)$ for some
$z\in F_n^l$; we just need to check that $z\in 5B_j$.
Let us first check that the condition in (\ref{eqn5.21}) holds. 
By (\ref{eqn5.10}), $|z-y| \leq 2^{-n-37}$; also $|y| \leq |x_i| + 5 r_i \leq 
\rho_{n+1} + 5 r_i$ because $y\in 5B_i$, so 
$|z| \leq \rho_{n+1} + 2^{-n-37} + 5 \cdot 2^{-n-60} < \rho_n$. 
Since $|x_{j(z)}-z| \leq 3 r_{j(z)} \leq 3 \cdot 2^{-n-40}$, 
we also get that $|x_{j(z)}| < \rho_n$, and (\ref{eqn5.21}) can be 
applied. So  $|y-z| = |g_n(z)-z| \leq C \varepsilon 2^{-n}$, and 
$z$ lies very close to $5B_i$. Now
$6B_i \i 5 B_j$, either because $x_i \in \overline B_j$
and $r_j = 2r_i$ (when $\in I_1(n)$), or because 
$x_i \in \frac{7}{4}B_j$ and $r_j = 2^{21}r_i$ (when $j \in I_2(n)$).
So $z\in 5B_j$, as needed for (\ref{eqn5.23}).

\medskip
We continue with the proof of (\ref{eqn5.13}) for $i\in I_1(n+1)$,
and now consider any point $z$ such that
\begin{equation}\label{eqn5.24} 
z \in F_n \cap 5B_j
\hbox{ and } g_n(z) \in 10B_i.
\end{equation}
In particular, we assume for the moment that such a point exists.
Let $l\in \{ 1,2,3\}$ be such that $z\in F_n^l$.
Observe that 
\begin{equation}\label{7.25} 
|z| \leq |x_j| + 5 r_j
\leq |x_i| + 7 r_j \leq \rho_{n+1} + 7 \cdot 2^{-n-40}
\end{equation}
by (\ref{eqn5.22}) and because $x_i \in B(0,\rho_{n+1})$.
Then $|x_k| \leq \rho_{n+1} + 10 \cdot 2^{-n-40} < \rho_n$
for every $k \in I_z$. Hence $z$ and the $x_k$, $k\in I_z$, lie in $B(0,\rho_n)$
and we can use the description near (\ref{eqn5.16})--(\ref{eqn5.21}).

Pick any $k\in I_z$, and let us start with the most delicate case a
priori when
\begin{equation}\label{eqn5.25} 
k \in I_2(n).  
\end{equation}

\smallskip \noindent{\bf First case.} We assume that 
we can find $z$ satisfying (\ref{eqn5.24}) and $k \in I_z \cap I_2(n)$.

Recall that $B_i$ is a ball of radius $r_i=2^{-n-1-60}$ centered at some 
$x_i\in E'_1 \i E \cap B(0,197/100)$; see near (\ref{eqn3.3}). Also recall 
from (\ref{eqn3.14}) and (\ref{eqn3.16}) that
\begin{equation}\label{7.26}   
\dist\Big(100B_i,\frac{7}{4}B_{i_0} \cup \bigcup_{j\in I_2(n+1)} 
\frac{3}{2} B_j
\Big) \geq 2^{-n-1-50} > 1000 r_i,
\end{equation}
Then (\ref{eqn3.7}) says that $\dist(z,E_2) \geq 1000 r_i$,
where we also use the fact that $x_i \in B(0,\rho_{n+1}) \i  B(0,195/100)$.

By (\ref{eqn5.25}), $r_k = 2^{21}r_i$, so $B_k$ is much larger than $B_i$. 
Also, $z\in 3B_k$ because $k\in I_z$, so $B(z,1000r_i) \i 4 B_k$.
Recall that $Y_k = Y(x_k,10r_k)$ is a set of type 2 such that 
(\ref{eqn4.4}) holds. In particular, its spine $L_k$ is such that
$D_{x_k,10r_k}(E_2,L_k) \leq 600\varepsilon$,
and hence
\begin{equation}\label{eqn5.26}  
\dist(z,L_k) \geq 400r_i \, .
\end{equation}

We checked just above (\ref{eqn5.25}) that $|x_k| < \rho_{n}$, so we can 
use the description of the three $F_n^m \cap 5B_k = T_{n,k}^m \cap 5B_k$
near (\ref{eqn5.15}). Notice also that $\dist(z,E) \leq C_1 \varepsilon 2^{-n}$
by (\ref{eqn5.12}), applied to $z \in F_n \cap 5B_j$ (we checked in
(\ref{7.23}) that $x_j \in B(0,\rho_n)$).
By (\ref{eqn4.4}), $\dist(z,Y_k) \leq C \varepsilon 2^{-n}$.
By (\ref{eqn5.26}) it lies far from $L_k$, so there is only one face
of $Y_k$ that gets close to $z$.
Since $z\in F_n^l \cap 5B_k$ by definition of $l$ above,
$z$ belongs to the corresponding $T_{n,k}^l$ alone, and the face of $Y_k$ that
is close to $z$ is $V_k^l$. This also says that $l$ above is unique, 
i.e., $z$ only lies on one $F^l_n$.

Notice that $|z-x_i| \leq |z-g_n(z)|+|g_n(z)-x_i| 
\leq C \varepsilon 2^{-n} + 10r_i\leq 11r_i$ because $g_n(z) \in 10 B_i$ 
and we can use (\ref{eqn5.21}). In $50B_i$ we know that $E$ is very 
close to $P_i$, by (\ref{eqn5.4}), and at the same time to the face 
$V_k^l$ of $Y_k$ that gets close to $z$ (by (\ref{eqn4.4}) and (\ref{eqn5.26})).
So 
\begin{equation}\label{7.30}  
D_{x_i,40r_i}(P_i,V_k^l) \leq C \varepsilon,
\end{equation}
and also
\begin{equation}\label{eqn5.27}  
|D\pi_k^l - D \pi_i| \leq C\varepsilon,
\end{equation}
where $\pi_k^l$ is the orthogonal projection onto the plane $P_k^l$ 
that contains $V_k^l$, and $\pi_i$ is the orthogonal projection onto 
$P_i$. In addition, $C$ in (\ref{7.30}) and (\ref{eqn5.27})
is a geometric constant (i.e., does not
depend on the $C_i$'s of the induction hypothesis).

Recall from (\ref{eqn5.15}) that $F_n^l \cap 5B_k = T_{n,k}^l \cap 5B_k$,
where $T_{n,k}^l$ is a $C_6\varepsilon$-Lipschitz graph over a set
$S_k^l \i P_k^l$ that is very close to $V_k^l$ (see (\ref{eqn5.14})). Call
$A=A_k^l$ a $C_6\varepsilon$-Lipschitz function from the whole $P_k^l$
to an orthogonal line, whose graph contains $T_{n,k}^l$.
Since $z$ is far from the spine $L_k$ by 
(\ref{eqn5.26}) and $B(z,25r_i) \i 4B_k$
(because $k\in I_z\cap I_2(n)$, hence $z\in 3B_k$),
\begin{equation}\label{eqn5.28}  
F_n^l \hbox{ and } T_{n,k}^l
\hbox{ coincide with the graph of $A$ in }
B(z,20r_i).
\end{equation}

Recall from (\ref{eqn5.19}) that
$\psi_k(w) = \pi_k^l(\eta_k(w))$ for $w\in T_{n,k}^l \cap 5B_k$
(which includes all points of $F_n^l\cap B(z,18r_i)$). The proof also
says that $\psi_k(w) = \pi_k^l(\eta_k(w))$ near $w$ if $w$ does not 
lie on the boundary $\Gamma_n$. Then
\begin{equation}\label{eqn5.29} 
|D\psi_k(w) - D\pi_k^l| \leq C \varepsilon,
\end{equation}
where again $C$ is a geometric constant, because (\ref{eqn5.7}) and
(\ref{eqn5.8}) say that $|D\eta_k-I|\leq \varepsilon$.

We shall also need information on the indices $m\in I(n)$
such that $m\in I_y$ for some $y\in F_n^l\cap B(z,20r_i)$.
When $m\in I_2(n)$, we can repeat the discussion above, 
replace the information that $z\in 3B_k$ with the slightly weaker
fact that $y\in 3B_m$, which still implies that $B(z,20r_i) \i 4B_m$, for 
instance, and get that
\begin{equation}\label{eqn5.30} 
|D\pi_m^l - D \pi_i| \leq C\varepsilon
\end{equation}
(as in (\ref{eqn5.27})), then $\psi_m(w) = \pi_m^l(\eta_m(w))$ for
$w\in F_n^l \cap 5B_m$, and
\begin{equation}\label{eqn5.31}  
|D\psi_m(w) - D\pi_m^l| \leq C \varepsilon.
\end{equation}
The index $l$ stays the same, because it is determined by the fact
that $z\in F_n^l$ and $y$ lies close to $z$, so $y\in F_n^l$. 
Also, $C$ is a geometric constant.

When $m\in I_1(n)$, we can use the fact that by (\ref{eqn5.4}),
$E$ is $100r_m$-close to $P_m$ in $50B_m$. It is also
$10\varepsilon r_k$-close to $V_k^l$ and $P_k^l$ in $B(y,50r_i)$, 
because $z\in 3B_k$, $y \in B(z,20r_i)$, and by (\ref{eqn5.26}).
Now $50B_m$ contains $B(y,45 r_i)$ because $m\in I_y$ and $r_m=r_i$,
so we get that $D_{y,40r_i}(P_m,V_k^l) \leq C \varepsilon$.
We also get that $|D\pi_m - D \pi_k^l| \leq C\varepsilon$ as in (\ref{eqn5.27}),
hence $|D\pi_m - D \pi_i| \leq C\varepsilon$ by (\ref{eqn5.27}), as in 
(\ref{eqn5.30}). In this case, we have the simpler formula $\psi_m = \pi_m$, 
and (\ref{eqn5.31}) holds as well.

We may now collect the estimates (\ref{eqn5.30}), (\ref{eqn5.31}), 
and their analogues for $m\in I_1(n)$, to get that for every 
$y\in F_n^l\cap B(z,20r_i)$,
\begin{eqnarray}\label{eqn5.32} 
|Dg_n(y)- D\pi_i|
&\leq &\sum_{m\in I(n)} \theta_m(y) |D\psi_m(y) - D\pi_i|
+ \quad \Big|\sum_{m\in I(n)} D\theta_m(y) \psi_m(y)\Big| \nonumber\\
&\leq & 
C \varepsilon + \sum_{m\in I(n)} |D\theta_m(y)| |\psi_m(y) - g_n(y)|
\leq  C \varepsilon, 
\end{eqnarray}
by (\ref{eqn5.1}), and where the last inequality comes from (\ref{eqn3.12}), 
(\ref{eqn5.17}) and the fact that
$g_n(y)$ is an average of the $\psi_m(y)$, $m\in I_y$. 

We may now return to the description of $F_n^l \cap B(z,20r_i)$
given by (\ref{eqn5.28}), and see what happens when we apply $g_n$.
If $w_1$, $w_2 \in F_n^l \cap B(z,19r_i)$, there is a curve in
$F_n^l \cap B(z,20r_i)$, with length at most
$(1+CC_6\varepsilon)|w_2-w_1|$, which goes from $w_1$ to $w_2$
(just project on $P_k^l$ and join by a line segment).
We can integrate $Dg_n$ on this curve to estimate $g_n(w_2)-g_n(w_1)$;
this yields
\begin{eqnarray}\label{eqn5.33} 
|\pi_i^\perp (g_n(w_2)) -\pi_i^\perp (g_n(w_1))| & \leq & C \varepsilon |w_2-w_1|
\nonumber \\
\hbox{ and } 
|\pi_i(g_n(w_2))-\pi_i(g_n(w_1))| & \geq & \frac{99}{100}|w_2-w_1|,
\end{eqnarray}
where $\pi_i^\perp = I - \pi_i$, and by (\ref{eqn5.32}). Again $C$ is a 
geometric constant (we seem to use $C_6$ in the definition of the
curve, but the main point is that $CC_6 \leq 10^{-3}$, say).

Already (\ref{eqn5.33}) forces $g_n(F_n^l \cap B(z,19r_i))$ to be 
contained in a Lipschitz graph $\Lambda$ with constant $C \varepsilon$ over $\pi_i$. 
Let us check that there is no hole. For $w\in F_n^l \cap B(z,19r_i)$, we can still 
apply (\ref{eqn5.21}), because 
$|w| \leq |z| + 19r_i < \rho_{n+1} + 7 \cdot 2^{-n-40} + 19r_i < \rho_n$
by (\ref{7.25}), and even $|x_{j(w)}| \leq |w| + 3 \cdot 2^{-n-40} <\rho_n$.
Then 
\begin{equation} \label{7.35} 
|g_n(w)-w| \leq C \varepsilon 2^{-n}
\hbox{ for } w\in F_n^l \cap B(z,19r_i).
\end{equation}
We do not care that $C$ here is not a geometric constant; the main 
point is that $C \varepsilon 2^{-n}$ is much smaller than $r_i$, so
we can use a a little bit of degree theory (or monodromy, or an 
inversion theorem) to deduce from the description of (\ref{eqn5.28}) 
and (\ref{7.35}) that $\pi_i(g_n(F_n \cap B(z,19r_i)))$ actually 
contains $P_i \cap B(\pi_i(g_n(z)),18r_i)$. So we
found a $C\e$-Lipschitz graph $\Lambda$ over $\pi_i$, such that,
if we set $D=B(g_n(z),18r_i)$,
\begin{equation}\label{eqn5.34} 
g_n(F_n^l \cap B(z,19r_i)) \cap D = \Lambda \cap D.
\end{equation}

Recall from (\ref{eqn5.24}) that $g_n(z)\in 10B_i \,$; thus 
$D$ contains $5B_i$ and (\ref{eqn5.34}) implies that
\begin{equation}\label{eqn5.35} 
g_n(F_n^l \cap B(z,19r_i)) \cap 5B_i
= \Lambda \cap 5B_i.
\end{equation}
Let us check that
\begin{equation}\label{7.38} 
g_n(F_n^l \cap 5B_j) \cap 5B_i\i g_n(F_n^l \cap B(z,19r_i)).
\end{equation}
Observe that $|g_n(w)-w| \leq C \varepsilon 2^{-n}$
when $w\in F_n^l \cap 5B_j$, by the proof of (\ref{7.35}) or 
by (\ref{eqn5.21}). If in addition $g_n(w) \in 5B_i$, 
then $|z-w| \leq |g_n(w)-g_n(z)|+ C\varepsilon 2^{-n} 
\leq 15 r_i + C\varepsilon 2^{-n} < 16r_i$
(because $g_n(w) \in 5B_i$ and $g_n(z)\in 10B_i$, by (\ref{eqn5.34})),
so $w \in B(z,19r_i)$. This proves (\ref{7.38}).

Now (\ref{eqn5.23}) and (\ref{7.38}) say that 
$F_{n+1}^l \cap 5B_i \i g_n(F_n^l \cap 5B_j) \cap 5B_i 
\i g_n(F_n^l \cap B(z,19r_i)) \cap 5B_i$.
Incidentally, the converse inclusion is trivial. We compare
with (\ref{eqn5.35}) and get that
\begin{equation}\label{eqn5.36} 
F_{n+1}^l \cap 5B_i = \Lambda \cap 5B_i.
\end{equation}
This is the same thing as (\ref{eqn5.13}), but we also need
to check the small additional properties that go with (\ref{eqn5.13}).

First we need to check that $\Lambda$ meets $B(x_i,C_3\varepsilon 2^{-n-1})$.
Recall that $g_n(z)$ lies in $\Lambda \cap 10B_i$ by (\ref{eqn5.24})
and (\ref{eqn5.34}); since $\Lambda$ is a $C\varepsilon$-Lipschitz 
graph over $P_i$ and the center $x_i$ lies in $P_i$, it is enough 
to check that $\dist(g_n(z),P_i) \leq C\varepsilon 2^{-n-1}$ for some 
geometric constant $C$. By (\ref{eqn5.12}) (for $n+1$), 
$\dist(g_n(z),E) \leq C \varepsilon 2^{-n-1}$.
(We can safely consider now that the constant in (\ref{eqn5.12}) is 
geometric, since we proved (\ref{eqn5.12}) already.)
Then $g_n(z)$ lies very close to $P_i$, because $g_n(z) \in 10B_i$ and
by (\ref{eqn5.4}). Thus $\Lambda$ meets $B(x_i,C_3\varepsilon 
2^{-n-1})$, as needed.

We also need to check that that two of the $F_{n+1}^l \cap 5B_i$
are empty and the third one is a graph. With our assumptions, we 
already found $l$ such that $F_{n+1}^l \cap 5B_i$ is a graph, we just
need to check that the other two $F_{n+1}^m \cap 5B_i$ are empty. 

Let $j$ be as in (\ref{eqn5.22}) and (\ref{eqn5.23});  
If $j\in I_1(n)$, then by induction assumption (relative to (\ref{eqn5.13})), 
only one $F_n^m \cap 5B_j$ is nonempty, and (\ref{eqn5.23}) says that at only 
one $F_{n+1}^l \cap 5B_i$ is nonempty, as needed.
If $j\in I_2(n)$, we know that $x_i \in \frac{7}{4}B_j$.
If $y\in F_{n+1}^m \cap 5B_i$, (\ref{eqn5.23}) says that
$y=g_n(z)$ for some $z\in F_n^m \cap 5B_j$; then $z$ satisfies
(\ref{eqn5.24}), and so $|g_n(z)-z| \leq C\varepsilon 2^{-n}$ by 
(\ref{eqn5.21})). Then $z\in 11 B_i \i 3B_j$ (recall that
$r_j = 2^{21}r_i$ because $j\in I_2(n)$), so $j$ lies in $I_z$.
Then we can do the argument above (from (\ref{eqn5.25}) down) with $k=j$. 
We get that $z\in F_n^l$ for the only $l$ such that 
$|D\pi_j^l - D \pi_i| \leq C\varepsilon$
(as in (\ref{eqn5.27})), for instance. In other words, $m=l$.
This completes proof of (\ref{eqn5.13}) and the related facts, but
only in first case when we can find $z$ as in (\ref{eqn5.24}) and 
$k\in I_z \cap I_2(n)$.

\smallskip\noindent{\bf The other case.}
First we should check that $F_{n+1} \cap B_i$ is not empty.
Let  $j$  be as in (\ref{eqn5.22}), and first assume that $j \in I_2(n)$; 
by induction assumption, we have a description of $F_{n} \cap 5B_j$
as a union of three pieces of Lipschitz graphs that implies that
every point of $Y_j \cap 4B_j$ lies $C\varepsilon 2^{-n}$-close to
some point of $F_{n}$ (here $C$ depends on $C_5$ and $C_6$, but this 
will not matter). Since $x_i \in E \cap \frac{7}{4} B_j$, it is
$10\varepsilon r_j$-close to $Y_j = Y(x_j,10r_j)$, so it is also
$C\varepsilon 2^{-n}$-close to $F_n$. Let $z\in F_n$ be such that
$|z-x_i| \leq C\varepsilon 2^{-n}$; then
$|z| \leq |x_i| + C\varepsilon 2^{-n} \leq \rho_{n+1} + C\varepsilon 
2^{-n}$, $|x_{j(z)}| \leq |z| + 3 \cdot 2^{-n-40} < \rho_n$, 
and (\ref{eqn5.21}) holds. That is, $|g_n(z)-z| \leq C\varepsilon 2^{-n}$,
hence $|g_n(z)-x_i| \leq C\varepsilon 2^{-n}$, and even though $C$ 
depends on the $C_i$'s, we get that $g_n(z) \in B_i$.

If $j \in I_1(n)$, then by (\ref{eqn5.13}) $F_{n} \cap 5B_j$ 
is a Lipschitz graph that stays close to $P_j$. Since now
$x_i \in E \cap \overline  B_j$, $x_i$ is also very close to $P_j$,
and we can find $z\in F_n$ such  that $|z-x_i| \leq C\varepsilon 2^{-n}$; 
the same argument as before shows that $g_n(z) \in F_{n+1} \cap B_i$.

So we can find $y \in F_{n+1} \cap B_i$. By (\ref{eqn5.23}),
$y = g_{n}(z)$ for some $z\in F_n \cap 5 B_j$. In particular, we can 
find $z$ such that (\ref{eqn5.24}) holds. 

Now let $z$ be any point such that (\ref{eqn5.24}) holds.
As before, (\ref{7.25}) shows that (\ref{eqn5.21}) holds, i.e.,
$|g_n(z)-z| \leq C\varepsilon 2^{-n}$, and hence $z\in 11B_j$. 

Pick $k\in I_z$. Since we are no longer in our first case,  
$k\in I_1(n)$. Let $P_k$ be as in (\ref{eqn5.4}) and the definition 
of $\psi_k$. Thus  $\psi_k=\pi_k$, where $\pi_k$ is the orthogonal projection 
onto $P_k$. Notice that $z\in 11B_i \cap 3B_k$ (because $k\in I_z$), so
(\ref{eqn5.4}) says that $E$ is $C\varepsilon r_k$-close to both $P_i$ and $P_k$ 
in $B(z,11r_k)$, and $D_{z,10r_k}(P_i,P_k) \leq C \varepsilon$.
Hence $|D\psi_k - D \pi_i|=|D\pi_k - D \pi_i| \leq C\varepsilon$
for $k\in I_z$ and, 
after averaging,
\begin{equation}\label{eqn5.37} 
|Dg_n(z) - D \pi_i| \leq C\varepsilon.
\end{equation}
See (\ref{eqn5.32}) for a little more detail about the proof.

Suppose in addition that $g_n(z) \in 6B_i$ (instead of $10B_i$ in
(\ref{eqn5.24})). Recall that $|g_n(z)-z| \leq C\varepsilon 2^{-n}$ by
(\ref{eqn5.21}), so $z$ lies very close to $6B_i$ too. Let $j$  be as in (\ref{eqn5.22}).
If $j\in I_2(n)$, then $B_j$ is much larger than $B_i$, hence
$z\in 3B_j$; this is impossible, because we assumed that $I_z$
does not meet $I_2(n)$. So $j\in I_1(n)$, $r_j=2r_i$, 
$x_i \in \overline B_j$, $6B_i \i 4B_j$, and $z$ lies in $4B_j$ or very 
close to it.

A first consequence of this is that $z$ lies in the only $F_n^l$,
$1 \leq l \leq 3$, such that $F_n^l \cap 5B_j$ is not empty.
This proves that
\begin{equation}\label{7.44}   
F_{n+1} \cap 5B_i = g_n(F_n^l \cap 5B_j) \cap 5B_i
=F_{n+1}^l \cap 5B_i
\end{equation}
for that $l$. Indeed, if $y\in F_{n+1} \cap 5B_i$,
(\ref{eqn5.23}) says that $y=g_n(z)$ for some $z\in 5B_j$, and then 
$z\in F_n^l$, so $y\in F_{n+1}^l$.

Return to $z\in F_n$ such that $g_n(z) \in 6B_i$, and consider
$w\in F_n \cap B(z,3r_i)$. As usual, (\ref{eqn5.21}) holds for $w$ because 
$|x_{j(w)}| \leq |w| + 3 \cdot 2^{-n-40} \leq |z| + 3r_i + 3 \cdot 2^{-n-40} 
\leq \rho_{n+1} + 7 \cdot 2^{-n-40} + 3r_i + 3 \cdot 2^{-n-40} < \rho_n$
by (\ref{7.25}). Then $g_n(w) \in 10 B_i$, and $w$ satisfies (\ref{eqn5.24}).

For $k\in I_w$, we know that $k\in I_1(n)$ because we are no longer 
in our first case, and (\ref{eqn5.13}) gives a Lipschitz graph $A_{n,k}$ 
such that $F_n \cap 5B_k 
= A_{n,k} \cap 5B_k$.
Then also $F_n \cap B(z,3r_i) = A_{n,k} \cap B(z,3r_i)$
(because $z\in 3B_k$, so $B(z,3r_i) \i 5B_k$). 
We may now repeat the argument near (\ref{eqn5.33}):
given $w_1$, $w_2 \in F_n \cap B(z,2r_i)$, we can find a short curve
in $F_n \cap B(z,3r_i)$ from $w_1$ to $w_2$; then we use (\ref{eqn5.37})
to prove that (\ref{eqn5.33}) holds; this proves that 
$g_n(F_n \cap B(z,3r_i))$ is contained in the graph of a Lipschitz 
function with constant $\leq C\varepsilon$. A monodromy or
degree argument then shows that $g_n(F_n \cap B(z,3r_i))$ coincides with
the whole graph in $B(g_n(z),r_i)$, as in (\ref{eqn5.34}). 

So we showed that for every $z\in F_n \cap 5B_j$ such that 
$g_n(z) \in 6B_i \,$, 
there is a Lipschitz function defined on $P_i$, with Lipschitz constant 
$\leq C\varepsilon$, whose graph $\Gamma_z$ satisfies
\begin{equation}\label{eqn5.38}  
g_n(F_n \cap B(z,3r_i)) \cap B(g_n(z),r_i) = \Gamma_z \cap B(g_n(z),r_i).
\end{equation}

Also  observe that $g_n(w)$ lies out of $B(g_n(z),r_i)$ when
$w$ lies out of $5B_j$, by the proof of (\ref{eqn5.23}), and when
$w$ lies in $5B_j\setminus B(z,3r_i)$, this time by (\ref{eqn5.21}).
Hence (\ref{eqn5.38}) also says that
\begin{equation}\label{eqn5.39} 
F_{n+1} \cap B(g_n(z),r_i) = \Gamma_z \cap B(g_n(z),r_i).
\end{equation}
So we got a description like the one in (\ref{eqn5.13}), with
with the smaller radius $r_i$, but for every
ball centered on $F_{n+1} \cap 6B_i$ (i.e., at a point $g_n(z)$,
$z$ as above).
We also know that each $\Gamma_z$ contains the point $g_n(z)$, 
and $\dist(g_n(z),E) \leq C \varepsilon 2^{-n}$ by (\ref{eqn5.20})
or (\ref{eqn5.12}). Then $\dist(g_n(z),P_i) \leq C \varepsilon 2^{-n}$,
by (\ref{eqn5.4}), and $C$ is even a geometric constant.

It is now easy to see that $F_{n+1}$ coincides with a Lipschitz graph 
$A_{n+1,i}$ in $5B_i$, as needed for (\ref{eqn5.13}). To prove that 
$A_{n+1,i}$ meets $B(x_i,C_3\varepsilon 2^{-n-1})$, we can use the point 
$z\in F_n$ such that $|z-x_i| \leq C\varepsilon 2^{-n}$, which we used 
to prove that $F_{n+1} \cap B_i$ is not empty. This point lies in 
$F_n \cap  B_i$, so $g_n(z)$ lies on the graph, and we just said that
$\dist(g_n(z),P_i) \leq C \varepsilon 2^{-n}$ for some geometric 
constant $C$. We now deduce that $A_{n+1,i}$ meets 
$B(x_i,C_3\varepsilon 2^{-n-1})$
from the fact that the Lipschitz constant for $A_{n+1,i}$ is small.

We already know from (\ref{7.44}) that 
$F_{n+1} \cap 5B_i= F_{n+1}^l \cap 5B_i$ 
for a single $l$, so we completed our proof of (\ref{eqn5.13}) by induction.

\smallskip\noindent{\bf Proof of (\ref{eqn5.14})--(\ref{eqn5.15}).}
Now we prove the description near (\ref{eqn5.14})--(\ref{eqn5.15}).
Let $i\in I_2(n+1)$ be given. This time, (\ref{eqn3.7}) says that
$x_i \in \overline B_j$ for some $j\in I_2(n)$. We still have that
\begin{equation}\label{eqn5.40} 
F_{n+1}^l \cap 5B_i  = g_n(F_n^l \cap 5B_j) \cap 5B_i,
\end{equation}
with the same proof as for (\ref{eqn5.23}). Thus we may restrict to
points $z\in 5B_j$. By (\ref{7.23}), we can apply the induction 
assumption to get a description of $F_n \cap 5B_j$ in terms of 
three graphs $T_{n,j}^l$ over domains $S_j^l \i P_j^l$.
We want to study $g_n$ on these sets.

Let us check once more that for $z\in 5B_j$ and $k\in I_z$,
$x_k$ lies in $B(0,\rho_k)$. Indeed
$|x_k| \leq |z| +  3r_k \leq |x_j| + 5r_j + 3r_k \leq |x_i| + 6r_j + 3r_k 
\leq \rho_{k+1} + 9 \cdot 2^{-n-40} < \rho_n$, by (\ref{eqn5.11}).
This will allow us to use the estimates (\ref{eqn5.16})--(\ref{eqn5.21})
if needed.

Let us pick an index $l$ and restrict to $z\in F_n^l \cap 5 B_j
= T_{n,j}^l \cap 5 B_j$. Call $P_j^l$ the corresponding
plane. Recall that $x_i \in \overline B_j$ and $r_j = 2 r_i$,
so $D_{x_i,9r_i}(Y_i,Y_j) \leq C \varepsilon$ because both sets
approximate $E$ well near $B(x_i,9r_i)$. Then there is a (unique) face $V_i^l$ 
of $Y_i$ such that
\begin{equation}\label{7.48}   
D_{x_i,8r_i}(V_i^l,V_j^l) \leq C \varepsilon \ \hbox{ and } \ 
|D\pi_j^l - D\pi_i^l| \leq C \varepsilon,
\end{equation}
where $\pi_i^l$ denotes the orthogonal projection onto
the plane $P_i^l$ that contains $V_i^l$.

If $k\in I_z$ lies in $I_2(n)$, we can do the same
argument with $i$ replaced with $k$, and we get that
$|D\pi_k^l - D\pi_j^l| \leq C \varepsilon$. There is
no ambiguity with the index $l$: we decided to look at
$F_n^l$, we know that near $z$, $T_{n,j}^l$
(or equivalently $F_n^l$) is a small
Lipschitz graph over $P_j^l$, and for each $k$, $P_k^l$
corresponds to the only face of $Y_k$ that makes a small
angle with $P_j^l$ or with the tangent planes to $T_{n,j}^l$ near $z$.

If $k\in I_z$ lies in $I_1(n)$, the function $\psi_k$ will
merely be the projection $\pi_k$ onto the plane $P_k$
that shows up in (\ref{eqn5.4}). We shall abuse notation from time 
to time and set $P_k^l = P_k$, to avoid distinguishing cases.
Since $E$ is close to both $P_k$ and $Y_j$ near $z$, 
$P_k$ is very close to one of the planes $P_j^m$. But then
$m=l$, because near $z$, $F_n^l$ is both a piece of small
Lipschitz graph over $P_j^m$ and a small Lipschitz graph over
$P_k$ (we can apply the induction hypothesis to describe 
$F_n \cap 5B_k$, because $x_k \in B(0,\rho_n)$). 
Altogether, we also have that
\begin{equation}\label{7.49} 
|D\pi_k - D\pi_i^l| \leq C \varepsilon \ \hbox{ when } 
k\in I_z \cap I_1(n)
\end{equation}

By the same computation as for (\ref{eqn5.32}), we get that
\begin{equation}\label{eqn5.41} 
|Dg_n(w)- D\pi_i^l| \leq C \varepsilon
\ \hbox{ for } z\in F_n^l \cap 5B_j.
\end{equation} 

Recall that in $5B_j$, $F_n^l$ coincides with a $C_6\e$-Lipschitz
graph $T_{n,j}^l$ over a Lipschitz domain $S_j^l \i P_j^l$.
Then let $w_1$, $w_2 \in F_n^l \cap 5B_j$ be given. We can find a
path in $F_n^l \cap 5B_j$ that goes from $w_1$ to $w_2$, with
length $\leq 2 |w_1-w_2|$. Then we can follow this path to estimate
$g_n(w_1)-g_n(w_2)$. That is, if $\gamma : [0,1] \to F_n^l$ is this
path, we write that
$g_n(w_2)-g_n(w_1) = \int_0^1 Dg_n(\gamma(t)).\gamma'(t)dt$,
observe that $\int_0^1 D\pi_i^l.\gamma'(t)dt = \pi_i^l(w_2)-\pi_i^l(w_1)$,
and then use (\ref{eqn5.41}) to show that
\begin{eqnarray}\label{eqn5.42} 
|\pi_i^{l\perp} (g_n(w_2))- \pi_i^{l\perp} (g_n(w_1))| & \leq & 
C \varepsilon |w_2-w_1|, \\
|\pi_i^l(g_n(w_2))-\pi_i^l(g_n(w_1))| & \geq & \frac{99}{100}|w_2-w_1|.\nonumber
\end{eqnarray}
In other words, $g_n(F_n^l \cap 5B_j)$ is contained
in a Lipschitz graph $G$ over $P_i$ (with constant less than
$C\varepsilon$).

We still need to show some surjectivity, and for
this we plan to use a little bit of degree theory (or at least an
index). We shall use the fact that in $5B_j$, $F_n^l$ coincides with the
graph $T_{n,j}^l$ over $S_j^l \i P_j^l$, that $S_j^l$ is a Lipschitz domain
with small constant, and that the boundary of $T_{n,j}^l$ (i.e., the
image of the boundary of $S_j^l$ by $\varphi_j^l$, the standard 
parameterization of the graph by the base) is $\Gamma_n \cap 5B_j$ 
(which indeed is the graph of a Lipschitz function with constant 
$\leq C\varepsilon$, by Section \ref{first}).

We want to construct a simple curve in $T_{n,j}^l \cap 5B_j$.
Set $D_0=B(\pi_j^l(x_j),48r_j/10) \cap P_j^l$, call $D$ the
interior of $D_0\cap S_j^l$, and call
$\gamma_0$ the boundary of $D$. This is a nice simple
curve, composed of a piece of $\pi_j^l(\Gamma_n)$ (which is small
Lipschitz graph, because $\Gamma_n$ is a small Lipschitz graph over
the spine $L_j\i P_j^l$), completed by an arc of the circle
$\partial D_0$. Then let $\gamma_1$ be the image of $\gamma_0$
by the standard parameterization $\varphi_j^l$. It is easy
to see that $\gamma_1 \i T_{n,j}^l \cap 5B_j$, and it is still
a simple Lipschitz curve.
Finally set $\gamma_2 = \pi_i^l(g_n(\gamma_1))$.
Once again,  $\gamma_2$ is a simple curve, because (\ref{eqn5.42})
says that $\pi_i^l \circ g_n$ is injective on $\gamma_{1} \i F_n^l \cap 5B_j$.
It is composed of an arc of $\pi_i^l(\Gamma_{n+1})$ (which is still a 
Lipschitz curve with small constant) and a curve in $P_i^l$ that looks a lot like 
an arc of circle and that stays out of $B(\pi_i^l(x_j),47r_j/10) \cap P_i^l$
because $\pi_i^l \circ g_n$ does not send points of $\gamma_{0}$ too 
far from where they are.

Consider the two domains in $P_i^l$ bounded by the arc of
$\pi_i^l(\Gamma_{n+1})$ mentioned above and
$\partial B(\pi_i^l(x_j),47r_j/10)$, and call them
$D_1$ and $D_2$. Also call $i(\xi)$ the index of $\gamma_2$
with respect to a point $\xi \in P_i^l$. It is easy to deduce from the
description of $\gamma_{2}$ above that $i(\xi) = 0$ on one of the 
two domains (call it $D_1$), because we can
deform $\gamma_2$ into a point outside of $D_1$.
Then $i(\xi) \neq 0$ on $D_2$, for instance because the index
changes when we cross the arc of $\pi_i^l(\Gamma_{n+1})$
near $x_j$, and does not change inside $D_2$ because we never
meet $\gamma_2$. 

Now set $h=\pi_i^l \circ g_n \circ \varphi_j^l$ on $D$.
We can deform $\gamma_0$ into a point in $\overline D$, and when
we take the image by $h$ we get a deformation of $\gamma_2$
into a point in $h(\overline D)$. If $\xi \in P_i^l$ lies out of
$h(\overline D)$, then the deformation of $\gamma_2$ never
crosses $\xi$, the index of the curve respect to $\xi$ does not
change along the deformation, and so $i(\xi)=0$. This does not happen
when $\xi \in D_2$, so $D_2 \i h(\overline D)$.

Now $D_1$ does not meet $h(D)$. Indeed suppose that $\xi \in D_1$
can be written $\xi = h(x)$ for some $x\in D$, and select any
$y\in D$ such that $h(y) \in D_2$. This is easy to do, for instance
pick $y$ near the middle of $D$.
Call $\zeta$ a path from $x$ to $y$ in $D$.
Observe that $h(\zeta) \i B(\pi_i^l(x_j),49r_j/10)$, because
$|h(x)-x| \leq C \varepsilon 2^{-n}$ in $D$ (by (\ref{eqn5.21}),
(\ref{eqn5.14}) for $B_j$, and because $P_i^l$ is very close to $P_j^l$).
Then $h(\zeta)$ crosses $\pi_i^l(\Gamma_{n+1})$,
because it goes from $D_{1}$ to $D_{2}$ and stays in 
$B(\pi_i^l(x_j),49r_j/10)$.  The intersection lies in 
$\pi_i^l \circ g_n(\Gamma_{n}) \cap B(\pi_i^l(x_j),49r_j/10) \i
\pi_i^l \circ g_n(\Gamma_{n} \cap 5B_j)$,
because $\pi_i^l \circ g_n$ barely moves points of $\Gamma_n$.
In other words, we found $w\in \zeta \i D$  and
$s \in \Gamma_{n} \cap 5B_j$ such that $h(w)=\pi_i^l \circ g_n(s)$.
Recall that $D$ is the interior of $D_0\cap S_j^l$, so
$\varphi^l_j(w)$ lies on $T_{n,j}^l \cap 5B_j \setminus\Gamma_n
= F_n^l \cap 5B_j \setminus\Gamma_n$
because $w \in D$ and $\varphi_j^l$ is injective, and we just said that
$\pi_i^l \circ g_n(\varphi_j^l(w))=\pi_i^l \circ g_n(s)$ (by definition of
$h$). This contradicts (\ref{eqn5.42}), which says that $\pi_i^l \circ g_n$
is injective on $F_n^l \cap 5B_j$.

Return to the graph $G$ and $g_n(T_{n,j}^l)$. Set 
$B=B(x_j, 46r_j/10)$. The curve $\Gamma_{n+1}$
crosses $G \cap B$ neatly (because $\Gamma_{n+1}\cap 5B_j$
is a Lipschitz graph with small constant over $L_j$, and of course
$\Gamma_{n+1}\cap B \i g_n(F_n^l \cap 5B_j) \i G$).
So $\Gamma_{n+1}$ splits $G \cap B$ into two halves.
The image by $\pi_i^l$ of the first one (call it $T_1$) is contained
in $D_1$, and the image of the other one (call it $T_2$) is contained
in $D_2$. Since $D_1$ does not meet $h(D)$ and $D_2 \i h(\overline D)$,
this means that $T_1$ does not meet $g_n \circ \varphi_i^l(D)$ and
$T_2 \i g_n \circ \varphi_i^l(\overline D)$.

The second part of the last sentence says that $T_2 \i F_{n+1}^l$. 
The first part only says that $T_1$ does not meet $g_n \circ \varphi_j^l(D)$, 
but if $y\in T_1 \cap 5B_i$ were to belong to $F_{n+1}^l$, 
it would have to lie in $g_n(F_n^l \cap 5B_j)$ by (\ref{eqn5.40}). 
But if $x\in F_n^l \cap 5B_j$ is such that $g_n(x)=y$, 
(\ref{eqn5.21}) says that $|x-y| \leq C \varepsilon 2^{-n}$,
and then $\pi_{i}^l(x) \in D$ and $y\in g_n \circ \varphi_j^l(D)$, 
a contradiction because then $\pi_{i}^l(y) \in h(D)$.

Observe that $5B_i \i B$ because $r_{i}=r_{j}/2$. Thus we showed that
$F_{n+1}^l \cap 5B_i = T_2 \cap 5B_i$, which is (\ref{eqn5.15}) with
$T_{n+1,i}^l = T_2$. We also know that $T_2$ is bounded by
an arc of $\Gamma_{n+1}$, and since the domain $D_2$
is bounded by $\pi_i^l(\Gamma_{n+1})$, we also get the desired
$C_6\e$-Lipschitz bound on the domain $S_i^l$.
Next, $T_{n+1,i}^l$ meets $B(x_i,C\varepsilon 2^{-n-1})$ trivially,
because it contains an arc of $\Gamma_{n+1}$ that meets
$B(x_i,C\varepsilon 2^{-n-1})$ by Lemma \ref{lem4.13}. 
Similarly, (\ref{eqn5.14}) holds by Lemma \ref{lem4.13}.

Finally, we should say that the only intersection of the three
$T_{n+1,i}^l \cap 5B_i$, $l=1,2,3$, is $\Gamma_{n+1} \cap 5B_i$.
This is because the $T_{n+1,i}^l$'s leave from $\Gamma_{n+1}$ in
directions that almost make $120^{\circ}$ angles.

This completes our verification of (\ref{eqn5.10})--(\ref{eqn5.15}) by induction.

\medskip
It is easy to see that the restriction of $f_n$ to each face $F^l$
is $C^1$ (and even smoother), but this may not be true of the restriction 
of $f_n$ to the whole $Y$. Let us rapidly explain why. Since we know that 
$f_n$ is $C^1$ on each face (and all the way up to the spine $L$), 
we would just need to know that the three derivatives of the three 
restrictions, at a given point of $L$, all come from a same {\it linear} 
mapping that is defined on $\R^3$. That is, we take a point $x$ 
in $f_n(L)$ and ask whether the three derivatives of $g_n$, that are defined 
on the three tangent planes at $x$ (to the faces $F_n^l$), come from a same 
linear mapping defined on $\R^3$. It would even be enough to consider each 
piece $\theta_i(x) \psi_{i}(x)$ separately (by (\ref{eqn5.1}) and linearity).
Near $x$, the $\psi_i$ for which $\theta_i \neq 0$ are given by
a formula like (\ref{eqn5.9}), and since $\eta_i(x)$ is nicely defined and
smooth on $\R^3$, we would just need to consider $h_i$ (now defined
on three faces of the $\eta_i(F_n^l)$).

Nothing much happens in the direction of $L_i$ (where tangent
maps are the identity), but in the transverse directions,
each branch is sent to a corresponding $P_i^l$ by orthogonal
projection, and so lengths of vectors are multiplied by a cosine.
Now even if we started from a $Y_n$ whose faces made $120^\circ$
angles everywhere, this property would probably be upset slightly by
composing with $\eta$, and then it could be that three vectors whose
sum is zero are sent by the tangent maps to three vectors whose sum
is not zero. So we should not expect to get a $C^1$ mapping directly.
We could probably ameliorate this with some extra work, but we don't 
think that this is worth the trouble: the fact that the $f_{n}$ are $C^1$
on each face, and coincide on the spine, is already quite good, and
probably enough for most purposes.
See the end of Section \ref{estimates} for additional comments.

\medskip
Now consider $f$, the limit of the $f_n$, and let us prove
(\ref{1.8}) and the biH\"older estimates in (\ref{eqn4}).
Set $B=B(0,19/20)$; notice that $B \i B(0,\rho_n)$ for all $n$.

The fact that $f(Y\cap B) \i E$ comes from (\ref{eqn5.12}). 
The fact that $E \cap B \i f(Y)$ comes from (\ref{eqn3.8}) 
(which says that for every $n$, $E \cap B(0,197/100)$ is covered 
by the $7B_i/4$, $i\in I(n)$), plus our description of the $F_n^l$ 
in the $5B_i$ (and in particular the fact that at least one $F_n^l$ 
meets $B_i$). Since $f$ does not move points much 
(by (\ref{eqn5.21}) or even (\ref{eqn5.10})), we get (\ref{1.8}) easily.

Let us also prove that 
\begin{equation}\label{eqn5.44} 
(1-C\varepsilon) \, |y-z|^{1+C\varepsilon}
\leq |f(y)-f(z)| \leq
(1+C\varepsilon) \, |y-z|^{1-C\varepsilon}
\ \hbox{ for } y,z \in Y \cap B.
\end{equation}
This is only a part of (\ref{eqn4}), we shall take care
of $B(0,18/20) \setminus Y$ in Section \ref{extension},
after we modify the values of $f$ away from $Y$.

We shall prove (\ref{eqn5.44}) as we did for $f^\ast$. 
We start with the analogue of (\ref{eqn4.22}),
i.e., the fact that
\begin{equation}\label{eqn5.43} 
(1-C\varepsilon) \dist(y,z)
\leq \dist(g_n(y),g_n(z))
\leq (1+C\varepsilon) \dist(y,z)
\end{equation}
when $y$, $z\in f_n(Y\cap B)$ are such that $|y-z| \leq 2^{-n-60}$.

By (\ref{eqn5.12}), $\dist(z,E) \leq C \varepsilon 2^{-n}$; 
by (\ref{eqn3.8}), we can find $i\in I(n)$ such that $z\in 2B_i$,
and then $y\in 3 B_i$. Moreover, we can use the description of the $F_n^l$
in (\ref{eqn5.12})--(\ref{eqn5.15}), because $B$ lies well inside 
$B(0,\rho_n)$.
If $i\in I_1(n)$, $y$ and $z$ lie in $A_{n,i}$, there is a short
curve $\gamma$ in $A_{n,i}$ that goes from $y$ to $z$, and we can use
the differential of $g_n$ on $\gamma$ to compute $g_n(y)-g_n(z)$.
This differential is $C\varepsilon$-close to $D\pi_i$,
by (\ref{eqn5.37}) (or even part of its proof, since here we do not even need
to compare with the projection on some other $P_{i'}$, $i'\in I(n+1)$).
This gives (\ref{eqn5.43}) in this first case.

If $i\in I_2(n)$ and $y, z$ lie in the same $F_n^l$, we can do the
same thing with an arc in $T_{n,i}^l$, using (\ref{eqn5.32}) instead of
(\ref{eqn5.37}).
Finally, if $i\in I_2(n)$ and $y, z$ lie in different faces $F_n^l$,
we introduce a third point $x\in \Gamma_n \cap 3B_i$.
Let us choose $x$ so that $|x-y|+|x-z|$ is minimal, but the
precise choice won't matter much.
Let $l$ and $m$ be such that $y\in F_n^l$ and $z\in F_n^m$.
We know that
$|g_n(x)-g_n(y) - D\pi_i^l \cdot (x-y)| \leq C\varepsilon |x-y|$
by the same trick as before (integrate on a path from $x$ to $y$).
Also, $|\pi_i^l(x)-\pi_i^l(y)-x+y| \leq  C\varepsilon |x-y|$
(because both points lie on the small Lipschitz graph $T_{n,i}^l$).
Thus $|g_n(x)-g_n(y)-x+y| \leq  C\varepsilon |x-y|$, and similarly
$|g_n(x)-g_n(z)-x+z| \leq  C\varepsilon |x-z|$. Finally,
$|g_n(z)-g_n(y)-z+y| \leq  C\varepsilon (|x-y|+|x-z|)
\leq  C\varepsilon (|y-z|)$ (by choice of $x$), and (\ref{eqn5.43}) follows.

We may now continue the discussion as we did between (\ref{eqn4.22}) and
(\ref{eqn4.27}), and we obtain (\ref{eqn5.44}) just like we proved 
(\ref{eqn4.27}).

\section{The case of a tetrahedron} \label{case}

In this section we show how to modify the construction above when
(\ref{eqn29}) holds (and hence $E_3 \cap B(0,199/200)=\{ 0 \}$).
As was hinted near the end of Section \ref{decomp}, the argument 
would work the same way if $E_3 \cap B(0,199/200)=\{ x_0 \}$ instead,
for some $x_0 \in B(0,195/100)$; we would just need to center the
balls $B_{i_0}$ at $x_0$. We shall use the same sort of formulas 
as in Sections \ref{first} and \ref{param}, except that for each 
$n\geq 0$, we have to define the $\psi_i^\ast$ and the $\psi_i$ 
differently when $i=i_0$ (i.e., when $x_i=0$). 

For $\psi_i^\ast$, we want to define approximate
projections on spines of sets of type $3$, a little like  
what we did for $h$ near (\ref{eqn5.9}).
Let $T_0$ denote a minimal set of type $3$ centered at the
origin; for instance, take $T_0$ as in Definition \ref{toro-defn2.3}.
Denote by $p_0$ a function on $\R^3$ which is positively
homogeneous of degree 1, Lipschitz (or even $C^\infty$ if we want)
on the unit sphere, and such that for each of the four branches $D_j$
of the spine of $T_0$, $p_0$ coincides with the orthogonal
projection on (the line containing) $D_j$ in the region
where $\dist(z,D_j) \leq |z|/10$. We can also arrange symmetry
with respect to the isometries that fix $T_0$, even though this is
not needed (it just allows us to have a standard map
that we use all the time). For a general set $T$ of type 3, we 
write $T = R(T_0)$ for some composition $R$ of a rotation and 
a translation, and we set $p_T = R \circ p_0 \circ R^{-1}$.

For each $n \geq 1$, Corollary \ref{cor11} says that 
$c(0,2^{-n-10}) \leq 15 \varepsilon$, so we can select
a set $Z_n$ of type 3, centered at $0$, such that
$D_{0,2^{-n-10}}(E,Z_n) \leq 15\varepsilon$.
For $n=0$, take $Z_0 = Z(0,2)$. Then set $\psi_{i_0}^\ast = p_{Z_n}$
($n$ will often stay implicit in our notation).

With this choice of $\psi_{i_0}^\ast$, we can construct 
a parameterization of (a good piece of) $E_2 \cup \{ 0 \}$ by 
the spine of $Z_0 = Z(0,2)$ (intersected with $B(0,197/100)$, say), 
as we did in Section \ref{first}. The argument is a simpler version of 
Section \ref{param}, and we shall not repeat it all. 
Notice that we never need to worry about the radii $\rho_n$ 
when the balls $B_{i_0}$ come into play, because $B_{i_0}$ is
far from $\partial B(0,197/100)$ (one less worry). Also observe 
that all the $\theta_i$, $i\neq i_0$,
vanish near $\frac{3}{2} B_{i_0}$, by (\ref{eqn3.10}) and 
(\ref{eqn3.13}). So 
\begin{equation}\label{8.1}   
g_n^\ast(x) = \psi_{i_0}^\ast(x) = p_{Z_n}(x)
\ \hbox{ for } x\in \frac{3}{2} B_{i_0}.
\end{equation}
Finally notice that the $f_n^\ast$ and $g_n^\ast$ preserve
the origin, by (\ref{8.1}).

So we get maps $f^\ast_n$, $g_n^\ast$, and sets $\Gamma_n$.
Most of Section \ref{first} stays the same, except that 
now each $\Gamma_n$ is composed of four branches $\Gamma_n^m$
which meet at the origin.

Next consider the mapping $f$ of Section \ref{param}. Now we
need to define the mappings $\psi_{i_0}$ for $n \geq 0$.
Again start from a simpler mapping $h_n$, which is 
positively homogeneous of degree 1, smooth on the unit sphere,
and such that, for each
of the six faces $V_l$ that compose $Z_n$, 
\begin{equation}\label{8.2}   
h_n(x) = \pi_l(x) 
\ \hbox{ for } x \in {\cal C}_l = \{x\in \R^3 \, ; \, 
\dist(x,V_l) \leq |x|/100 \},
\end{equation}
where $\pi_l$ denotes the orthogonal projection on the plane that 
contains $V_l$. The simplest is to construct such a map for $T_0$
and get the other ones by conjugation with isometries as before,
but this is not needed (provided that we keep some uniform estimates
on the derivatives). As before, the values of $h_n$ out of the six
blade-like sectors ${\cal C}_l$ of (\ref{8.2}) do not matter too much, 
because the $F_n^l \cap 5B_{i_0}$ will stay in the ${\cal C}_l$.

Denote by $L^{m}$, $1 \leq m \leq 4$, the four half lines that
compose the spine of $Z_n$. We want to precompose with a mapping 
$\eta$ like the one in (\ref{eqn5.8}), whose effect near the origin 
is to map the four branches $\Gamma_n^m$ to the corresponding half 
lines $L^{m}$.

For each $m$, we have a Lipschitz parameterization
of $\Gamma_n^m \cap 5B_{i_0}$ by a segment of $L^{m}$, which we
can write $x \to (x,\zeta^m(x))$ (this time, with $x\geq 0$)
in appropriate coordinates (depending on $m$).
We want $\eta$ to be Lipschitz, with
\begin{equation}\label{eqn6.1}   
\eta(0) = 0 \ \hbox{ and }
|D\eta - I| \leq C \varepsilon
\end{equation}
and such that
\begin{equation}\label{eqn6.2}  
\eta(x,y) = (x,y-\zeta^m(x)),
\end{equation}
in the little cone where $\dist((x,y),L^m) \leq 10^{-1}|(x,y)|$.
Of course we want such a description simultaneously for all $m$,
but there is no problem because the little cones near the $L^m$ 
are far from each other. We can even require that $\eta(x,y)=(x,y)$ out of 
the four cones. Then we set $\psi_{i_0} = h_n \circ \eta$, as before.

Let us check that then $f_n$ coincides with $f_n^\ast$ on $\Gamma$.
First, $f_0(x) = f_0^\ast(x)=x$, so $n=0$ is all right.
If we already know that $f_n=f_n^\ast$ on $\Gamma$ and want
to check that the same thing holds for $f_{n+1}$, we just need to see
that $\psi_{i}=\psi^\ast_{i}$ on $\Gamma_n \cap 5B_i$ for $i\in I(n)$. 
Indeed, since we use the same partition of unity for $g_n$ and $g^\ast_n$,
$g_n$ will coincide with $g^\ast_n$ on $\Gamma_n$. See (\ref{eqn4.3}) and 
(\ref{eqn5.1}). We checked this in Section \ref{param} for $i \neq i_0$, 
so it is enough to see that $\psi_{i_0}=\psi^\ast_{i_0}$ on $\Gamma_n$.

For $z\in \Gamma_n \cap 3B_{i_0}$, (\ref{eqn6.2}) says that 
$\eta(z) = \pi_m(z)$, where $1 \leq m \leq 4$ is such that 
$z\in \Gamma_n^m$ and $\pi_m$ denotes the projection of $z$ on 
the half line $L^m$. So $\eta(z) = p_{Z_n}(z) = \psi^\ast_{i_0}(z)$,
because $p_{Z_n}$ coincides with $\pi_m$ in a small
cone around $L^m$; see the definitions above (\ref{8.1}).
Then $\psi_{i_0}(z) = h_n(\pi_m(z))=\psi^\ast_{i_0}(z)$, as needed,
because $h_n(w)=w$ on $L^m$, by (\ref{8.2}).

\medskip
We need to revise slightly our description of the $F_n^l \cap 5B_i$.
The first difference is that now $Y$ is composed of $6$ faces, so
$l$ should vary between $1$ and $6$. Here is what we need to prove now.

When $i\in I_1$, not much changes; there is still only one
$l$ for which $F_n^l$ meets $5B_i$, and for this one we have the same
description as in (\ref{eqn5.13}).

When $i\in I_2$, only three of the $F_n^l$ meets $5B_i$
(and their names depend on which of the four branches of $E_2$
contains the center $x_i$). And for these three $F_n^l$
we have the same description with Lipschitz graphs as in
(\ref{eqn5.14})--(\ref{eqn5.15}).

Then we need a description of $F_n^l \cap 5B_{i_0}$. This time, the
six faces $F_n^l$ meet $5B_{i_0}$, and the intersections are
$C\varepsilon$-Lipschitz graphs $T^{l}$ over sets $S^l$ (as before),
where $S^l$ is a Lipschitz domain in the plane $P^l$
that contains the corresponding face $V^l$ of $Z_n$, and $S^l$ is very 
close to $V^l$ (as in (\ref{eqn5.14})). In addition, the boundary 
of the graph $T^l$ is the intersection of $5B_{i_0}$ with the two curves 
of $\Gamma_n$ that bound the face $F^l$.

We need a proof by induction of these descriptions. For 
$i\in I_1 \cup I_2$, not much needs to be changed. Recall from 
(\ref{eqn3.13}) that $100 B_i$ does not meet $\frac{3}{2} B_{i_0}$
when $i \neq i_0$. We repeat the proof of Section \ref{param},
and, whenever $B_{i_0}$ shows up, we observe that we only need the 
description of the $F_n^l \cap 5B_{i_0}$ (at order $n$) away from
$B_{i_0}$. We can use this description just as the description
of the $F_n^l \cap 5B_{j}$, $j\in I_2$ (notice that in a small ball
away from $B_{i_0}$, our description of $F_n \cap 5B_{i_0}$ is
of the same type as for the $F_n \cap 5B_{j}$, $j\neq i_0$).

For $i=i_0 \in I(n+1)$, we use the description of the 
$F_n^l \cap 5B_{i_0}$ at order $n$, and check that
$g_n$ does not destroy things. Away from the origin (say,
out of $\frac{1}{2}B_i$), the proof is the same as before, 
because the maps and the description of $F_n$ are of the same type
as in Section \ref{param}.
Near $\frac{1}{2}B_i$, we can use the fact that by (\ref{eqn3.10}) and 
(\ref{eqn3.13}) only $\theta_{i_0}$ in $\frac{3}{2}B_i$,
is nonzero, so $g_n$ coincides with $\psi_{i_0}$ near $B_{i_0}$.
We use the fact that by induction assumption, the various faces 
$F_n^l$ are bounded by the arcs $\Gamma_n^m$, their images by $\eta$ 
are still reasonably small Lipschitz graphs and lie in the sectors 
where $h_n$ is a projection, so their images are locally Lipschitz 
graphs with small constants. Then composing with $h_n$ makes things 
even better, and in fact $g_n(F^l_n \cap \frac{3}{2}B_{i_0})$ 
is contained in a face of $Z_n$.

We also need to make sure that the six $F_n^l \cap 5B_{i}$, 
with $i=i_0$, are equal to the full graph $T_n^l$ (the argument
for $i \neq i_0$ can stay the same). This is again done with
an index argument like the one below (\ref{eqn5.42}); there are small 
differences, because $S_j^l$ and $D$ now look like triangular sectors, 
and the piece of $\pi_j^l(\Gamma_n)$ that bounds $D$ now comes from 
two arcs $\Gamma_n^m$. This does not upset the argument.

\medskip  
Now we need to check the analogue of (\ref{eqn5.44}), and we start with
the analogue of (\ref{eqn5.43}). We proceed essentially as for (\ref{eqn5.43}).
By (\ref{eqn3.8}), we can find $i\in I(n)$ such that $y$ and $z$ lie in $3B_i$.
We can use the description of the $F^l_n$, $1\le l\le 6$
(i.e., the analogue of (\ref{eqn5.12})--(\ref{eqn5.15})),
because $B$ lies well inside $B(0,\rho_n)$. The cases where 
$i\in I_1(n)\cup I_2(n)$ are dealt with exactly as in the proof of
(\ref{eqn5.43}), so we can assume that $i=i_0$.

If $y$ and $z$ lie in the same face $F^l_n$, we can proceed as in the proof of 
(\ref{eqn5.44}). So let us assume that $y\in F^l_n$ and $z\in F^m_n$,
with $l \neq m$. We claim that we can find 
$x\in F^l_n \cap F^m_n \cap 4B_{i_0}$, with $|x-y|+|x-z| \leq 100 |y-z|$
(notice that the first condition forces $x\in \Gamma_n$, and even 
$x=0$ if the faces $F^l_n$ and $F^m_n$ are not adjacent).

We can always try $x=0$, which works when $|y|+|z| \leq 100 |y-z|$. 
So let us assume that $|y-z| < (|y|+|z|)/100$. 
Set $y'=\pi^l(y)$, where $\pi^l$ denotes the projection onto the plane 
$P_l$ that contains $V^l$. Notice that $|y'-y| \leq C \varepsilon |y|$ 
because $F^l_n \cap 5B_{i_0} \i T^l$,
$T^l$ is a contained in a $C\varepsilon$-Lipschitz graph the plane $P^l$, 
and $T^l$ contains the origin. By construction, $y'$ lies is the Lipschitz 
domain $S^l \i P^l$. This domain is the projection of $F^l_n \cap 5B_{i_0}$,
and it is bounded by two arcs $\pi^l(\Gamma_n^s)$, where the two
values of $s$ correspond to the two half lines $L^s$ that bound $V^l$.
These two  arcs in $P^l$ are $C\varepsilon$-Lipschitz over $L^s$,
because the $\Gamma_n^s$ (in space) are $C\varepsilon$-Lipschitz over $L^s$.
Then $\dist(y',V^l) \leq C \varepsilon |y'|$. Pick $y'' \in V^l$, with
$|y''-y'| \leq C \varepsilon |y'|$. Then $|y''-y| \leq C \varepsilon |y|$
(because $|y'-y| \leq C \varepsilon |y|$).
Similarly, we can find $z''\in V^m$ such that $|z''-z| \leq C \varepsilon |z|$.
Next 
\begin{eqnarray} \label{8.5} 
|y''-z''| &\leq& |y-z| + C \varepsilon (|y|+|z|)
\leq (10^{-2}+ C \varepsilon) (|y|+|z|)
\\
&\leq& (|y|+|z|)/95\leq (|y''|+|z''|)/90. \nonumber
\end{eqnarray}
This is impossible if $y''$ and $z''$ lie in faces of $Z_n$ that
are not adjacent, so $V^l$ and $V^m$ are adjacent. Denote by 
$L^s$ the half line $V^l \cap V^m$. Notice that 
$|y''-z''| \geq \dist(y'',P^m) \geq \dist(y'',L^s)/2$ because
$V^l$ and $V^m$ make $120^{\circ}$ angles. That is,
$\dist(y'',L^s) \leq 2|y''-z''| \leq (|y|+|z|)/45$, by (\ref{8.5}). 
Then $\dist(y,L^s) \leq (|y|+|z|)/44$, 
and $\dist(y,\Gamma_n^s) \leq (|y|+|z|)/43$ because $\Gamma_n^s$ is a 
$C\varepsilon$-Lipschitz graph over $L^s$
that goes through the origin (recall that $0$ is the only point of 
type 3). Notice also that $|y|$ and $|z|$ are both pretty close to
$(|y|+|z|)/2$, because $|y-z| < (|y|+|z|)/100$. Then we can use the
description of $F^l_n = T^l$ and $F^m_n = T^m$ near $y$ and $z$, as 
Lipschitz graphs over $P^l$ and $P^m$  with a common Lipschitz boundary
$\Gamma_n^s$, to find a point $x\in \Gamma_n^s$ such that 
$|x-y|+|x-z| \leq 3|y-z|$. This proves the claim.

Once we have $x$, we can use it as in Section \ref{param}.
We use a short path $\gamma$ in $T^l$ from $x$ to $y$, obtained
for instance by lifting a line segment from $\pi^l(x)$ to $\pi^l(y)$,
to show that 
\begin{equation}  
|y-x-(\pi^{l}(y)-\pi^{l}(x))| \le C\e|y-x|,
\end{equation}
just because $x$ and $y$ lie on the small Lipschitz graph $T^l$, and
\begin{equation}  
|g_n(y)-g_n(x) - D\pi^l(y)\cdot (y-x)| 
= \Big|\int_{\gamma} [Dg_n - D\pi^l] \cdot \gamma' \Big|
\le C\e|y-x|,
\end{equation}
because $|Dg_n - D\pi^l| \leq C \varepsilon$ on $4B_{i_0}\cap F^l_n$.
Thus $|g_n(y)-g_n(x) -(y-x)| \leq C \varepsilon |y-x|$.
Similarly, $|g_n(z)-g_n(x) -(z-x)| \leq C \varepsilon |z-x|$,
and 
\begin{equation}  
|g_n(z)-g_n(y) -(z-y)| \leq C \varepsilon (|y-x|+|z-x|)
\leq 100C \varepsilon (|y-x|),
\end{equation}
by  definition of $x$. 
This proves (\ref{eqn5.43}); (\ref{eqn5.44}) follows as before.

\section{Extension to $\R^3$} \label{extension}

The mapping $f$ that was constructed in the previous sections is not
the final one for Theorem~\ref{T1.1} and Remark 1.1. 
So far, we focused on the restriction of $f$ 
to $Z(0,2)$, but in this section we need to modify $f$ away from $Z(0,2)$,
in particular to make sure that the result is injective. 
Notationally it will be more convenient to leave $f$ and the $f_n$ as 
they are, and denote by $\widetilde f$ and $\widetilde f_n$ the 
required modifications. Thus the true mapping for 
Theorem \ref{T1.1} and Remark 1.1 
is the mapping $\widetilde f$ constructed in this section.

We want to  define $\widetilde f_n$ and $\widetilde f$ near 
$B(0,18/10) \i \R^3$; we shall use the same formula as before, i.e.,
we start from $\widetilde f_0(x) = x$ and set 
\begin{equation}\label{9.1}   
\widetilde f_{n+1} = \widetilde g_{n} \circ \widetilde f_{n} 
\ \hbox{ for } n\geq 0, \hbox{ with }
\widetilde g_n(x) = \sum_{i\in I(n)} \theta_i(x) \widetilde\psi_i(x),
\end{equation}
with the same functions $\theta_i$ as before, but slightly different
$\widetilde\psi_i$. This time we shall need all the $\theta_i$, 
even with $i\in I_0$. 

Since we want $\widetilde f_n$ to coincide with $f_n$ on $Y$, 
we require that
\begin{equation}\label{eqn7.1} 
\widetilde\psi_i(x) = \psi_i
\ \hbox{ on }  3 B_i \cap F_n,
\end{equation}
where $F_n$ is the union of the three or six faces $F^l_n = f_n(F^l)$, 
and the $F^l$ are the faces of $Z(0,2) \cap B(0,195/100)$. 

When $i\in I_0(n)$, the simplest is to set $\psi_i(x)=x$;
the requirement (\ref{eqn7.1}) is empty because $3B_i$ does not meet 
$F_n$, by (\ref{eqn3.23}) and since the $F_n^l$ stay close to $E$ 
(by (\ref{eqn5.12})).

Next take $i \in I_1(n)$. Recall from (\ref{eqn5.13}), or its analogue 
when $Z(0,2)$ is of type 3, that only one of the $F_n^l$ meets $5B_i$, 
and that for this one $F_n^l \cap 5B_i=A_{n,i} \cap 5B_i$ for some 
$C_2\e$-Lipschitz graph $A_{n,i}$ over $P_i$. Call $a_{n,i} : P_i \to P_i^\perp$ 
the corresponding Lipschitz function, and set
$\widetilde\psi_i(x,y) = (x,y-a_{n,i}(x))$ for $(x,y)\in P_i \times P_i^\perp$
(we identify $\R^3$ with $P_i \times P_i^\perp$). Here the requirement
(\ref{eqn7.1}) holds by definition of $a_{n,i}$ and because $\psi_i$ is the
orthogonal projection onto $P_i$.

Now consider $i \in I_2(n)$. This time the $F_n^l \cap 5B_i$ coincide
with three half-Lipschitz graphs $T^l_{n,i}$ that meet along
$\Gamma_n$ (see near (\ref{eqn5.15})), and on these half-Lipschitz graphs
we need to take $\widetilde\psi_i = \psi_i = h_i \circ \eta_i$, 
as in (\ref{eqn5.9}). Let us check that we can extend this mapping to 
$4B_i$, in such a way that $\widetilde\psi_i$ is $(1+C \varepsilon)$-bilipschitz
on $4B_i$, and 
\begin{equation}\label{eqn7.2}  
|\widetilde\psi_i(x)-x| \leq C \varepsilon 2^{-n}
\hbox{ on } 4B_i \ \hbox{ and } \ 
|D\widetilde\psi_i(x) - I| \leq C \varepsilon
\hbox{ on } 4B_i\setminus \Gamma_n.
\end{equation} 
Recall from (\ref{eqn5.8}) that $\eta_i$, whose role is to move 
$\Gamma_n$ to the spine $L_i$ of $Y_i$, is a $(1+C \varepsilon)$-bilipschitz
mapping that satisfies (\ref{eqn7.2}). It is also smooth, because $\Gamma_n$ 
is smooth. 
Thus it is enough to find an extension $\widetilde h_i$ to 
$\frac{9}{2}B_i$, of the restriction of $h_i$ to $\eta_i(F_n \cap 5B_i)$, 
with the required properties. Recall that near each 
$\eta_i(F_n^l \cap 5B_i)$, $h_i$ coincides with the orthogonal
projection onto the plane $P^l$ that contains the face $V^l$ of $Y_i$
that lies close to $\eta_i(F_n^l \cap 5B_i)$. Set ${\cal C}^l = 
\{ z \in \R^3 \, ; \, \dist(z,V^l) \leq \frac{1}{10} |z| \}$;
we claim that we can even choose $\widetilde h_i$ so that 
$\widetilde h_i(w)=w$ when $w$ lies out of the three ${\cal C}^l$.
To define $\widetilde h_i$ on ${\cal C}^l$, pick coordinates such
that $V^l = \{(x,y,z) \, ; \, x \geq 0 \hbox{ and } y=0 \}$; 
notice that $\eta_i(F_n^l \cap 5B_i)$ is contained in the graph
of some $C\varepsilon$-Lipschitz function $F : \R^2 \to \R$ such that
$F(0,y) = 0$ for $y \in \R$ (because $\eta_i$ sends $\Gamma_n$ to 
$L_i$), and $h_i(x,y,z) = (x,y)$ on that set.
We pick a smooth function $\varphi$ on $\R$ such that $\varphi(t)=1$ for
$|t|\leq 1/20$, $\varphi(t)=0$ for $|t| \geq 1/10$, and 
$|\varphi'(t)| \leq 30$ everywhere, and set
\begin{equation}\label{9.4}  
\widetilde h_i(x,y,z) = (x,y,z) - \varphi\big(y/\sqrt{x^2+y^2}\big) F(x,y)
\ \hbox{ on } {\cal C}^l.
\end{equation}
It is easy to see that $\widetilde h_i(x,y,z)$ has all the required 
properties and that $\widetilde \psi_i=\widetilde h_i\circ\eta_i$ satisfies
(\ref{eqn7.2}). [Notice that $|F(x,y)| \leq C \varepsilon \sqrt{x^2+y^2}$
because $F(0,y) = 0$, that $|y/\sqrt{x^2+y^2}| \leq C \varepsilon < 1/20$ 
on $\eta_i(F_n^l \cap 5B_i)$, and that the $(1+C \varepsilon)$-bilipschitzness
follows from the second part of (\ref{eqn7.2}).]

The definition of $\widetilde \psi_i$ for $i=i_0$ when (\ref{eqn29}) holds 
and $Z(0,2)$ is of type 3 is similar. We have a good description of 
$F_n \cap 5B_i$, this time as a union of six Lipschitz faces $T_n^l$ 
bounded by four arcs of $\Gamma_n$, which themselves meet nicely at the origin, 
as described in the last section. 
Again, we know $\widetilde \psi_i$ on these six faces, and we claim 
that there is a $(1+C \varepsilon)$-bilipschitz extension to $4B_i$
such that (\ref{eqn7.2}) holds. Call $\Gamma_n^s$, $1 \leq s \leq 4$, 
the four arcs of $\Gamma_n$, and denote by $L^s$ the corresponding
branch of the spine of $Z_n$; thus the spine of $Z_n$ is the union of
the four $L^s$. Let $\eta$ be as in (\ref{eqn6.1}) and (\ref{eqn6.2});
thus $\eta$ is $(1+C \varepsilon)$-bilipschitz and satisfies the analogue 
of (\ref{eqn7.2}).  
Notice that $\eta$ sends $\Gamma_n^s \cap 5B_i$ to $L^s$, by (\ref{eqn6.2}).
By (\ref{eqn6.1}), it sends $T_n^l$ to a subset of a $C\varepsilon$-Lipschitz 
graph over the plane $P^l$ that contains the corresponding face $V^l$ of $Z^n$, 
and this graph goes through the two $L^s$ that bound $V^l$ (because 
$\eta$ sends the $\Gamma_n^s \cap 5B_i$ to the $L^s$). That is,
we have the same sort of description of $\eta(F_n \cap 5B_i)$ as for
$F_n \cap 5B_i$ itself, except that now $\Gamma_n$ is replaced with the
spine of $Z_n$.

We shall take $\widetilde \psi_i = \psi_i^{\sharp} \circ \eta$
for some $\psi_i^{\sharp} \,$. 
Since $\psi_i = h_n \circ \eta$ (see below (\ref{eqn6.2})),
(\ref{eqn7.1}) will hold if $\psi_i^{\sharp} = h_n$ on
$\eta(3B_i \cap  F_n)$. By (\ref{8.2}) and our description
of $\eta(F_n \cap 5B_i)$, $h_n$ coincides with the orthogonal
projection $\pi^l$ onto $P^l$, on the set $\eta(F_n^l \cap 5B_i)$. 
So  we just have to find an extension $\psi_i^{\sharp}$ to $\frac{9}{2}B_i$, 
say, of the map defined on the union of the six Lipschitz graphs 
(that contain the) $\eta(F_n^l \cap 5B_i)$, that is equal to 
$\pi^l$ on $\eta(F_n^l \cap 5B_i)$.
We can even require that $\psi_i^{\sharp}(z)=z$ out of the sectors 
${\cal C}_l$ of (\ref{8.2}), while in ${\cal C}_l$ we use a formula 
like (\ref{9.4}). This gives a bilipschitz $\psi_i^{\sharp}$ that satisfies 
(\ref{eqn7.2}), and then $\widetilde \psi_i = \psi_i^{\sharp} \circ \eta$
also satisfies (\ref{eqn7.2}) because $\eta$ does. 

This completes our definition of $\widetilde \psi_i$ in all cases.
So we also have a full definition of $\widetilde f_n$ and 
$\widetilde f$. Let us check that $\widetilde f$ satisfies the 
conclusions of Remark 1.1. 

First we have that $|\widetilde g_n(x)-x| \leq C \varepsilon$ 
everywhere, by (\ref{9.1}) and (\ref{eqn7.2}). Then 
$|\widetilde f(x)-x| \leq C \varepsilon$, again by (\ref{9.1}).
This proves (\ref{1.6}). The second half of (\ref{1.7}) follows 
immediately; the first part of (\ref{1.7}) follows from (\ref{1.6})
and a little bit of degree theory. That is, for $x\in B(0,17/10)$,
the mapping $F_1: \partial B(0,1) \to \partial B(0,1)$ defined by
$\displaystyle F_1(\xi) = 
\frac{\widetilde f(18\xi/10)-x}{|\widetilde f(18\xi/10)-x|}$, 
is homotopic to the identity (by (\ref{1.6})), so it has degree 1, 
and if $x$ did not lie in $\widetilde f(B(18/10))$, 
the mappings $F_t$, $0 \leq t \leq 1$, defined by 
$\displaystyle F_t(\xi) = 
\frac{\widetilde f(18t\xi/10)-x}{|\widetilde f(18t\xi/10)-x|}$,
would provide a homotopy from $F_1$ to a constant.

We already proved (\ref{1.8}) a little above (\ref{eqn5.44}),
so we just need to check (\ref{eqn4}), i.e., that 
\begin{equation}\label{eqn7.3}  
(1-C\varepsilon) \, |y-z|^{1+C\varepsilon}
\leq |\widetilde f(y)-\widetilde f(z)| \leq
(1+C\varepsilon) \, |y-z|^{1-C\varepsilon}
\end{equation}
for $y,z\in B(0,19/10)$. As before, it is enough to prove the 
analogue of (\ref{eqn4.22}) or (\ref{eqn5.43}), i.e., that
\begin{equation}\label{9.6} 
(1-C\varepsilon) \dist(y,z)
\leq \dist(\widetilde g_n(y),\widetilde g_n(z))
\leq (1+C\varepsilon) \dist(y,z)
\end{equation}
when $y$, $z\in B(0,185/100)$ are such that $|y-z| \leq 2^{-n-80}$.
By (\ref{9.1}), it is enough to check that for $j\in I(n)$, 
\begin{equation}\label{9.7} 
|\widetilde \psi_j(y)-\widetilde \psi_j(z) - (y-z)| 
\leq C\varepsilon |z-y|
\ \hbox{ for } y,z \in 4B_j.
\end{equation}
This last follows from (\ref{eqn7.2}).

We finally completed our verification of Theorem \ref{T1.1}.

\section{$C^1$ estimates when $E$ is asymptotically flatter}
\label{estimates}

All the estimates above were done with $\varepsilon$ fixed, but
we can do better if our basic assumption (\ref{eqn1}) improves when 
$r$ tends to $0$.
In this section, we assume that there is a sufficiently small nondecreasing 
positive function $\varepsilon(r)$ such that for each $x\in E\cap B(0,2)$ 
and $r>0$ such that $B(x,r)\i B(0,2)$, there is a minimal cone $Z(x,r)$ 
that contains $x$ such that
\begin{equation}\label{eqn8.1}  
D_{x,r}(E,Z(x,r)) \leq \varepsilon(r).
\end{equation}
In particular, we still have (\ref{eqn1}) with $\varepsilon = \varepsilon(2)$,
and we assume that $\varepsilon(2)$ is so small that we can apply the
construction above. We want to say that when $\varepsilon(r)$ tends to
$0$, the construction above comes with slightly better estimates. We
shall assume that 
\begin{equation}\label{10.2}  
\sum_{k \geq 0} \varepsilon(2^{-k}) < + \infty,
\end{equation}
and then show that the restriction of our mapping $f$ to each of the
three or six faces of $Z(0,2)  \cap B(0,17/10)$, is of class $C^1$.

Set $\varepsilon_k = \varepsilon(2^{-k})$. For technical reasons, it 
will be more convenient to assume that $\varepsilon_k$ never drops 
too brutally, i.e., that 
\begin{equation}\label{10.3}  
\varepsilon_{k} \leq C_0 \varepsilon_{k+1}
\ \hbox{ for $k \geq 0$}
\end{equation}
and some constant $C_0 > 0$. This is not really an additional 
assumption, because we can always replace $\{\varepsilon_k\}$ above 
with a larger sequence that satisfies (\ref{10.3}) and such that 
$\sum_{k \geq 0} \varepsilon_k$ still converges.

The point of assuming (\ref{10.3}) is that when we apply the
construction and arguments above, we can replace $\varepsilon$
with $C \varepsilon_n$ in our various intermediate statement 
(such as Lemmas~\ref{lem4.13} and the Lipschitz descriptions in 
(\ref{eqn5.12})-(\ref{eqn5.15})). If we did not assume (\ref{10.3}),
our proofs by induction would force us to carry additional 
terms that control the larger scales.

In particular, the proof above shows that
\begin{equation}  \label{10.4} 
(1-C\varepsilon_n) \dist(y,z)
\leq \dist(\widetilde g_n(y),\widetilde g_n(z))
\leq (1+C\varepsilon_n) \dist(y,z)
\end{equation}
for $y$, $z\in B(0,185/100)$ such that $|y-z| \leq 2^{-n-80}$,
instead of (\ref{9.6}).

We need estimates on the first and second derivatives of $f$ near
the faces. We start with the mappings $g_n^\ast$ and $f_n^\ast$
from Section \ref{first}. We see them as mappings defined on $\R^3$,
because this will make it easier to differentiate, even though we shall
later restrict to a small neighborhood of the spine of $Z(0,2)$. When
we differentiate (\ref{eqn4.3}), we get that
\begin{equation}  \label{10.5} 
Dg_n^\ast = \sum_{i\in I(n)} D\theta_i \psi_i^\ast
+ \sum_{i\in I(n)} \theta_i D\psi_i^\ast
= \sum_{i\in I(n)} D\theta_i \, [\psi_i^\ast - \psi_j^\ast]
+ \sum_{i\in I(n)} \theta_i \, D\psi_i^\ast, 
\end{equation}
where $j \in I(n)$ is any fixed index and we used the fact that
$\sum_{i\in I(n)} \theta_i = 1$. Here we just need to differentiate
near $x\in \Gamma_n = f_n^\ast(\Gamma)$, so we can restrict to 
$i\in I_2(n)$ because the other $\theta_i$ vanish near $x$, and then
$\psi_i^\ast$ is the orthogonal projection onto the line $L_i$. We
select $j \in I_2(n)$ such that $\theta_j(x) \neq 0$, and then
\begin{equation}  \label{10.6} 
|\psi_i^\ast - \psi_j^\ast| \leq C \varepsilon_n 2^{-n}
\ \hbox{ and } \ 
|D\psi_i^\ast - D\psi_j^\ast| \leq C \varepsilon_n
\ \hbox{ near $x$,}
\end{equation}
by (\ref{eqn4.5}) and (\ref{eqn4.6}). Thus (\ref{10.5}) yields
\begin{equation}  \label{10.7} 
|Dg_n^\ast - D \psi_j^\ast| \leq C \varepsilon_n
\ \hbox{ near $x$,}
\end{equation}
where we used the fact that 
$D\psi_j^\ast = \sum_{i\in I(n)} \theta_i D\psi_j^\ast$,
and (\ref{eqn3.12}) to control the $|D\theta_i|$. In particular,
\begin{equation}  \label{10.8} 
|Dg_n^\ast| \leq 1 + C \varepsilon_n
\ \hbox{ near $x$,}
\end{equation}
When we differentiate a second time in (\ref{10.5}), we get that
\begin{equation}  \label{10.9} 
D^2g_n^\ast 
= \sum_{i\in I(n)} D^2\theta_i \, [\psi_i^\ast - \psi_j^\ast] 
+ 2\sum_{i\in I(n)} D\theta_i \, [D\psi_i^\ast - D\psi_j^\ast],
\end{equation}
because $\sum_{i\in I(n)} D\theta_i = 0$ and the $D^2\psi_i^\ast$ vanish 
(each $\psi_i^\ast$ is a linear projection). Thus
\begin{equation}  \label{10.10} 
|D^2g_n^\ast| \leq C \varepsilon_n 2^{n}
\ \hbox{ near $x$,}
\end{equation}
by (\ref{10.6}) and (\ref{eqn3.12}). Next $f_n^\ast$ is the
composition of the $g_k^\ast$, $k<n$ (see (\ref{eqn4.2})), so
\begin{equation}  \label{10.11} 
Df_n^\ast = Dg_{n-1}^\ast \circ Dg_{n-2}^\ast\circ \cdots \circ Dg_0^\ast,
\end{equation} 
where we do not write the arguments to save space. When we 
differentiate once more, we get a sum of terms with one second
derivative of some $g_k^\ast$, composed with a collection of
first-order differentials; thus
\begin{equation}  \label{10.12} 
|D^2f_n^\ast| \leq C \varepsilon_n 2^{n} 
\prod_{0 \leq k < n} (1 + C \varepsilon_n)
\leq C \varepsilon_n 2^{n}, 
\end{equation}
by (\ref{10.8}), (\ref{10.10}), and (\ref{10.2}).

Let us now consider the restriction of $f_n^\ast$ to
$\Gamma$, the intersection of the spine of $Z(0,2)$ with
$\overline B(0,197/100)$. When $Z(0,2)$ is of type $3$ we
can still use the discussion above, we just need restrict our
attention to $\Gamma \setminus \{ 0 \}$. The derivative of this
restriction is (locally, and modulo identification of $\Gamma$
with an interval in $\R$) a vector-valued function  ${f_n^\ast}'$, 
and $|{f_n^\ast}'(x)| \geq \prod_{0 \leq k < n} (1 - C \varepsilon_n)
\geq C^{-1}$ by multiple applications of (\ref{eqn4.17}),
or equivalently of (\ref{10.7}) and the fact that $\Gamma_n$ is a small
Lipschitz graph over $L_i$. 

We need these estimates to control the parameterization
$x \to (x,\zeta(x))$ of the graph $G_{n,i}$ that we use in 
(\ref{eqn5.7})-(\ref{eqn5.9}). Observe that 
$|\zeta''(x)| \leq C \varepsilon_n 2^{n}$ because 
$|{f_n^\ast}'(x)| \geq C^{-1}$ and 
$|{f_n^\ast}''(x)| \leq |D^2f_n^\ast| \leq C \varepsilon_n 2^{n}$
(by (\ref{10.12})). Then
\begin{equation}  \label{10.13} 
|D\eta_i - I| \leq C \varepsilon_n
\ \hbox{ and } \ 
|D^2\eta_i|\leq C \varepsilon_n 2^{n}, 
\end{equation}
by (\ref{eqn5.7}) and (\ref{eqn5.8}).

Next fix a face $F_n^l$ and a point $z$ in the interior of 
$F_n^l$. Assume that $z\in B(0,195/100) \i B(0,\rho_n)$,
so that we can use the estimates below (\ref{eqn5.15}).
Let $i \in I(n)$ be such that $z\in 3B_i$ 
(i.e., with the notation of Section \ref{param}, $i\in I_z$). 
If $i\in I_1(n)$, (\ref{eqn5.5}) says that $\psi_i=\pi_i$, 
the orthogonal projection onto a plane $P_i$. If $i\in I_2(n)$,
$\psi_i=h_i \circ \eta_i$ by (\ref{eqn5.9}), and (\ref{eqn5.19}) says 
that 
\begin{equation}  \label{10.14} 
\psi_i= \pi_i^l \circ \eta_i
\ \hbox{ near $z$,}
\end{equation}
because $z$ lies in the interior of $F_n^l$. 
Here $\pi_i^l$ is the orthogonal projection onto
the plane that contains the appropriate face of $Y_i$.
If $Z(0,2)$ is of type $3$ and $i\in I_3(n) = \{ i_0 \}$,
the same argument shows that $\psi_i= \pi_i^l \circ \eta_i$ near $z$, 
where $\pi_i^l$ is the orthogonal projection onto the plane that contains 
some face of $Z_n$. The case when $i\in I_0(n)$ does not happen, 
as noted below (\ref{eqn5.15}), by (\ref{eqn5.12}) and (\ref{eqn3.23}). 
Set $\pi_i^l = \pi_i$ when $i\in I_1(n)$, to uniformize the notation.

If $j \in I(n)$ is some other index such that $z\in 3B_j$, then 
$|D\pi_j^l-D\pi_i^l| \leq C \varepsilon_n$, by the proof of 
(\ref{eqn5.18}) (recall, we used the fact that near $z$, $F_n^l$ is 
at the same time a small Lipschitz graph over $P_i^l$ and $P_i^l$).
We easily deduce from this, (\ref{10.13}), (\ref{10.14}) and its analogue 
for $i\in I_3(n)$ (both applied to $j$), that
\begin{equation}  \label{10.15} 
|D\psi_j - D\pi_i^l| \leq C \varepsilon_n
\hbox{ and } |D^2\psi_j| \leq C \varepsilon_n 2^n
\ \hbox{ near $z$,}
\end{equation}
We also know from (\ref{eqn5.17}) and (\ref{eqn5.18}) 
that $|\psi_j(z) - \pi_i^l(z)| \leq C \varepsilon_n 2^{-n}$.
The same computation as for (\ref{10.5}) and (\ref{10.7})
yields
\begin{equation}  \label{10.16} 
Dg_n = \sum_{j\in I(n)} D\theta_j \psi_j
+ \sum_{j\in I(n)} \theta_j D\psi_j
= \sum_{j\in I(n)} D\theta_j \, [\psi_j - \pi_i^l] 
+ \sum_{i\in I(n)} \theta_j \, D\psi_j,
\end{equation}
by (\ref{eqn5.1}), and then 
\begin{equation}  \label{10.17} 
|Dg_n - D \pi_i^l| \leq C \varepsilon_n 
+ \Big|\sum_{i\in I(n)} \theta_j \, [D\psi_j -D\pi_i^l] \Big|
\leq C \varepsilon_n
\end{equation}
near $z$. Next 
\begin{equation}  \label{10.18} 
D^2g_n = \sum_{j\in I(n)} D^2\theta_j \, [\psi_j - \pi_i^l] 
+ 2\sum_{j\in I(n)} D\theta_j \, [D\psi_j - D\pi_j^l] 
+ \sum_{j\in I(n)} \theta_j D^2\psi_j \, ,
\end{equation}
so $|D^2g_n| \leq  C \varepsilon_n 2^n$ near $z$, by (\ref{10.15})
and the line below it. The same proof as for (\ref{10.12}) now shows 
that
\begin{equation}\label{10.19}  
|D^2 f_n| \leq C \varepsilon_n 2^{n}
\ \hbox{ on the interior of } F^l \cap B(0,194/100).
\end{equation}
(we only consider $x\in B(0,194/100)$ to make sure that all the
$f_n(x)$ lie in $B(0,195/100)$.)

Let $x$ in the interior of $F^l \cap B(0,194/100)$ be given,
and let $u$ be a unit tangent vector to $F^l$ at $x$.
Set $z_n = f_n(x)$ and $u_n = Df_n(x)(u)$. Then $u_n$ is a tangent 
vector to $F_n^l$ at $z_n$, and 
\begin{eqnarray}\label{10.20}  
|u_{n+1}-u_n| &=& |Dg_n(z_n)(u_n)-u_n|
\leq |Dg_n(z_n)(u_n)- \pi_i^l(u_n)| + |\pi_i^l(u_n)- u_n| 
\nonumber \\
&\leq& C \varepsilon_n |u_n| + |\pi_i^l(u_n)- u_n|
\leq C \varepsilon_n |u_n| 
\end{eqnarray}
where $\pi_i^l$ is associated to $z_n$ as above, by (\ref{10.17}),
and because $F_n^l$ is a small Lipschitz graph over $P_i^l$ near $z_n$.
Multiple applications of this yield $C^{-1} \leq |u_n| \leq C$, by 
(\ref{10.2}), and $|u_{n+1}-u_n| \leq C \varepsilon_n$. 

We just proved that the derivatives $Df_n(x)$ converge, uniformly
in $x$ in the interior of $F^l \cap B(0,194/100)$, to some limit $G$.
Since we already know that the $f_n$ converge to $f$, we get that
$f$ is continuously differentiable on the interior of 
$F^l \cap B(0,194/100)$, with $Df=G$. In addition,
$|Df(u)-Df_n(u)| \leq \sum_{k \geq n} |u_{k+1}-u_k| 
\leq C\widehat\varepsilon_n$, where we set
$\widehat\varepsilon_n = \sum_{k \geq n} \varepsilon_k$.
Then, if $x$ and $y$ both lie in the interior of $F^l \cap B(0,194/100)$,
\begin{eqnarray}\label{10.21}  
\,\hskip0.9cm
|Df(x)-Df(y)| &\leq& |Df(x)-Df_n(x)|+|Df_n(x)-Df_n(y)|+|Df_n(y)-Df(y)|
\nonumber \\
&\leq& C \widehat\varepsilon_n + |Df_n(x)-Df_n(y)|
\leq C \widehat\varepsilon_n + C \varepsilon_n 2^n |x-y|,
\end{eqnarray}
by (\ref{10.19}). We choose $n$ such that $2^{-n} \leq |x-y| < 2^{-n+1}$,
and we get that 
\begin{equation}  \label{10.22} 
|Df(x)-Df(y)| \leq C \widehat\varepsilon_n 
\leq C \sum_{k \geq 0 \, ; \, |x-y| \leq 2^{-k}} \varepsilon(2^{-k}).
\end{equation}

Of course we can let $x$ or $y$ tend to the boundary of $F^l$.
So we proved that the restriction of $f$ to each face 
$F^l \cap B(0,194/100)$ is of class $C^1$, with a modulus of 
continuity for $Df$ that depends on the decay of the function 
$\varepsilon(r)$.
For instance, if $\varepsilon(r) \leq C r^{\alpha}$ for some
$\alpha \in (0,1)$, we get corresponding $C^{1,\alpha}$ estimates
for $f$ on $F^l$.

Now the reader should not pay too much attention to this more
precise result; we proved it because it is easy (once we have the
other estimates), but it is also easy to prove directly that on
each face, there is a tangent plane at $x$ that depends nicely on
$x$. Indeed, we can use the good approximation of $E$ near $x$
by planes $P(x,2^{-k})$ (or sets of type $2$ or $3$ for the scales 
that are larger than $\dist(x,E_2)$) to show that the $P(x,2^{-k})$ vary 
slowly with $k$, then tend to a plane $P(x)$. Then $P(x)$ is tangent to 
$E$ at $x$, and we can use the fact that $P(x)$ is close to $P(x,2^{-k})$
to estimate the distance from $P(x)$ to $P(y)$ (compare $P(x,2^{-k})$
to $P(y,2^{-k})$, with $2^{-k}$ a little larger than $|x-y|$).

Our construction does not give a mapping $f$ that is $C^1$
on $Z(0,2)$ near its spine. But on the other hand, each face of
$E$ is $C^1$, and the faces make the right $120^{\circ}$ angles with 
each other along $E_2$, so we could fairly easily re-parameterize 
$E$ locally in a $C^1$ way. Once again, if we know that $E$ is
fairly nice and each piece is smooth, getting good parameterizations 
is rather easy.

\bigskip\bigskip

\noindent G.\ David, Universit\'{e} de Paris-Sud, 
Batiment 425, 91405 Orsay, France.\\
\noindent E-mail: guy.david@math.u-psud.fr

\vskip 10mm

\noindent T. \ De Pauw, Universit\'{e} Catholique de Louvain,
Chemin du cyclotron 2, B-1348 Louvain-la-Neuve, Belgique.\\
\noindent E-mail: depauw@math.ucl.ac.be

\vskip 10mm

\noindent T.\ Toro, University of Washington, Department of Mathematics, Box
354350, Seattle, WA 98195--4350.\\
\noindent E-mail: toro@math.washington.edu

\end{document}